\documentclass[envcountsame,envcountsect]{svmult}

\usepackage{makeidx}         % allows index generation
\usepackage{graphicx}        % standard LaTeX graphics tool
\usepackage{multicol}        % used for the two-column index
\usepackage[bottom]{footmisc}% places footnotes at page bottom
\usepackage{amsmath,amssymb,latexsym} %AMSFonts and AMS and Latex symbols
\usepackage{eucal} %%inserted by E, produces scr D anstatt cal D

\usepackage[all]{xy}
\makeindex \smartqed
\spnewtheorem{construction}[theorem]{Construction}{\bfseries}{\rmfamily}

\newcommand\sm{\vskip\smallskipamount} %from Retro
\newcommand\ms{\vskip\medskipamount} %from Retro
\newcommand\bs{\vskip\bigskipamount} %from Retro
\newcommand\vfe{}

\newtheorem{eia}{}[section]
\def\subsec[#1]{\begin{eia}\label{#1}\rm }
\def\es{  \end{eia}}
\def\endsubsec{ \end{eia}}
\def\eq{\eqno}

\def\iiitem{\item}
\def\iitem{\indent}

\def\amen#1#2{
     \hsize=#1\hsize
     \vbox{
        \hbox{
             \vbox{ %\vskip \parindent
                      \relax
                  \noindent \it  #2 }%\smallskip} %
            }} }

\def\Corollary{{\bf Corollary.} }
\def\Theorem{{\bf Theorem.} }
\def\Lemma{{\bf Lemma.} }
\def\Proposition{{\bf Proposition.} }

\newcommand\an{^{\rm an}}
\newcommand\inpr{(\cdot|\cdot)}

\newcommand\Iff{\Longleftrightarrow}
\newcommand\lan{\langle} \newcommand\ran{\rangle}

\newcommand\ot{\otimes}
\newcommand\pa{\partial}

\newcommand\bull{\bullet}
\newcommand\tsum{\textstyle \sum}
\newcommand\ti{^\times}
\newcommand\wi{\widehat}
\newcommand\ch{\sp{\scriptscriptstyle\vee}}
\newcommand\tl{\tilde}
\def\re{^{\rm re}}
\def\im{^{\rm im}}

\def\int{^{\rm int}}
\def\ind{_{\rm ind}}
\def\div{_{\rm div}}
\def\sh{_{\rm sh}}
\def\lg{_{\rm lg}}
\def\ts{\textstyle}
\def\bfm{{\bf m}}
\def\INV#1{{\rm (IA#1)}}
\def\g*{^{{\rm gr}*}}
\def\EA#1{{\rm (EA#1)}}

 \DeclareMathOperator\ad{ad}
 
 \DeclareMathOperator\Aut{Aut}
 \DeclareMathOperator\rmC{C}
 \DeclareMathOperator\grrmC{grC}
 \DeclareMathOperator\cent{Cent}
 \DeclareMathOperator\CDer{CDer}
 \DeclareMathOperator\Der{Der}
 \DeclareMathOperator\ev{ev}
 \DeclareMathOperator\End{End}
 \DeclareMathOperator\GL{GL}
 \DeclareMathOperator\GIF{GIF}
 \DeclareMathOperator\grcent{grCent}
\DeclareMathOperator\grEnd{grEnd}
 \DeclareMathOperator\grDer{grDer}
 \DeclareMathOperator\grCDer{grCDer}
 \DeclareMathOperator\grSCDer{grSCDer}
 \DeclareMathOperator\Hom{Hom}
 \DeclareMathOperator\HC{HC}
 \DeclareMathOperator\Span{span}
 \DeclareMathOperator\supp{supp}
 \DeclareMathOperator\rmZ{Z}
 \DeclareMathOperator\Id{Id}
 \DeclareMathOperator\IDer{IDer}
 \DeclareMathOperator\Ker{Ker}
 \DeclareMathOperator\Mat{Mat}
 \DeclareMathOperator\Rad{Rad}
 \DeclareMathOperator\rank{rank}
 \DeclareMathOperator\SCDer{SCDer}
 \DeclareMathOperator\tr{tr}

\newcommand{\al}{\alpha}
\newcommand{\be}{\beta}
\newcommand{\ga}{\gamma} \def\Ga{\Gamma}
\newcommand\de{\delta}
\def\eps{\epsilon}
\newcommand\io{\iota}
\newcommand\ka{\kappa}
\newcommand{\la}{\lambda}
\newcommand{\La}{\Lambda}
\newcommand\ph{\varphi} %\def\phi{\varphi}

\newcommand\si{\sigma}
 \newcommand\boldsi{{\boldsymbol{\sigma}}}
\newcommand\ta{\tau}
 \newcommand\vth{\theta} %\newcommand\vth{\theta}

\newcommand\ze{\zeta}

\newcommand\Om{\Omega}

\newcommand\CC{{\mathbb C}}
\newcommand\KK{{\mathbb K}}
\newcommand\NN{{\mathbb N}}
\def\RR{{\mathbb R}}
\def\SS{{\mathbb S}}
\newcommand\QQ{{\mathbb Q}}
\newcommand\ZZ{{\mathbb Z}}

\newcommand\frf{{\mathfrak f}}
\newcommand\frg{{\mathfrak g}}

\newcommand\frh{{\mathfrak h}}
\newcommand\frk{{\mathfrak k}}

\newcommand\frL{{\mathfrak L}}
\newcommand\frslI{{\mathfrak s}{\mathfrak l}_I}
\newcommand\frpsl{{\mathfrak p}{\mathfrak s}{\mathfrak l}}
\newcommand\frsln{{\mathfrak s}{\mathfrak l}_n}

\newcommand\fru{{\mathfrak u}}

\newcommand\uce{{\mathfrak u}{\mathfrak c}{\mathfrak e}}
\newcommand\frgl{{\mathfrak g}{\mathfrak l}}
\newcommand\frst{{\mathfrak s}{\mathfrak t}}

\newcommand\scB{{\mathcal  B}}
\newcommand\scC{{\mathcal  C}}
\newcommand\scD{{\mathcal  D}}
\newcommand\scE{{\mathcal E}}
\newcommand\scM{{\mathcal  M}}
\newcommand\scQ{{\mathcal  Q}}
\newcommand\scq{{\mathcal Q}}

\newcommand\rma{\dot{\rm A}}
\newcommand\rmfa{{\rm A}}
\newcommand\rmb{{\rm B}}
\newcommand\rmbc{{\rm BC}}
\newcommand\rmc{{\rm C}}
\newcommand\rmd{{\rm D}}
\newcommand\rme{{\rm E}}
\newcommand\rmf{{\rm F}_4}
\newcommand\rmg{{\rm G}_2}

\begin{document}

\title*{Extended affine Lie algebras and other generalizations
of affine Lie algebras -- a survey}
\titlerunning{EALAs and generalizations}
\author{Erhard Neher\inst{1}\thanks{Support of the Natural
Sciences and Engineering Research Council of Canada is gratefully
acknowledged.}} \institute{Department of Mathematics and Statistics,
University of Ottawa, Ottawa, Ontario K1N 6N5, Canada, email:
\texttt{neher@uottawa.ca} }

%\date{\today}

\maketitle

\abstract This is a survey on extended affine Lie algebras and
related types of Lie algebras, which generalize affine Lie algebras.
{\keywords{affine Lie algebras, extended affine Lie algebra, Lie
torus, AMS subject classification: Primary 17B65, 17B67. Secondary
16W50, 17B70.}}

%Subject classifications
%16W50 Graded rings and modules
%17B60 Lie (super)algebras associated with other structures (associative, Jordan, etc.)
%17B65 Infinite-dimensional Lie (super)algebras
%17B67 Kac-Moody (super)algebras
%17B70 Graded Lie (super)algebras

\section{Introduction} \label{perintro}

{\it Motivation.} The theory of affine (Kac-Moody) Lie algebras has
been a tremendous success story. Not only has one been able to
generalize essentially all of the well-developed theory of
finite-dimensional simple Lie algebras and their associated groups
to the setting of affine Lie algebras, but these algebras have found
many striking applications in other parts of mathematics. It is
natural to ask -- why? What is so special about affine Lie algebras?
What really makes them ``tick''? One way to understand this, is to
generalize affine Lie algebras and see where things go wrong or
don't. After all, it is conceivable that there is a whole theory of
Lie algebras, ready to be discovered, for which affine Lie algebras
are just one example.
\smallskip

Of course, one immediately thinks about Kac-Moody Lie algebras, for
which finite-dimensional simple and affine Lie algebras are the
basic examples. But is it the right generalization? This is not such
a sacrilegious question as it seems. After all, one the founders and
main contributor to the theory does not seem to be convinced himself
(\cite{kac:locality}). So where to go? There are several directions
that have been pursued: toroidal algebras, Borcherds algebras,
GIM-algebras, root-graded algebras, to name just a few. \ms

{\it The first generalization.} The generalizations that are of
interest for this paper arose in the work of Saito and Slodowy on
elliptic singularities and in the paper by the physicists
H{{\o}}egh-Krohn and Torr{\'e}sani \cite{HKT} on Lie algebras of
interest to quantum gauge field theory. The latter paper, written by
physicists, did however not stand up to mathematical scrutiny.
Fortunately, a mathematical sound foundation of the theory was later
provided by Allison, Azam, Berman, Gao and Pianzola in their AMS
memoirs \cite{aabgp}. They also describe the type of root systems
appearing in these algebras and give many examples. Last not least
they coined the name extended affine Lie algebras for this new class
of algebras. \sm

Before we give the definition, it is convenient to introduce the
basic concept of a toral pair over $\KK$. Throughout this paper
$\KK$ denotes a field of characteristic $0$. A {\it toral pair\/}
$(E,T)$ over $\KK$ consists of a Lie algebra $E$ over $\KK$ and a
subspace $T\subset E$ acting $\ad$-diagonalizably on $E$, i.e.,
putting
 \begin{align*}
    E_\al &= \{ e\in E : [t,e] = \al(t) e
                \hbox{ for all }
       t\in T\} \quad \hbox{for $\al \in T^*$ and} \cr
 R &= \{ \al \in T^*: E_\al \ne 0\}, \end{align*}
we have $E = \bigoplus_{\al \in R} E_\al$. It turns out that $T$ is
abelian, so $T \subset E_0$ and hence $0 \in R$ (unless $T=\{0\}=
E$). In general $T$ is a proper subspace of $E_0$. If, however,
$T=E_0$ we will call $T$ a {\it splitting Cartan subalgebra\/} and,
following common usage, will use the letter $H$ instead of $T$. The
following definition is essentially the one given in \cite{aabgp}:
\sm

An {\it extended affine Lie algebra\/}, or EALA for short, is a
toral pair $(E,H)$ over $\KK=\CC$ satisfying the conditions
(E1)--(E6) below:
\begin{description}
\item{(E1)} $E$ has an invariant nondegenerate symmetric bilinear form
$\inpr$.\sm

\item{(E2)} $H$ is a non-trivial, finite-dimensional splitting Cartan
subalgebra. \end{description} \noindent  Since $H=E_0$, it is
immediate that the restriction of $\inpr$ to $H$ is nondegenerate.
We can therefore represent any $\al \in H^*$ by a unique vector
$t_\al \in H$ via $\al(h) = (t_\al | h)$. Hence
\begin{equation}
\label{perintro0}
  R= R^0 \cup R\an %\eq{perintro0}
\end{equation}
 for $R^0 = \{ \al \in R : (t_\al | t_\al) = 0 \}$ ({\it null
roots\/}) and $ R\an = \{ \al \in R : (t_\al | t_\al) \ne 0 \}$
({\it anisotropic roots\/}). We can now state the remaining axioms
(E3)--(E6):
\begin{description}
 \item{(E3)} For any $\al \in R\an$ and $x \in E_\al$, the adjoint
map $\ad x$ is locally nilpotent on $E$. \sm

\item{(E4)} $R\an$ is connected in the sense that any partition
$R\an = R_1 \cup R_2$ with $(R_1 | R_2) = 0$ either has $R_1=
\emptyset$ or $R_2 = \emptyset$. \sm

\item{(E5)} $R^0 \subset R\an + R\an$. \sm

\item{(E6)} $R$ is a discrete subset of $H^*$, equipped with the
natural topology. \end{description}

%In general there will be more roots in $R^0$ than just $0\in R^0$.
%But the essence of axiom (E5) below is that $R\an$ controls $R^0$.

Let $(E,H)$ be an EALA. By (E6), the $\ZZ$-span of $R^0$ in $H^*$ is
free of finite rank, called the {\it nullity\/} of $(E,H)$ (it turns
out that this is independent of $\inpr$). The {\it core of\/
$(E,H)$\/} is the ideal $E_c$ of $E$ generated by all subspaces
$E_\al$, $\al \in R\an$. We then have a Lie algebra homomorphism
$E\to \Der E_c$, given by $e \mapsto \ad e\mid E_c$, and one says
that $(E,H)$ is {\it tame\/} if the kernel of this map lies in
$E_c$. The idea of tameness is that one can recover $E$ from $E_c$
and its derivations. Although there are examples of non-tame EALAs,
it is natural to assume tameness for the purpose of classification.
\sm

We can now convince the reader that extended affine Lie algebras
indeed do what they are supposed to do: A tame EALA of nullity $0$
is the same as a finite-dimensional simple Lie algebra (exercise),
and a tame EALA of nullity $1$ is the same as an affine Lie algebra,
proven by Allison-Berman-Gao-Pianzola in \cite{abgp}. Tame EALAs of
nullity $2$ are closely related to the Lie algebras studied by Saito
and Slodowy, and toriodal Lie algebras provide examples of tame
EALAs of arbitrarily high nullity. What makes EALAs an interesting
class of Lie algebras is that there are other, non-toroidal examples
requiring noncommutative coordinates. \ms

{\it The revisions.} The system of axioms (E1)--(E6) has been very
successful. It describes an important and accessible class of Lie
algebras, about which much is known today. But there are drawbacks.
Foremost among them is the restriction to the ground field $\CC$,
which is necessary for the axiom (E6). Looking at the basic examples
of EALAs given in \cite{aabgp}, one sees that they make sense over
arbitrary fields $\KK$. Hence, in \cite{n:eala} the author proposed
a new definition of EALAs, valid for Lie algebras over $\KK$. It
used the same set of axioms, but replaced (E6) with its essence:
\begin{description}
\item{(E6$'$)} the $\ZZ$-span of $R$ is free of finite rank.
\end{description}
\noindent (To be more precise, only tame EALAs were considered in
\cite{n:eala}.) It turned out that one can develop a satisfactory
theory for EALAs over $\KK$, partially described in \cite{n:eala}
and \cite{n:tori}. \sm

Some drawbacks however still remained: The theory only covers the
``split case'' because of the condition in (E2) that $H$ be a
splitting Cartan subalgebra, and it does not include
finite-dimensional semisimple Lie algebras because of (E4), nor
infinite rank affine Lie algebras because of the requirement in (E2)
that $H$ be finite-dimensional. Some of these issues were addressed
in papers by Azam \cite{az:grla}, Azam-Khalili-Yousofzadeh
\cite{AKY} and Morita-Yoshii \cite{morita-yoshii} (see
\ref{LEALAdef}--\ref{stinvd}). \sm

A radically new look at the system of axioms seems to be called for:
Delete the less important axioms (E4)--(E6) and modify (E2) so that
all the aforementioned examples can be included. The purpose of this
survey and research announcement is to propose such an approach by
defining a new class of Lie algebras and to explain its main
features. \ms

{\it The new definition.} An {\it invariant affine reflection
algebra\/}, an IARA for short, is a toral pair $(E,T)$ over $\KK$
satisfying the axioms (IA1)--(IA3) below (recall $E= \bigoplus_{\al
\in R} E_\al$):

\begin{description}
\item{(IA1)} $E$ has an invariant nondegenerate symmetric bilinear
form $\inpr$ whose restriction to $T$ is nondegenerate and which
represents every root $\al\in R$: There exists $t_\al \in T$ such
that $\al(t)= (t_\al | t)$ for all $t\in T$. \ms

\item{(IA2)} For every $\al \in R$, $\al \ne 0$,  there are $e_{\pm
\al} \in E_{\pm \al}$ such that $0\ne [e_\al,\, e_{-\al}] \in T$.\sm

\item{(IA3) = (E3)} For every $\al \in R$ with $(t_\al | t_\al ) \ne 0$
and for all $x_\al \in E_\al$ the adjoint map $\ad x_\al$ is locally
nilpotent on $E$. \end{description}

\noindent The following are examples of invariant affine reflection
algebras: finite-dimensional reductive Lie algebras, infinite-rank
affine Lie algebras, Lie algebras satisfying (E1)--(E3), in
particular EALAs, and all the generalizations of EALAs that had been
considered before. \sm

To describe the structure of IARAs, we will use the inclusions of
categories
\begin{equation}
 \label{perintro1} \hbox{affine Lie algebras} \qquad
\subset \qquad \hbox{EALAs}
 \qquad \subset \qquad
 \hbox{IARAs}
\end{equation}
 as a leitmotif. We will first discuss the combinatorics, viz.
the roots appearing in IARAs, and then describe a procedure to
determine their algebra structure. \ms

{\it The roots.} The roots for the algebras in (\ref{perintro1})
form, respectively, an
$$
\hbox{affine root sytem} \qquad \subset \qquad
 %\hbox{extended affine root systems} \cr &
 \hbox{EARS} \qquad \subset \qquad
 \hbox{ARS}
$$
where the affine root systems are the roots appearing in affine Lie
algebras, EARS is the acronym for an extended affine root system,
the roots of EALAs, and ARS stands for affine reflection system
(\ref{affdef}), of which the roots of an IARA are an example. To
explain the general structure of an ARS, recall (\ref{afrs}) that to
any affine root system $R$ one can associate a finite irreducible,
but not necessarily reduced root system $S$, a torsion-free abelian
group $\La$ and a family $(\La_\al  : \al \in S)$ of subsets
$\La_\al \subset \La$ such that
\begin{equation}
  R = \ts \bigcup_{\xi \in S} \, (\xi \oplus \La_\xi)
   \quad \subset \quad
   \Span_\KK(R) = \Span_\KK (S) \oplus \Span_\KK(\La). \label{perintro2}
\end{equation}
 Indeed, $\La = \ZZ \de$ where $\de$ is the
basic null root. Now (\ref{perintro2}) also holds for EARS if one
takes $\La$ to be free of finite rank. Since the toral subalgebra of
an IARA may be infinite-dimensional, it is not surprising that root
systems of infinite rank will occur. Indeed, if one allows $S$ to be
a locally finite root system (\ref{affdef}) and $\La$ to be any
torsion-free abelian group, then (\ref{perintro2}) also holds for
ARSs. Moreover, the basic relations between the subsets $\La_\xi$,
$\xi \in S$, developed in \cite[Ch.II]{aabgp} for EARS, actually
also hold in the more general setting of affine reflection systems.
Hence, their theory is about as easy or complicated as the theory of
EARS.  \ms

{\it The core and centreless core.} Recall that the core of an EALA
$(E,H)$ is the ideal generated by all $E_\al$, $\al \in R\an$. It is
a perfect Lie algebra. Hence $E_{cc} = E_c / \rmZ(E_c)$ has centre
$\{0\}$ and is justifiably called the {\it centreless core}. By
(IA1) the decomposition $R=R^0 \cup R\an$ is also available for
IARAs, and we can therefore define $E_c$ and $E_{cc}$ as in the
EALA-case. Hence, from any IARA $E$ we can go down to its centreless
core $E_{cc}$:
$$
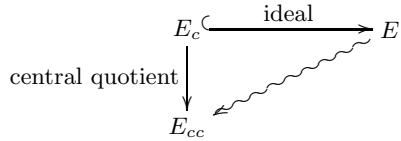
\xymatrix@C=60pt{ E_c \ar@{^{(}->}[r]^{\hbox{ideal}}
         \ar[d]_{\hbox{central quotient}}
 &  E  \ar@{~>}[dl]\cr
 E_{cc}}
$$
This diagram describes the fundamental approach to the study of the
Lie algebras $(E,T)$ in (\ref{perintro1}):
\begin{description}
 \item{(A)} Describe the class $\scC$ of Lie algebras that arise as
centreless cores $E_{cc}$ in a way which is independent of $(E,T)$,
\sm

\item{(B)} for every Lie algebra $L$ in $\scC$ provide a
construction which yields all the Lie algebras $E$ with centreless
core $L$, and \sm

\item{(C)} classify the Lie algebras in $\scC$. \end{description}

\noindent That this approach has good chances of success, can
already be witnessed by looking at affine Lie algebras. For them,
the core and the centreless core are the basic ingredients of the
structure theory: $E_c$ is the derived algebra of the affine Lie
algebra $E$ and $E_{cc}$ is the (twisted or untwisted) loop algebra
that one uses to construct $E$ from $E_{cc}$. The centreless cores
of EALAs were described by Allison-Gao \cite{AG2} and then put into
an axiomatic framework as centreless Lie $\ZZ^n$-tori by Yoshii
\cite{y:lie} and the author \cite{n:tori}. One can show that the
centreless core of an IARA is a certain type of root-graded Lie
algebra, a so-called invariant predivision-$(S,\La)$-graded Lie
algebra for $(S,\La)$ as above: $S$ is a locally finite root system
and $\La$ is a torsion-free abelian group. Conversely, every such
Lie algebra arises as the centreless core of some IARA (\ref{rogr},
\ref{invth}, \ref{invcop}). Thus, for the Lie algebras in
(\ref{perintro1}) the classes $\scC$ are
 \begin{eqnarray*}
  \hbox{loop algebras} &\quad &\subset \quad \hbox{centreless Lie
   $\ZZ^n$-tori}\quad \subset  \\
 &\quad& \hbox{invariant predivision-root-graded Lie algebras, $\La$
torsion-free.} \end{eqnarray*}
 This solves problem (A). \sm

As for (B), we describe in \ref{invcon} a construction, which starts
with any invariant predivi\-sion-root-graded Lie algebra $L$ with a
torsion-free $\La$ and produces an invariant affine reflection
algebra $(E,T)$ with centreless core $L$. It maintains the main
features of the construction of an affine Lie algebra from a loop
algebra, i.e. take a central extension and then add on some
derivations:
$$
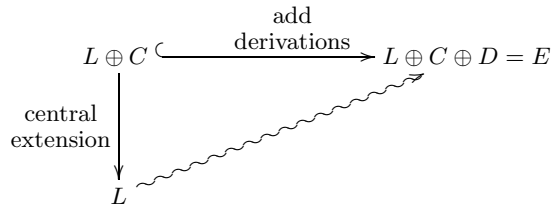
\xymatrix@C=80pt@R=40pt{ L \oplus C \;
 \ar@{^{(}->}[r]^{\txt{add \\ derivations}}
         \ar[d]_{\txt{central \\ extension}}
 &  L \oplus C \oplus D = E\cr
 L \ar@{~>}[ur]}
$$
If $L$ is a Lie torus this construction yields all EALAs with
centreless core $L$, and hence settles problem (B) for them. \sm

Finally, the efforts of many people have now lead to a complete
solution of problem (C) for Lie tori, see the survey by
Allison-Faulkner in \cite[\S 12]{AF:isotopy}. The classification of
predivision-root-graded Lie algebras is not known in general, but in
several important cases, see the review in \ref{perrsgrno}. What
certainly is fascinating about this class and root-graded Lie
algebras in general, is that there are many examples requiring
noncommutative coordinates. For example, for $S=\rmfa_l$, $l\ge 3$,
one has to consider central quotients of ${\mathfrak s}{\mathfrak
l}_{l+1}(A)$ where $A$ is an associative algebra, a crossed product
algebra to be more precise. Moreover, all the important classes of
non-associative algebras appear as coordinates: alternative algebras
for $S=\rmfa_2$, Jordan algebras for $\rmfa_1$, structurable
algebras for $\rmbc_1$ and $\rmbc_2$. \ms

In the end, what have we gained? What is emerging, is the picture
that EALAs and IARAs share many structural features with affine Lie
algebras, for example they all have concrete realizations (unlike
arbitrary Kac-Moody algebras). But they allow many more
possibilities and hence possible applications, than affine Lie
algebras. We are just at the beginning of this new theory and many
problems have not been solved (representation theory, quantization,
to name just two). \ms

{\it Contents.} The theory of affine reflection systems is described
in \S  \ref{peraffrs} after we have have laid the background for
more general reflection systems in \S  \ref{perrootsys}. Both
sections are based on the author's joint paper with Ottmar Loos
\cite{prs}. We then move to graded algebras. In \S  \ref{pergralg}
we introduce what is needed from the theory of arbitrary graded
algebras, e.g. the centroid and centroidal derivations. Following
this, we describe in \S  \ref{perrsgr} the necessary background from
the theory of root-graded Lie algebras and the subclasses of
predivision-root-graded algebras and Lie tori. These two sections
contain few new results. The new theory of invariant affine
reflection algebras and even more general algebras is presented in
\S\ref{pereala}. Finally, in \S\ref{perexass} we consider the
example ${\mathfrak s}{\mathfrak l}_n(A)$ in detail. In particular,
we describe central extensions and derivations, the important
ingredients of the EALA- and IARA-construction.
%\endsec

\smallskip\noindent \emph{Acknowledgements:}
The author thanks Bruce Allison for stimulating discussions on the
topics of this paper and the referee for helpful comments on a
preliminary version of this article.

\vfe

\section{Root systems and other types of reflection systems}
\label{perrootsys}
 %\bs

In this section and the following section \S\ref{peraffrs} we will
describe the combinatorial background for the sets of roots that
appear in the theory of extended affine Lie algebras, Lie tori and
their generalizations. See \ref{perrootsysnotes} for references. \sm

Throughout we work in vector spaces over a field $\KK$ of
characteristic $0$. For a vector space $X$ we denote by $\GL(X)$ the
general linear group of $X$, by $\Id_X$ the identity map on $X$, by
$X^* = \Hom(X,\KK)$ the dual space of $X$ and by $\lan \cdot , \cdot
\ran\colon\, X \times X^* \to \KK$ the canonical pairing given by
$\lan x,\la\ran = \la(x)$.

\subsec[prers]{\bf Pre-reflection and reflection systems.} A {\it
pre-reflection system\/}\index{pre-reflection system} is a triple
$(R,X,s)$ where $R$ is a subset of a vector space $X$ and $s\colon\,
R \to \GL(X)$, $ \al \mapsto s_\al$ is a map satisfying the
following axioms:
\begin{description} \item{(ReS0)} $0\in R$, $R$ spans $X$ and for all
$\al \in R$ we have ${\rm codim} \{ x\in X : s_\al x = x \} \le 1$
and $s_\al^2 = \Id_X$. Thus, either $s_\al = \Id_X$ or $s_\al$ is a
reflection with fixed point space of codimension $1$, leading to a
partition $R = R\re \cup R\im$ where
 \begin{eqnarray*}
  R\re & = &\{ \al \in R : s_\al \ne \Id_X\} \quad
       (\hbox{{\it reflective} or {\it real roots}})
  \index{root!reflective}
 \index{reflective root}
  \index{root!real}
 \index{real root}
 \index{R@$R\re$ (real roots of $R$)}
\\
 R\im &=& \{ \al \in R : s_\al = \Id_X\}
        \quad(\hbox{\it imaginary roots}).
  \index{root!imaginary}
 \index{imaginary root}
  \index{R@$R\im$ (imaginary roots of $R$)}
\end{eqnarray*}

\item{(ReS1)} If $\al \in R\re$ then $s_\al (\al) = -\al \ne \al$.
\sm

\item{(ReS2)} For all $\al \in R$ we have $s_\al(R\re) = R\re$ and
$s_\al (R\im) = R\im$. \end{description}

A {\it reflection system\/}\index{reflection system} is a
pre-reflection system $(R,X,s)$ for which, in addition to
(ReS0)--(ReS2), we also have
\begin{description}
\iiitem{(ReS3)} $s_{c\al} = s_\al$ whenever $0 \ne c\in \KK$ and
both $\al$ and $c\al$ belong to $R\re$. \sm

\iiitem{(ReS4)} $s_\al s_\be s_\al = s_{s_\al \be} $ for all $\al$,
$\be \in R$. \end{description}

Some explanations are in order. First, following recent practice
(e.g. \cite{aabgp} and \cite{lfrs}) we assume that $0$ is a root, in
fact $0\in R\im$ because of (ReS1). While this is not really
important, it will turn out to be convenient later on. The
terminology of ``real'' and ``imaginary" roots is of course borrowed
from the theory of Kac-Moody algebras, whose roots provide an
example of a reflection system. It is immediate from (ReS1) and
(ReS2) that
$$
  R \re = - R\re\quad \hbox{and} \quad
   s_\al (R) = R \quad\hbox{for all }\al \in R.
$$
But we do not assume $R= -R$ (this would exclude some examples in
Lie superalgebras, see \ref{inro}). \end{eia}

\subsec[sbc]{\bf Some basic concepts.} Let $(R,X,s)$ be a
pre-reflection system, often just denoted $(R,X)$ or even simply
$R$. For all pre-reflection systems we will use the same letter $s$
to denote the reflections. \sm

Let $\al \in R\re$. Because of (ReS1) there exists a unique
$\al\ch\in X^*$ such that $s_\al$ is given by the usual formula
$$ \label{prers1}
 s_\al (x) = x - \lan x,\al\ch\ran \al \quad \hbox{for }x \in X.
         \eqno(1) %\eq{prers1}
 \index{s@$s_\al$ (reflection in $\al$)}
$$
Indeed, $\al\ch$ is determined by $\al\ch(\al) = 2$ and $\Ker \al
\ch = \{x\in X : s_\al x = x\}$. So that (1) holds for all $\al \in
R$, we put $\al\ch = 0$ for $\al \in R\im$. We thus get a map ${}\ch
\colon\, R \to X^*$ which uniquely determines the map $s\colon\, R
\to \GL(X)$.  \sm

One obtains a category of pre-reflection systems, and a fortiori a
full subcategory of reflection systems, by defining a {\it
morphism\/}\index{pre-reflection system!morphism}
 $f\colon\, (R,X,s) \to (S,Y,s)$ of pre-reflection
systems as a $\KK$-linear map $f\colon\, X \to Y$ such that $f(R)
\subset S$ and $f(s_\al \be) = s_{f(\al)} f(\be)$ hold for all
$\al,\be \in R$. It is immediate that any morphism $f$ satisfies
$$ f(R\im ) \subset S\im \quad\hbox{and}\quad f(R\re) \subset
 S\re \cup \{0\}. \eqno(2)%{prers2}
$$
Note that the axiom (ReS4) of reflection systems is equivalent to
the condition that $s_\al$ be an automorphism for every $\al \in R$.
\sm

The subgroup of $\GL(X)$ generated by all reflections $s_\al$ is
called the {\it Weyl group\/}\index{Weyl group} of $R$ and denoted
$W(R)$\index{W@$W(R)$}.  \sm

A {\it subsystem\/}\index{subsystem} of $R$ is a subset $S\subset R
$ with $0\in S$ and $s_\ga \de \in S$ for all $\ga,\de \in S$. It
follows that $(S,\Span(S))$ together with the obvious restriction of
$s$ is a pre-reflection system such that the inclusion
$S\hookrightarrow R$ is a morphism. For any subsystem $S$ we have
$S\im = S\cap R\im$ and $S\re = S\cap R\re$. An example of a
subsystem is
$$ {\rm Re}(R) = R\re \cup \{0\}\index{R@${\rm Re}(R)$ (the real part of $R$)},$$
called the {\it real part of $R$}\index{pre-reflection system!real
part of}\index{real part of a pre-reflection system}. Also, for any
subspace $Y$ of $X$ the intersection $R\cap Y$ is a subsystem of
$R$. \sm

A pre-reflection system $R$ is called
\begin{description}
 \item{(i)} {\it reduced\/}
   \index{pre-reflection system!reduced}
    \index{reduced pre-reflection system}
 if $\al \in R\re$, $0 \ne c \in \KK$ and $c\al
\in R\re$ imply $c = \pm 1$ (note that we do not require this
condition for roots in $R\im$); \sm

\item{(ii)} {\it integral\/}
  \index{pre-reflection system!integral}
    \index{integral pre-reflection system}
 if $\lan R, R\ch \ran \subset \ZZ$; \sm

 \item{(iii)} {\it nondegenerate\/}
   \index{pre-reflection system!nondegenerate}
    \index{nondegenerate pre-reflection system}
if $\bigcap_{\al \in R} \Ker (\al\ch) = \{0\}$; \sm

\item{(iv)} {\it symmetric\/}
 \index{pre-reflection system!symmetric}
    \index{symmetric pre-reflection system}
if $R = -R$, equivalently, $R\im = -R\im$; \sm

\item{(v)} {\it coherent\/}
  \index{pre-reflection system!coherent}
    \index{coherent pre-reflection system}
if for $\al,\beta \in R\re$ we have $\lan \al, \beta\ch\ran = 0 \Iff
\lan \beta, \al\ch\ran = 0$, \sm

\item{(vi)} {\it tame\/}
  \index{pre-reflection system!tame}
    \index{tame pre-reflection system}
if $R\im \subset R\re + R\re$. \ms
 \end{description}
 \es

\subsec[dirsumco]{\bf Direct sums, connectedness.} Let
$(R_i,X_i)_{i\in I}$ be a family of pre-reflection systems. Put $R =
\bigcup_{i\in I} R_i$, $X=\bigoplus_{i\in I} X_i$ and define
$s\colon\, R \to \GL(X)$ by extending each $s_{\al_i}$, $\al_i \in
R_i$ to a reflection of $X$ by $s_{\al_i} | X_j = \Id$ for $i\ne j$.
Then $(R,X,s)$ is a pre-reflection system, called the {\it direct
sum of the pre-reflection systems $(R_i,X_i)$}\index{pre-reflection
system!direct sum}\index{direct sum of pre-reflection systems}.
Obviously, $R\re = \bigcup_i R\re_i$ and $R\im = \bigcup_i R\im_i$.
We will say that a pre-reflection system with a non-empty set of
real roots is {\it indecomposable\/}\index{pre-reflection
system!indecomposable}\index{indecomposable pre-reflection system}
 if it is not isomorphic to a direct sum of two
pre-reflection systems, each of which has a non-empty set of real
roots. \sm

Suppose $R$ is coherent. Two roots $\al,\be\in R\re$ are called {\it
connected\/}\index{connected roots} if there exist finitely many
roots $\al = \al_0, \al_1, \ldots, \al_n = \be$ such that $\lan
\al_i , \al_{i-1}\ch\ran \ne 0$ for $i=1, \ldots , n$, in particular
all $\al_i \in R\re$. It is easily seen that connectedness is an
equivalence relation on $R\re$. A {\it connected component of\/
${\rm Re}(R)$\/}\index{pre-reflection system!real part of!connected
component}\index{connected component of ${\rm Re}(R)$} is by
definition the union of $\{0\}$ and an equivalence class of $R\re$,
and ${\rm Re}(R)$ is called {\it connected\/}\index{pre-reflection
system!real part of!connected}\index{connected real part} if ${\rm
Re}(R)$ has only one connected component. Equivalently, ${\rm
Re}(R)$ is connected if and only if for any decomposition $R\re= R_1
\cup R_2$ with $\lan R_1 ,R_2\ch\ran = 0$ necessarily $R_1 =
\emptyset$ or $R_2 = \emptyset$. Each connected component of ${\rm
Re}(R)$ is a subsystem of $R$.
\endsubsec

\subsec[nondeg]{\bf Lemma.} {\it If $(R,X,s)$ is a nondegenerate
coherent pre-reflection system, the subsystem ${\rm Re}(R)$ is the
direct sum of its connected components. In particular, ${\rm Re}(R)$
is connected if and only if it is indecomposable.}
 \endsubsec

\subsec[ibf]{\bf Invariant bilinear forms.} Let $(R,X)$ be a
pre-reflection system. A symmetric bilinear form $b\colon\, X\times
X \to \KK $ is called {\it invariant\/}\index{pre-reflection
system!invariant bilinear form}\index{invariant!bilinear
form}\index{bilinear form!invariant} if $b(wx,wy) = b(x,y)$ holds
for all $w\in W(R)$ and $x,y\in X$. It is of course equivalent to
require $b(s_\al x,s_\al y) = b(x,y)$ for all $x,y\in X$ and $\al
\in R\re$. Moreover, invariance of $b$ is also equivalent to
$$
    2\,b(x,\al) = \lan x,\al\ch\ran \, b(\al,\al) \quad
      \hbox{for all $x\in X$ and $\al\in R\re$.} \eq{(1)}%\eq{ibf1}
$$
We will call an invariant bilinear form {\it strictly
invariant\/}\index{pre-reflection system!strictly invariant bilinear
form}\index{strictly invariant bilinear form} if (1) not only holds
for $\al \in R\re$ but for all $\al\in R$. It easily follows that an
invariant symmetric bilinear form $b$ is strictly invariant if and
only if $R\im \subset \{ x\in X : b(x,X)=0\} = : \Rad
b$\index{R@$\Rad$ (radical of a bilinear form)}, the {\it radical of
$b$\/}\index{bilinear form!radical of}. For such a bilinear form $b$
we have
$$   R \cap \Rad b = \{ \al \in R : b(\al, \al) = 0 \}
         \quad \hbox{and}\quad
    \bigcap_{\al \in R} \, \Ker \al\ch \subset \Rad b. \eq{(2)} %{ibf2}
$$
We note that we need not distinguish between invariant and strictly
invariant forms if $R= {\rm Re}(R)$, for example if $R$ is a locally
finite root system as defined in \ref{resexama}. %
\es%\endsubsec

\subsec[rostr]{\bf Root strings.} Let $R$ be a pre-reflection
system, let $\al \in R\re$, $\beta \in R$, and assume that $a:=
-\lan \beta,\al\ch\ran \in \ZZ$. The {\it $\al$-string through
$\beta$\/}\index{root string}\index{pre-reflection system!root
string} is defined by
$$
 \SS(\beta,\al):= (\beta + \ZZ\al) \cap R.
 \index{S@$\SS(\beta,\al)$ ($\al$-root sting through $\beta$)}
$$
The reflection $s_\al$ leaves\/ $\SS(\beta,\al)$ invariant and
corresponds to the reflection $i \mapsto a - i$ of
$$
         \ZZ(\beta, \al) := \{i \in \ZZ: \beta + i\al \in R\}
  \index{Z@$\ZZ(\beta,\al)$}
$$
about the point $a/2$. The set $\ZZ(\beta,\al)$ is bounded if and
only if it is bounded on one side. In this case, we put $-q=\min
\ZZ(\beta,\al)$ and $p = {\rm max}\, \ZZ(\beta,\al)$ and then have
$p,q\in \NN$ and $p - q = a = -\lan\beta,\al\ch\ran$. Moreover, the
following conditions are equivalent:
\begin{description}
\item{(i)} the {\it root string $\SS(\beta,\al)$ is unbroken\/},
\index{unbroken root string}\index{root
string!unbroken}\index{pre-reflection system!root string!unbroken}
i.e., $\ZZ(\beta,\al)$ is a finite interval in $\ZZ$ or equals
$\ZZ$,

\item{(ii)} for all $\ga \in \SS(\beta,\al)$ and all integers  $i$
between $0$ and $\lan \ga,\al\ch\ran$ we have $\ga - i \al\in
\SS(\beta,\al)$,

\item{(iii)} for all $\ga \in \SS(\beta,\al)$, $\lan \ga,\al\ch\ran
> 0$ implies $\ga - \al \in R$ and $\lan \ga, \al\ch\ran < 0$
implies $\ga + \al \in R$. \end{description}

\ms

We will now present some examples of (pre-)reflection systems. \es

\subsec[resexama]{\bf Examples: Finite and locally finite root
systems.} A root system \`a la Bourbaki (\cite[VI, \S   1.1]{brac})
with $0$ added will here be called a {\it finite root
system\/}\index{root system!finite}\index{finite root system}. Using
the terminology introduced above, a finite root system is the same
as a pre-reflection system $(R,X)$ which is integral, coincides with
its real part ${\rm Re}(R)$, and moreover satisfies the finiteness
condition that $R$ be a finite set.\sm

Replacing the finiteness condition by the local-finiteness
condition, defines a locally finite root system. Thus, a {\it
locally finite root system\/}\index{root system!locally
finite}\index{locally finite root system} is a pre-reflection system
$(R,X)$ which is integral, coincides with its real part ${\rm
Re}(R)$, and is locally-finite in the sense that $|R \cap Y|<
\infty$ for any finite-dimensional subspace $Y$ of $X$. \sm

We mention some special features of locally finite root systems,
most of them well-known from the theory of finite root systems. By
definition, a locally finite root system $R$ is integral and
symmetric. One can show that it is also a reflection system whose
root strings are unbroken and which is nondegenerate and coherent,
but not necessarily reduced. Moreover, $R$ is the direct sum of its
connected components, see \ref{nondeg}. In particular, $R$ is
connected if and only if $R$ is indecomposable, in which case $R$ is
traditionally called {\it
irreducible}\index{irreducible}\index{locally finite root
system!irreducible}. \sm

A root $\al\in R$ is called {\it divisible\/}\index{divisible
root}\index{locally finite root system!divisible
root}\index{root!divisible} or {\it indivisible\/}\index{indivisible
root}\index{locally finite root system!indivisible
root}\index{root!divisible} according to whether $\al/2$ is a root
or not. We put
$$
   R\ind  = \{0\} \cup \{ \al \in R : \al \hbox{ indivisible} \}
   \index{R@$R\ind$ (indivisible roots of $R$)}
     \quad\hbox{and} \quad
   R\div = \{ \al \in R : \al \hbox{ divisible} \}. \eq{(1)}%{resexama1}
   \index{R@$R\div$ (divisible roots of $R$)}
$$
Both $R\ind$ and $R\div$ are subsystems of $R$. \sm

A {\it root basis\/}\index{root basis}\index{locally finite root
system!root basis} of a locally finite root system $R$ is a linearly
independent set $B\subset R$ such that every $\al\in R$ is a
$\ZZ$-linear combination of $B$ with coefficients all of the same
sign. A root system $R$ has a root basis if and only if all
irreducible components of $R$ are countable (\cite[6.7, 6.9]{lfrs}).
 \sm

The set $R\ch= \{ \al\ch \in X^*: \al \in R\}$\index{R@$R\ch$
(coroot system of $R$)} is a locally finite root system in $X\ch=
\Span_\KK(R\ch)$, called the {\it coroot system of $R$}\index{coroot
system}\index{locally finite root system!coroot system}. The root
systems $R$ and $R\ch{}\ch$ are
canonically isomorphic (\cite[Th.~4.9]{lfrs}). \es%\endsubsec

\subsec[lfrsclas]{\bf Classification\index{locally finite root
system!classification}\index{classification of locally finite root
systems} of locally finite root systems.} (\cite[Th.~8.4]{lfrs})
Examples of possibly infinite locally finite root systems are the
so-called {\it classical root systems\/}\index{classical root
systems}\index{locally finite root system!classical} $\rma_I$ --
$\rmbc_I$, defined as follows. Let $I$ be an arbitrary set and let
$X_I=\bigoplus_{i\in I}\, \KK\eps_i$ be the vector space with basis
$\{\eps_i : i \in I \}$. Then
\begin{eqnarray*}
  \rma_I &=&  \{ \eps_i -\eps_j : i,j\in
I\}, \index{A@$\rma_I$ (root system of type A)}  \\
 \rmb_I &=& \{0\} \cup \{ \pm \eps_i : i\in I \} \cup
 \{ \pm \eps_i \pm \eps_j : i,j \in I, i\ne j \}
  = \{\pm \eps_i : i\in I \} \cup \rmd_I,
  \index{B@$\rmb_I$ (root system of type B)} \\
   \rmc_I &=& \{ \pm \eps_i \pm \eps_j : i,j \in I \}
     = \{\pm 2\eps_i : i\in I \} \cup \rmd_I,
 \index{C@$\rmc_I$ (root system of type C)} \\
 \rmd_I &=& \{0\} \cup \{ \pm \eps_i \pm \eps_j : i,j \in I, i\ne
 j\},\hbox{ and}  \index{D@$\rmd_I$ (root system of type D)} \\
 \rmbc_I &=&  \{ \pm \eps_i \pm \eps_j : i,j \in I \}
 \cup \{ \pm \eps_i : i\in I \}= \rmb_I \cup \rmc_I
      \index{BC@$\rmbc_I$ (root system of type BC)}
 \end{eqnarray*}
are locally finite root systems which span $X_I$, except for
$\rma_I$ which only spans a  subspace $\dot X_I$ of codimension $1$
in $X_I$. The notation $\rma$ is supposed to indicate this fact. For
finite $I$, say $|I|=n\in \NN$, we will use the usual notation
$\rmfa_n = \rma_{\{0, \ldots, n\}}$\index{A@$\rmfa_n$ (finite root
system of type A)}, but $\rmb_n$\index{B@$\rmb_n$ (finite root
system of type B)}, $\rmc_n$\index{C@$\rmc_n$ (finite root system of
type C)}, $\rmd_n$\index{D@$\rmd_n$ (finite root system of type D)},
$\rmbc_n$\index{BC@$\rmbc_n$ (finite root system of type BC)} for
$|I|=n$. The root systems $\rma_I, \ldots , \rmd_I$ are reduced, but
$R=\rmbc_I$ is not.
 \sm

Any locally finite root system $R$ is the direct limit of its finite
subsystems, which are finite root systems, and if $R$ is irreducible
it is a direct limit of irreducible finite subsystems. It is then
not too surprising that an irreducible locally finite root system is
isomorphic to exactly one of
\begin{description}
\item{(i)}  the finite  exceptional root systems $\rme_6,  \rme_7,
\rme_8, \rmf, \rmg$ or

\item{(ii)} one of the root systems ${\dot {\rm A}}_I$ $(|I|\ge
1)$, $\rmb_I$ $(|I|\ge 2)$, $\rmc_I$ $(|I|\ge 3)$, $\rmd_I$ $(|I|\ge
4)$ or $\rmbc_I$ $(|I|\ge 1)$. \end{description} \noindent
Conversely, all the root systems listed in (i) and (ii) are
irreducible. \sm

The classification of locally finite root systems is independent of
the base field $\KK$, which could therefore be taken to be $\QQ$ or
$\RR$. But since the root systems which we will later encounter
naturally ``live'' in $\KK$-vector spaces, fixing the base field is
not convenient. \sm

Locally finite root systems can also be defined via invariant
nondegenerate symmetric bilinear forms: \endsubsec%\es

\subsec[lfrschar]{\bf Proposition} (\cite[2.10]{prs}) {\it For an
integral pre-reflection system $(R,X)$ the following conditions are
equivalent:
\begin{description}
\item{\rm (i)} $R$ is a locally finite root system;

\item{\rm (ii)} there exists a nondegenerate strictly invariant form on $(R,X)$,
and for every $\al\in R$ the set $\lan R, \al\ch\ran$ is bounded as
a subset of $\ZZ$.
\end{description}

\noindent In this case, $(R,X)$ has a unique invariant bilinear form
$\inpr$ which is {\rm normalized\/}\index{normalized
form}\index{locally finite root system!normalized form} in the sense
that $ 2 \in \{(\al|\al) : 0 \ne \al\in C\} \subset \{2,4,6,8\}$ for
every connected component $C$ of $R$. The normalized form is
nondegenerate in general and
positive definite for $\KK=\RR$.}  %\es%
\es%\endsubsec

\subsec[resexamd]{\bf Example: Reflection systems associated to
bilinear forms.}\index{reflection system!associated to bilinear
forms} Let $R$ be a spanning set of a vector space $X$ containing
$0$, and let $\inpr$ be a symmetric bilinear form on $X$. For $\al
\in X$ we denote the linear form $x \mapsto (\al|x)$ by $\al^\flat$.
Let $\Phi\subset \{\al\in R : (\al | \al) \ne 0\}$, define $\ch
\colon\, R \to X^*$ by
%$$
%\al\ch := \begin{cases}{ {2\al^\flat \over (\al| \al)} & if $\al\in
%\Phi$\cr \noalign{\medskip } \quad 0 & otherwise},
%  \end{cases}
% \eq{aff3}
%$$
$$% \begin{equation*}
 \al\ch = \begin{cases}\frac{2 \al^\flat}{(\al | \al)}
           & \text{if $\al \in \Phi$,} \\
            0 & \text{otherwise}
           \end{cases}   \eq{(1)}
$$ % \end{equation*}
and then define $s_\al$ by \ref{prers1}. Thus
$s_\al$ is the orthogonal reflection in the hyperplane $\al^\perp$
if $\al\in \Phi$, and the identity otherwise. If $s_\al(R) \subset
R$ and $s_\al(\Phi) \subset \Phi$ for all $\al \in \Phi$, then
$(R,X,s)$ is a coherent reflection system with the given subset
$\Phi$ as set of reflective roots. The bilinear form $\inpr$ is
invariant, and it is strictly invariant if and only if $R\im = R
\cap \Rad \, \inpr$, i.e., $\inpr$ is affine in the sense of
\ref{affo}. \sm

By \ref{lfrschar} every locally finite root system is of this type.
But also the not necessarily crystallographic finite root systems,
see for example \cite[1.2]{hum:cox}, arise in this way.
The latter are not necessarily integral reflection systems. %
\es%\endsubsec

\subsec[inro]{\bf Example: Roots in Lie algebras with toral
subalgebras.} Let $L$ be a Lie algebra over a ring $k$ containing
$\frac{1}{2}$ (this generality will be used later). Following
\cite[Ch.VIII, \S  11]{bou:lie78} we will call a non-zero triple
$(e, h, f)$ of elements of $L$ an {\it ${\mathfrak s}{\mathfrak
l}_2$-triple\/}\index{S@${\mathfrak s}{\mathfrak l}_2$-triple} if
$$[e,f] = - h, \quad [h,\, e] = 2 e \quad
 \hbox{and} \quad [h,\, f] = -2f. $$
For example, in
$${\mathfrak s}{\mathfrak l}_2(k)   = \Big\{ \ts
       \begin{pmatrix} a & b \\ c & -a \end{pmatrix}
               : a,b,c\in k \Big\}
$$
the elements
$$
  e_{{\mathfrak s}{\mathfrak l}_2}
         = \begin{pmatrix}0 & 1 \\ 0 &0 \end{pmatrix}, \quad
 h_{{\mathfrak s}{\mathfrak l}_2}
         = \begin{pmatrix} 1 & 0 \\ 0 &-1 \end{pmatrix}, \quad
 f_{{\mathfrak s}{\mathfrak l}_2}
        = \begin{pmatrix} 0 & 0 \\ -1 &0 \end{pmatrix},
 \eq{(1)}
$$ form an ${\mathfrak s}{\mathfrak l}_2$-triple in ${\mathfrak s}{\mathfrak l}_2(k)$. In general, for any ${\mathfrak s}{\mathfrak l}_2$-triple
$(e,h,f)$ in a Lie algebra $L$ there exists a unique Lie algebra
homomorphism $\ph \colon\, {\mathfrak s}{\mathfrak l}_2(k) \to L$
mapping the
matrices in (1) %\ref{Heisl23}
onto the corresponding elements $(e,h,f)$ in $L$. \sm

Let now $k=\KK$ be a field of characteristic $0$, and let $L$ be a
Lie algebra over $\KK$ and $T\subset L$ a subspace. We call
$T\subset L$ a {\it toral subalgebra\/}\index{toral
subalgebra}\index{subalgebra!toral}, sometimes also called an
$\ad$-diagonalizable
subalgebra\index{subalgebra!$\ad$-diagonalizable}, if
$$
   L= \ts\bigoplus_{\al \in T^*} L_\al (T)\quad\hbox{where}\quad
     L_\al (T) = \{ x\in L : [t,x]=\al(t)x \hbox{ for all } t\in T\}
\eq{(2)} %{inro1}
$$ for any $\al \in T^*$. In this case, (2) %\ref{inro1}
is referred to as the {\it root space decomposition\/}\index{root
space decomposition} of $(L,T)$ and $R= \{ \al \in T^* : L_\al(T)
\ne 0\}$ as the {\it set of roots of $(L,T)$}. We will usually
abbreviate $L_\al = L_\al(T)$\index{L@$L_\al$ ($\al$-root space of
$L$)}. Any toral subalgebra $T$ is abelian, thus $T \subset L_0$ and
$0 \in R$ unless $T=\{0\}=E$ (which is allowed but not interesting).
We will say that a toral $T$ is a {\it splitting Cartan
subalgebra\/}\index{Cartan subalgebra!splitting}\index{splitting
Cartan subalgebra} if $T=L_0$. \sm

We denote by $R\int$\index{R@$R\int$ (integrable roots of $R$)} the
subset of {\it integrable
roots}\index{root!integrable}\index{integrable root} of $(L,T)$,
i.e., those $\al\in R$ for which there exists an ${\mathfrak
s}{\mathfrak l}_2$-triple $(e_\al,h_\al,f_\al) \in L_\al \times T
\times L_{-\al}$ such that the adjoint maps $\ad e_\al$ and $\ad
f_\al$ are locally nilpotent. \sm

Let $X= \Span_\KK(R) \subset T^*$ and for $\al \in R$ define $s_\al
\in \GL(X)$ by
$$
s_\al (x) = \begin{cases}
          x - x(h_\al)\al, & \al \in R\int, \\
                    x & \text{otherwise}. \end{cases}
$$
It is a straightforward exercise in ${\mathfrak s}{\mathfrak
l}_2$-representation theory to verify that with respect to the
reflections $s_\al$ defined above, {\it $(R,X)$ is an integral
pre-reflection system, which is reduced in case $T$ is a splitting
Cartan subalgebra.} A priori, $s_\al$ will depend on $h_\al$, but
this is not so for the Lie algebras which are of main interest here,
the invariant affine reflection algebras of \S\ref{pereala}. \sm

We note that several important classes of pre-reflection systems
arise in this way from Lie algebras, e.g., locally finite root
systems \cite{n:3g}, Kac-Moody root systems \cite{kac}, \cite{mp}
and extended affine root systems \cite{aabgp}, see \ref{affex}. An
example of the latter is described in detail in \ref{afrs}. We point
out that in general $R$ is not symmetric. For example, this is so
for the set of roots of classical Lie superalgebras \cite{kac:sup1}
(they can be viewed in the setting above by forgetting the
multiplication of the odd part).
\es%\endsubsec

\subsec[afrs]{\bf Example: Affine root systems.} Let $S$ be an
irreducible locally finite root system. By \ref{lfrschar}, $S$ has a
unique normalized invariant form, with respect to which we can
introduce {\it short roots $S\sh= \{\al \in S : (\al|\al) =
2\}$\/}\index{S@$S\sh$ (short roots)}\index{short
root}\index{root!short} and {\it long roots\/} $S\lg =\{ \al \in S :
(\al|\al) = 4 \hbox{ or } 6\} = S\setminus (S\sh\cup
S\div)$\index{root!long}\index{long root}\index{S@$S\lg$ (long
roots)}, where $S\div$ are the divisible roots of \ref{resexama}.
Both $S\sh\cup \{0\}$ and $S\lg \cup \{0\}$ are subsystems of $S$
and hence locally finite root systems. Note that $S = \{0\} \cup
S\sh$ if $S$ is simply-laced, $(S\div \setminus \{0\}) \ne \emptyset
\Leftrightarrow S = \rmbc_I$, and if $S\lg \ne \emptyset$ the set of
long roots is given by $S\lg = \{ \al \in S : (\al |\al) = k(S)\}$
for
$$
 k(S) = \begin{cases}
           2, & S=\rmf \text{ or of type }
                    \rmb_I, \rmc_I, \rmbc_I, |I|\ge 2, \\
                   3  &  S= \rmg. \end{cases} \eq{(1)}%{afrs0}
\index{kS@$k(S)$} $$ \sm

Let $X = \Span( S)\oplus \KK \de$ and let $\inpr$ be the symmetric
bilinear form which restricted to $\Span(S)$ is the normalized
invariant form of $S$ and which satisfies $(X \,|\,\KK \de) = 0$.
Furthermore, choose a {\it tier number\/}\index{tier number}
$t(S)=1$ or $t(S) = k(S)$, with $t(S) = 1$ in case $S = \rmbc_I$,
$|I|\ge 2$. Finally, put
%$$\eqalign{
\begin{eqnarray*}
 \Phi &=& \Big( \bigcup_{\al \in S\sh}\, (\al \oplus \ZZ \de) \Big)
 \cup \Big( \bigcup_{\al \in S\lg }\,
     \big(\al \oplus t(S) \,\ZZ \de \big) \Big)
 \cup \Big( \bigcup_{0\ne \al \in S\div}\, \big(
            \al + (1 + 2\ZZ) \de \big) \Big), \\
  R &=& \ZZ \de \cup \Phi.
\end{eqnarray*}
 %{afrs1}$$
Of course, if $S\lg=\emptyset$ or $S\div=\{0\}$ the corresponding
union above is to be interpreted as the empty set. Note $(\be | \be)
\ne 0$ for any $\be \in \Phi$. One can easily verify that every
reflection $s_\be$ defined in \ref{resexamd} leaves $\Phi$ and $\ZZ
\de$ invariant, hence $R= R(S, t(S))$ is a symmetric reflection
system, called the {\it affine root system associated to $S$ and
$t(S)$\/}.\index{affine root system}\index{root system!affine} \sm

The following table lists all possibilities for $(S, t(S))$ where
$S$ is locally finite. For finite $S$ of type ${\rm X}_I$, say
$|I|=l$, we also include the corresponding labels of the affine root
system $R(S, t(S))$ as defined \cite{mp} (third column) and
\cite{kac} (fourth column).
$$
\vcenter{\offinterlineskip
\def\strut{\vrule height 14 pt depth 6 pt width 0 pt}%
\everycr={\noalign{\hrule}}%
\halign{% Musterzeile:
\strut\vrule\hfil\enspace$#$\enspace\hfil\vrule\vrule % 1. Spalte
&\hfil\enspace$\;\;#\;\;$\enspace\hfil\vrule
&\hfil\enspace$\;\;#\;\;$\enspace\hfil\vrule
&&\hfil\enspace$\;\;#\;\;$\enspace\hfil\vrule% der && bewirkt, dass alle folgenden
                                       % Spalten genauso behandelt werden.
\cr                                    % Ende der Musterzeile
%%  1. Zeile: Ueberschriften:
 S &     t(S) &  \hbox{affine label \cite{mp}} &
         \hbox{affine label \cite{kac}} \cr  %  Ende der ersten Zeile
%%  Nun die eigentliche Tabelle:
 \noalign{\hrule height .8 pt}% zusaetzliche dicke horizontale Linie
 \hbox{reduced}&1 & S^{(1)} & S^{(1)} \cr %Ende der dritten Zeile
 \rmb_I\; (|I| \ge 2) & 2 & \rmb_l^{(2)} & \rmd_{l+1}^{(2)} \cr
 \rmc_I \; (|I|\ge 3) & 2 & \rmc_l^{(2)} & \rmfa_{2l-1}^{(2)} \cr
 \rmf & 2 & \rmf^{(2)} & \rme_6^{(2)} \cr
 \rmg & 3 & \rmg^{(3)} & \rmd_4^{(3)} \cr
 \rmbc_1 & - & \rmbc_1^{(2)} & \rmfa_2^{(2)}\cr
 \rmbc_I \, (|I| \ge 2) & 1 & \rmbc_l^{(2)} & \rmfa_{2l}^{(2)} \cr
}}
$$
If $S$ is reduced then $R= S \oplus \ZZ\de$, i.e., $t(S)=1$ in the
definition of $\Phi$, %\ref{afrs1},
is a so-called {\it untwisted affine root system}.\index{untwisted
affine root system}\index{affine root system!untwisted}\index{root
system!affine!untwisted}
\endsubsec%\sm

\subsec[afrreal]{\bf Affine Lie algebras.}\index{affine Lie
algebra}\index{Lie algebra!affine} Every affine root system appears
as the set of roots of a Lie algebra $E$ with a splitting Cartan
subalgebra $H$ in the sense of \ref{inro} (the letters $(E,H)$ are
chosen in anticipation of the definition of an extended affine Lie
algebra \ref{ealadef}). This is easily verified in the untwisted
case. \sm

Indeed, let $\frg$ be a finite-dimensional split simple Lie
$\KK$-algebra, e.g., a finite-dimensional simple Lie algebra over an
algebraically closed field $\KK$, or let $\frg$ be a centreless
infinite rank affine algebra \cite[7.11]{kac}. Equivalently, $\frg$
is a locally finite split simple Lie algebra, classified in
\cite{neeb-stu} and \cite{stu:struct}, or $\frg$ is the
Tits-Kantor-Koecher algebra of a Jordan pair spanned by a connected
grid, whose classification is immediate from \cite{n:3g}. In any
case, $\frg$ contains a splitting Cartan subalgebra $\frh$ such that
the set of roots of $(\frg, \frh)$ is an irreducible locally finite
root system $S$, which is finite if and only if $\frg$ is
finite-dimensional. Moreover, $\frg$ carries an invariant
nondegenerate symmetric bilinear form $\ka$, unique to a scalar (in
the finite-dimensional case $\ka$ can be taken to be the Killing
form). The Lie algebra $E$ will be constructed in three steps: \sm

(I) $L=\frg \ot \KK[t^{\pm 1}]$ is the so-called {\it untwisted loop
algebra of $\frg$\/}\index{untwisted loop algebra}\index{loop
algebra!untwisted}. Its Lie algebra product $[\cdot,\cdot]_L$ is
given by
$$[x \ot p, \, y \ot q]_L = [x,  y]_\frg \ot pq$$ for $x, y \in \frg$
and $p, q \in \KK[t^{\pm 1}]$. Although $L$ could be viewed as a Lie
algebra over $\KK[t^{\pm 1}]$, we will view $L$ as a Lie algebra
over $\KK$. \sm

(II) $K= L \oplus \KK c$ with Lie algebra product $[.,.]_K$ given by
$$[l_1 \oplus s_1 c,\, l_2 \oplus s_2 c]_K = [l_1,\, l_2]_L \oplus \si(l_1,
l_2) c $$ for $l_i \in L$, $s_i \in \KK$ and $\si\colon\, L\times L
\to \KK$ the $\KK$-bilinear map determined by $\si(x_1 \ot t^m, \,
x_2 \ot t^n) = \de_{m+n,0} \,  m\, \ka(x_1, x_2)$. \sm

(III) $E= K \oplus \KK d$. Here $d$ is the derivation of $K$ given
by $d(x\ot t^m) = m x \ot t^m$, $d(c) = 0$, and $E$ is the
semidirect product of $K$ and $\KK d$, i.e., $$[k_1 \oplus s_1 d, \,
k_2 \oplus s_2 d]_E = [k_1, \, k_2]_K + s_2 d(k_2)- s_1 d(k_1).$$
One calls $E$ the {\it untwisted affine Lie algebra associated to
$\frg$}. \sm

Denote the root spaces of $(\frg, \frh)$ by $\frg_\al$, so that
$\frh = \frg_0$. Then
  $$     H = (\frh \ot \KK1) \oplus \KK c \oplus \KK d
         %\eq{affreal1}
 $$
is a toral subalgebra of $E$ whose set of roots in $E$ is the
untwisted affine root system associated to the root system $S$ of
$(\frg, \frh)$: $E_0  = H$, $E_{m\de} = \frh \ot \KK t^m$ for $0\ne
m \in \ZZ$ and $E_{\al \oplus m\de}= \frg_\al \ot t^m $ for $0\ne
\al \in S$ and $m\in \ZZ$. \ms

One can also realize the {\it twisted affine root
systems\/},\index{root system!affine!twisted}\index{twisted affine
root system}\index{affine root system!twisted} i.e. those with
$t(S)>1$, by replacing the untwisted loop algebra $L$ in (I) by a
twisted version. To define it in the case of a finite-dimensional
$\frg$, let $\si$ be a non-trivial diagram automorphism of $\frg$ of
order $r$. Hence $\frg$ is of type $\rmfa_l \; (l\ge 2, r=2)$,
$\rmd_l \;(l\ge 4, r=2)$, $\rmd_4 \;(r=3)$ or $\rme_6 \;(r=2)$. For
an infinite-dimensional $\frg$ (where diagrams need not exist) we
take as $\si$ the obvious infinite-dimensional analogue of the
corresponding matrix realization of $\si$ for a finite-dimensional
$\frg$, e.g. $x \mapsto - x^t$ for type $\rmfa$. The {\it twisted
loop algebra $L$\/}\index{twisted loop algebra}\index{loop
algebra!twisted} is the special case $n=1$ of the multi-loop algebra
$\scM_r(\frg, \si)$ defined in \ref{perweakv}. The construction of
$K$ and $E$ are the same as in the untwisted case. Replacing the
untwisted $\frh$ in the definition of  $H$ above %\ref{affreal1}
by
the subspace $\frh^\si$ of fixed points under $\si$ defines a
splitting Cartan subalgebra of $E$ in the twisted case, whose set of
roots is a twisted affine root systems and all twisted affine root
systems arise in this way, see \cite[8.3]{kac} for details. \ms

In the later sections we will generalize affine root systems and
affine Lie algebras. Affine root systems will turn out to be
examples of extended affine root systems and locally extended affine
root systems (\ref{affex}), namely those of nullity $1$. Extended
and locally extended affine root systems are in turn special types
of affine reflection systems, which we will describe in the next
section. In the same vein, untwisted affine Lie algebras are
examples of (locally) extended affine Lie algebras, which in turn
are examples of affine reflection Lie algebras which we will study
in \S  \ref{pereala}. From the point of view of extended affine Lie
algebras, \ref{ealadef}, the Lie algebras $K$ and $L$ are the core
and centreless core of the extended affine Lie algebra $(E,H)$. As
Lie algebras, they are Lie tori, defined in \ref{rogr}. The
construction from $L$ to $E$ described above is  a special case of
the construction \ref{ealacon}, see the example in \ref{ealamt}.
\es%\endsubsec

\subsec[perrootsysnotes]{\bf Notes.} With the exceptions mentioned
in the text and below,  all results in this section are proven in
\cite{lfrs} and \cite{prs}. Many of the concepts introduced here are
well-known from the theory of finite root systems. The notion of
tameness (\ref{sbc}) goes back to \cite{ABY}. It is equivalent to
the requirement that $R$ not have isolated roots, where $\be \in R$
is called an isolated root (with respect to $R\re$), if $\al + \be
\not\in R$ for all $\al \in R\re$. In the setting of the Lie
algebras considered in \S  \ref{pereala}, tameness of the root
system is a consequence of tameness of the algebra, see \ref{invth}.
\sm

Locally finite root systems are studied in \cite{lfrs} over $\RR$.
But as already mentioned in \cite[4.14]{lfrs} there is a canonical
equivalence between the categories of locally finite root systems
over $\RR$ and over any field $\KK$ of characteristic $0$. In
particular, the classification is the ``same'' over any $\KK$. The
classification itself is proven in \cite[Th.~8.4]{lfrs}.  A
substitute for the non-existence of root bases are grid bases, which
exist for all infinite reduced root systems (\cite{n:cr}). \sm

Our definition of an ${\mathfrak s}{\mathfrak l}_2$-triple follows
Bourbaki (\cite[Ch.VII, \S  11]{bou:lie78}). It differs from the one
used in other texts by a sign. Replacing $f$ by $-f$ shows that the
two notions are equivalent. Bourbaki's definition is more natural in
the setting here and avoids some minus signs later. The concept of
an integrable root (\ref{inro}) was introduced by Neeb in
\cite{neeb:split}. \sm

For a finite $S$, affine root systems (\ref{afrs}) are determined in
\cite[Prop.~6.3]{kac}. The description given in \ref{afrs} is
different, but of course equivalent to the one in loc. cit. It is
adapted to viewing affine root systems as extended affine root
systems of nullity $1$. In fact, the description of affine root
systems in \ref{afrs} is a special case of the Structure Theorem for
extended affine root systems \cite[II, Th.~2.37]{aabgp}, keeping in
mind that a semilattice of nullity $1$ is a lattice by \cite[II,
Cor.~1.7]{aabgp}. The table in \ref{afrs} reproduces
\cite[Table~1.24]{abgp}. \sm

Other examples of reflection systems are the set of roots associated
to a ``root basis'' in the sense of H\'ee \cite{hee}. Or, let $R$ be
the root string closure of the real roots associated to root data
\`a la Moody-Pianzola \cite[Ch.~5]{mp}. Then $R$ is a symmetric,
reduced and integral reflection system. \sm

The Weyl group of a pre-reflection system is in general not a
Coxeter group. For example, this already happens for locally finite
root systems, see \cite[9.9]{lfrs}. As a substitute, one can give a
so-called presentation by conjugation, essentially the relation
(ReS4), see \cite[5.12]{lfrs} for locally finite root systems and
\cite{HofG} for Weyl groups of extended affine root systems (note
that the Weyl groups in \cite{HofG} and also in \cite{aabgp} are
defined on a bigger space than the span of the roots and hence our
Weyl groups are homomorphic images of the Weyl groups in loc. cit.).
\endsubsec %\endsec\endinput

\vfe

%\sec[peraffrs]{Affine reflection systems}
\section{Affine reflection systems}\label{peraffrs}

 In this section we
describe extensions of pre-reflection systems and use them to define
affine reflection systems, which are a generalization of extended
affine root systems. The roots of the Lie algebras, which we will be
studying in \S\ref{pereala}, will turn out to be examples of affine
reflection systems. Hints to references are given in
\ref{peraffrsnotes}.\bs

\subsec[parsecdef]{\bf Partial sections.} Let $f\colon\, (R,X) \to
(S,Y)$ be a morphism of pre-reflection systems with $f(R) = S$, and
let $S' \subset S$ be a subsystem spanning $Y$.  A {\it partial
section of $f$ over $S'$\/}\index{partial
section}\index{pre-reflection system!morphism!partial section} is a
morphism $g\colon\, (S',Y) \to (R,X)$ of pre-reflection systems such
that $f\circ g = \Id_Y$. The name is (partially) justified by the
fact that a partial section leads to a partial section of the
canonical epimorphism $W(R) \to W(S)$, namely a section defined over
$W(S')$. \sm

Let $f\colon\, (R,X) \to (S,Y)$ be a morphism of pre-reflection
systems satisfying $f(R)=S$. In general, a partial section of $f$
over all of $S$ need not exist. However, partial sections always
exist in the category of reflection systems.
\endsubsec

\subsec[ext]{\bf Extensions.} Recall \ref{sbc}: $f(R\im) \subset
S\im$ and $f(R\re) \subset \{0\} \cup  S\re$ for any morphism
$f\colon\, (R,X) \to (S, Y)$ of pre-reflection systems. We call $f$
an {\it extension\/}\index{extension!of pre-reflection
systems}\index{pre-reflection system!extension} if $f(R\im) = S\im$
and $f(R\re) = S\re$. We will say that {\it $R$ is an extension of
$S$\/} if there exists an extension $f\colon\, R \to S$.\sm

We mention some properties of an extension $f\colon\, R \to S$. By
definition, $f$ is surjective. Also, $R$ is coherent if and only if
$S$ is so, and in this case $f$ induces a bijection $C \mapsto f(C)$
between the set of connected components of ${\rm Re}(R)$ and of
${\rm Re}(S)$. In particular,
$$ \hbox{${\rm Re}(R)$ is connected $\Iff {\rm Re}(S)$ is connected
$\Iff {\rm Re}(R)$ is indecomposable.}$$
 Moreover, $R$ is integral if and
only if $S$ is so. Finally, $f$ maps a root string $\SS(\beta,\al)$,
$\beta \in R$, $\al\in R\re$, injectively to $\SS\big(f(\beta),
f(\al)\big)$. \sm

If $R$ is an extension of a nondegenerate $S$, e.g. a locally finite
root system, then $S$ is unique up to isomorphism and will be called
the {\it quotient pre-reflection system  of $R$\/}.\index{quotient
pre-reflection system}\index{pre-reflection system!quotient} On the
other hand, if $R$ is nondegenerate every extension $f\colon\, R \to
S$ is injective, hence an isomorphism. In particular, a locally
finite root system $R$ does not arise as a non-trivial extension of
a pre-reflection system $S$. But as we will see below, a locally
finite root system does have many interesting extensions. \sm

A first example of an extension is the canonical projection $R\to S$
for $R$ an affine root system associated to a locally finite
irreducible root system $S$, \ref{afrs}. In this case, $S$ is the
quotient pre-reflection system = quotient root system of $R$. %
\es%\endsubsec

\subsec[gexda]{\bf Extension data.} Let $(S,Y)$ be a pre-reflection
system, let $S'$ be a subsystem of $S$ with $\Span(S')=Y$ and let
$Z$ be a $\KK$-vector space. A family $\frL =(\La_\xi)_{\xi \in S}$
of nonempty subsets of $Z$ is called an {\it extension datum of type
$(S,S',Z)$\/}\index{extension!datum}\index{extension datum}  if
\smallskip

(ED1)\ \ for all $\xi,\eta \in S$ and all $\la \in \La_\xi$, $\mu
\in \La_\eta$ we have $ \mu - \lan \eta , \xi\ch\ran \la \in
   \La_{s_\xi(\eta)}$, \smallskip

(ED2)\ \ $0\in \La_{\xi'}$ for all $\xi' \in S'$, and \smallskip

(ED3)\ \ $Z$ is spanned by the union of all $\La_\xi$, $\xi\in S$.
\smallskip

The only condition on $\La_0$ is (ED2): $0 \in \La_0$. Moreover,
$\La_0$ is related to the other $\La_\xi$, $0\ne \xi \in S$, only
via the axiom (ED3). As (ED3) only serves to determine $Z$ and can
always be achieved by replacing $Z$ by the span of the $\La_\xi$,
$\xi \in S$, it follows that one can always modify a given extension
datum by replacing $\La_0$ by some other set containing $0$.
However, as \ref{affdef} shows, this may change the properties of
the associated pre-reflection system. \sm

Any extension datum $\frL$ has the following properties, where
$W_{S'} \subset W(S)$ is the subgroup generated by all $s_{\xi'}$,
$\xi' \in S'$.
 \begin{eqnarray*}
  \La_{-\xi}& = & - \La_\xi \quad\hbox{for all }\xi \in S\re,
          \\
 2\La_\xi - \La_\xi &\subset& \La_\xi \quad\hbox{for all } \xi \in
    S\re,  \\
  \La_\eta &=&  \La_{w'(\eta)} \quad\hbox{for all
    $\eta \in S$ and $w'\in   W_{S'}$}\,,  \\
   \La_\eta - \lan \eta, \xi'\phantom{}\ch\ran \La_{\xi'}
   &\subset& \La_\eta
        \quad \hbox{for } \xi'\in S', \;\eta \in S
     \quad\hbox{and}  \\
   \La_{\xi'} &=& \La_{-\xi'} \quad \hbox{for } \xi'\in S'.
  \end{eqnarray*}
%$$\eqalignno{
% \La_{-\xi}& =- \La_\xi \quad\hbox{for all }\xi \in S\re,
%         &\eq{gexda1}  \cr
% 2\La_\xi - \La_\xi &\subset \La_\xi \quad\hbox{for all } \xi \in
%    S\re, &\eq{gexda2}\cr
%  \La_\eta &= \La_{w'(\eta)} \quad\hbox{for all
%    $\eta \in S$ and $w'\in   W_{S'}$}\,, &\eq{gexda3}\cr
%   \La_\eta - \lan \eta, \xi'\phantom{}\ch\ran \La_{\xi'} &\subset \La_\eta
%        \quad \hbox{for } \xi'\in S', \;\eta \in S
%     \quad\hbox{and}  &\eq{gexda4}\cr
%   \La_{\xi'} &=\La_{-\xi'} \quad \hbox{for } \xi'\in S'.
%&\eq{gexda4a} }$$
In particular, $\La_\xi$ is constant on the $W_{S'}$-orbits of $S$.
However, in general, the $\La_\xi$ are not constant on all of $S$,
see the examples in \ref{affexd}. \sm

One might wonder if the conditions (ED1)--(ED3) are strong enough to
force the $\La_\xi$ to be subgroups of $(Z,+)$. But it turns out
that this is not the case. For example, this already happens for the
extension data of locally finite root systems, \ref{affexd}. On the
other hand, if $R$ is an integral pre-reflection system and $\La$ a
subgroup of $(Z,+)$ which spans $Z$, then $\La_\xi \equiv \La$ is an
example of an extension datum, the so-called {\it untwisted\/} case,
cf. the definition of an untwisted affine root system in \ref{afrs}.
\sm

The conditions above have appeared in the context of reflection
spaces (not to be confused with a reflection system). Recall
(\cite{l:sp}) that a {\it reflection space}\index{reflection space}
is a set $S$ with a map $S\times S \to S : (s,t) \mapsto s\cdot t$
satisfying $s\cdot s = s$, $s\cdot (s\cdot t) = t$ and $s\cdot
(t\cdot u) = (s\cdot t) \cdot (s\cdot u)$ for all $s,t,u\in S$. A
{\it reflection subspace\/}\index{reflection subspace} of a
reflection space $(S, \cdot)$ is a subset $T\subset S$ such that
$t_1\cdot t_2 \in T$ for all $t_1, t_2 \in T$. In this case, $T$ is
a reflection space with the induced operation. Any abelian group
$(Z,+)$ is a reflection space with respect to the operation $x\cdot
y=2x-y$ for $x,y\in Z$. Correspondingly, a reflection subspace of
the reflection space $(Z, \cdot)$ is a subset $A\subset Z$
satisfying $2a -b\in A$ for all $a,b\in A$, symbolically $2A-A
\subset A$. Hence all $\La_\xi$, $\xi\in S\re$, are reflection
subspaces of $(Z,\cdot)$. Moreover, for any subset $A$ of $Z$ it is
easily seen that
$$
  A-2A \subset A \;\Leftrightarrow\; A=-A \hbox{ and } 2A+A \subset
       A \;\Leftrightarrow\; A =
  -A \hbox{ and } 2A- A \subset A.   \eq{(1)} %{gexda5}
$$
A subset $A$ satisfying (1) %\no{gexda5}
will be called a {\it symmetric reflection
subspace\/}.\index{reflection subspace!symmetric}\index{symmetric
reflection subspace} We will consider $0$ as the base point of the
reflection space $Z$. Also, we denote by $\ZZ[A]$ the subgroup of
$(Z,+)$ generated by $A\subset Z$. Then the following are equivalent
(\cite[2.1]{NY}):
%\smallskip

 %\iitem{(i)}
 (i)  $0\in A$ and $A-2A \subset A$,
 %\iitem{(ii)}

  (ii) $0\in A$ and $2A-A \subset A$,

 %\iitem{(iii)}
 (iii) $2\ZZ[A]\subset A$ and $2\ZZ[A] -A \subset A$,

 %\iitem{(iv)}
 (iv) $A$ is a union of cosets modulo $2\ZZ[A]$, including
 the trivial coset $2\ZZ[A]$.\smallskip

\noindent In this case, $A$ will be   called   a {\it pointed
reflection subspace}.\index{reflection
subspace!pointed}\index{pointed reflection subspace} It is immediate
from the above that every $\La_{\xi'}$, $\xi'\in S'{}\re$, is a
pointed reflection subspace. We note that pointed reflection
subspaces are necessarily symmetric.
% (which follows for the $\La_{\xi'}$ from \no{gexda1} and \no{gexda4}).
It is obvious from (iv) above that a pointed reflection subspace is
in general not a subgroup of $(Z,+)$. \ms

 The following theorem characterizes extensions in terms of
extension data. \es

\subsec[extchar]{\bf Theorem} {\it Let $(S,Y)$ be a pre-reflection
system. \sm

{\rm (a)} Let $\frL=(\La_\xi)_{\xi \in S}$ be an extension datum of
type $(S,S',Z)$. Put $X := Y \oplus Z$, denote by $\pi\colon\, X \to
Y$ the projection with kernel $Z$, and define
$$
 R := \bigcup_{\xi\in S} \xi \oplus \La_{\xi}\subset X
         \quad \hbox{and} \quad
 s_\al(x) :=  s_\xi(y) \oplus \big( z - \lan y,\xi\ch\ran \la \big),
  \eq{(1)} %{gexdapdef}
$$
for all $\al = \xi \oplus \la \in \xi \oplus \La_\xi \subset R$ and
all $x = y \oplus z \in X$. Then $R$ is a pre-reflection
 system   in $X$, denoted  $\scE=\scE(S, S', \frL)$. Moreover, $\pi \colon\,
(R,X) \to (S,Y)$  is an extension of pre-reflection systems, and the
canonical injection $\io\colon\, Y \to X$ is a partial section of
$\pi$ over~$S'$.
\medskip

{\rm (b)}  Conversely, let $f\colon\, (R,X) \to (S,Y)$ be an
extension and let $g\colon\, S' \to R$ be a partial section of $f$,
cf.\/ {\rm \ref{parsecdef}}. For every $\xi\in S$ define $R_\xi
\subset R$ and $\La_\xi \subset Z:=\Ker(f)$ by
$$ \label{lxi}
R_\xi = R\cap f^{-1}(\xi) = g(\xi) \oplus \La_{\xi}.  \eq{(2)}%{lxi}
$$
Then $\frL = (\La_{\xi})_{\xi\in S}$ is an extension datum of type
$(S,S',Z)$, and the vector space isomorphism $\ph\colon\, Y \oplus Z
\cong X$ sending $y\oplus z$ to $g(y) \oplus z$ is an isomorphism
$\scE(S,S', \frL) \cong R$ of pre-reflection systems making the
following diagram commutative:
$$
\vcenter{ \xymatrix@R=1.5pc{
 & S'\ar[ld]_\iota \ar[rd]^g & \cr
   \scE \ar[rr]^\ph \ar[dr]_\pi && R \ar[dl]^f \cr
& S & \cr }}.  %\eq{extchar1}
$$

{\rm (c)}  In the setting of {\rm (b)}, the following are equivalent
for $g'\in \Hom_\KK(Y,Z)$: \sm

\begin{description} %{\leftmargin}{2cm}
%\iiitem{\rm (i)}
\item{\rm (i)} $g'\colon\, S' \to R$ is another partial section of
$f$,

%\iiitem{\rm (ii)}
\item{\rm (ii)}  there exists $\ph\in \Hom_\KK(Y,Z)$ such that $g'=g
+ \ph$ and $\ph(\xi') \in \La_{\xi'}$ for all $\xi'\in S'$.
\end{description}

\noindent In this case, the extension datum $\frL' =
(\La'_{\xi})_{\xi \in S}$ defined by {\rm (2)} with respect to $g'$
is related to the extension datum $\frL$  by
$$
    \La'_\xi = \La_\xi - \ph(\xi) \quad \hbox{for } \xi \in S.
       \eq{(3)}%{extchar2}
$$}%

In the context of root-graded Lie algebras, the partial sections $g$
and $g'$ in (c) lead to isotopic Lie algebras, see \ref{isot} and
\ref{arlath}.
\endsubsec

\subsec[ecor]{\bf Corollary.} {\it Let $(S,Y)$ be a pre-reflection
system and let $\frL= (\La_\xi)_{\xi\in S}$ be an extension datum of
type $(S,S',Z)$. Let $R=\scE(S,S',Z)$ be the pre-reflection system
defined in {\rm \ref{extchar}(a)}. \sm

{\rm (a)} $R$ is reduced if and only if, for all\/ $0 \ne \xi \in
S$,
$$
\xi, \; c \xi\in S  \hbox{ for } c\in \KK \setminus \{0, \pm 1\}
\quad\implies \quad \La_{c\xi} \cap c\La_{\xi} = \emptyset.
%\eq{gexdapred}
$$

{\rm (b)} $R$ is symmetric if and only if $S$ is symmetric and $
\La_{-\xi}=- \La_\xi$ for all  $\xi \in S\im$. \sm

{\rm (c)} $R$ is a reflection system if and only if $S$ is a
reflection system.}

\es

\subsec[affdef]{\bf Affine reflection systems.} A pre-reflection
system is called an {\it affine reflection system\/}\index{affine
reflection system}\index{reflection system!affine} if it is an
extension of a locally finite root system.\index{extension!locally
finite root system!} A {\it morphism\/}\index{morphism!of affine
reflection systems}\index{affine reflection system!morphism} between
affine reflection systems is a morphism of the underlying
pre-reflection systems. \sm

The name affine reflection system, as opposed to affine
pre-reflection system, is justified in view of \ref{ecor}(c), since
a locally finite root system is a reflection system and hence so is
any extension of it. \sm

A morphism between affine reflection systems is a morphism of the
underlying reflection systems (\ref{sbc}). Isomorphisms have the
following characterization.
$$\vcenter{
\amen{0.80}{A vector space isomorphism $f\colon\, X \to X'$ is an
isomorphism of the reflection systems $(R,X)$ and $(R',X')$ if and
only if $f(R\im) = R'{}\im$ and $f(R\re)= R'{}\re$.}} \eq{(1)} %{affdef0}
$$

By definition, the {\it nullity\/}\index{nullity!of an affine
reflection system}\index{affine reflection system!nullity} of an
affine reflection system $R$ is the rank of the torsion-free abelian
group $\lan R\im\ran$ generated by $R\im$, i.e., nullity $= \dim_\QQ
(\lan R\im\ran \ot_\ZZ \QQ)$. It is clear from the definitions that
a locally finite root system is an affine reflection system of
nullity $0$ and, conversely, any affine reflection system of nullity
$0$ is a locally finite root system. The affine root systems $R(S,
t(S))$ of \ref{afrs} are affine reflection systems of nullity $1$,
see also \ref{affex}(a). \sm

In the following let $(R,X)$ be an affine reflection system and let
$f\:(R,X)\to (S,Y)$ be an extension where $(S,Y)$ is a locally
finite root system. As mentioned in \ref{ext}, $R$ is then coherent
and integral. Also, non-degeneracy of $S$ implies that $S$ is
unique, up to a unique isomorphism.  We will call $S$ the {\it
quotient root system of $R$\/}\index{quotient root
system}\index{root system!quotient} in this context and refer to $f$
as the {\it canonical projection\/}.\index{canonical
projection}\index{extension!locally finite root system!canonical
projection} \sm

Let $Z = \Ker (f)$. One can show that the extension $f$ has a
partial section $g$ over $S\ind$. Let $\frL = (\La_\xi)_{\xi \in S}$
be the extension datum of type $(S, S\ind, Z)$ associated to $f$ and
$g$ in Th.~\ref{extchar}(b). Thus, {\it up to an isomorphism we may
assume that $R$ is given by\/} (\ref{extchar}.1). The extension
datum $\frL$ appearing there has some special properties besides the
ones mentioned in \ref{ecor}. Namely, for $0\ne \xi, \eta \in S$ and
$w\in W(S)$ we have \begin{eqnarray*} %$$\eqalignno{
  \La_\xi &=& \La_{-\xi} = - \La_\xi
                     %\quad \hbox{for all $\xi \in S$}, \cr
         = \La_{w(\xi)}
          \quad \hbox{for all $w\in W(S)$},  \\%&\eq{affdef1}\cr
    \La_\xi &\supset& \La_\xi - \lan \xi, \eta\ch\ran \La_\eta,
         \quad \hbox{and} \\
  \La_\xi   & \supset& \La_{2\xi}
      \quad \hbox{whenever $2\xi \in S$}.
\end{eqnarray*}%}$$
Moreover, define $\La_{\rm diff} := \bigcup_{0 \ne \xi \in S}
\big(\La_\xi - \La_\xi\big)$.
%\eq{exda9}$$
Then $\ZZ  \La_{\rm diff} = \La_{\rm diff}$ and we have: \sm

(a) $R$ is symmetric if and only if $\La_0=-\La_0$. \sm

(b)  All root strings $\SS(\beta, \al)$, $\beta \in R$, $\al \in
R\re$, are unbroken if and only if $\La_{\rm diff} \subset \La_0$.
\sm

(c) $R$ is tame (cf. \ref{sbc}) if and only if $\La_0 \subset
\La_{\rm diff}$. \sm

(d) $|\SS(\be,\al)|\le 5$ for all $\be \in R$ and $\al\in R\re$. \ms

\noindent
We now describe extension data in the irreducible case in
more detail. \es

\subsec[affexd]{\bf Extension data for irreducible locally finite
root systems.}\index{extension!datum! of locally finite root
system}\index{extension datum!of a locally finite root
system}\index{locally finite root system!extension datum} Let $S$ be
an irreducible locally finite root system. Recall the decomposition
$S= S\sh \cup S\lg \cup S\div$ of \ref{afrs}. Let $\frL =
(\La_\xi)_{\xi \in S}$ be an extension datum of type $(S, S\ind,
Z)$.  One knows that $W(S)$ operates transitively on the roots of
the same length (\cite[5.6]{lfrs}). By \ref{affdef} we can therefore
define $\La\sh$, $\La\lg$ and $\La\div$ by
$$ \La_\al = \begin{cases} \La\sh & \text{for $\al \in S\sh$,} \\
                     \La\lg & \text{for $\al \in S\lg$,} \\
                     \La\div & \text{for $0\ne \al\in S\div$.}
\end{cases}
 \eq{(1)}%{exd1}
$$
By convention, $\La\lg$ and $\La\div$ are only defined if the
corresponding set of roots $S\lg$ and $S\div \setminus \{0\}$ are
not empty (we always have $S\sh\ne \emptyset$). To streamline the
presentation, this will not be specified in the following. It will
(hopefully) always be clear from the context what is meant. \sm

(a) We have already seen in \ref{gexda} that in general $\La\sh$ and
$\La\lg$ are pointed reflection subspaces and that $\La\div$ is a
symmetric reflection subspace of $\La$. In the setting of this
subsection we have in addition the following: \sm

(i)  $\La\sh$ is a subgroup of $Z$, if $S\sh\cup\{0\}$ contains a
subsystem of type $\rmfa_2$, \sm

(ii) $\La\lg$ is a subgroup of $Z$, if $S\lg\cup \{0\}$ contains a
subsystem of type $\rmfa_2$, and \sm

(iii) the following relations hold for $S$ as indicated and $k=k(S)$
defined in\/ \ref{afrs}
$$  \begin{matrix}%$$\eqalign{
  \La\sh +  \La\lg &\subset  \La\sh, \\
  \La\sh + \La\div & \subset \La\sh, \\
  \La\lg + \La\div &\subset \La\lg ,
   \end{matrix}\qquad
 \begin{matrix} %\eqalign{
  \La\lg + k\La\sh &\subset \La\lg, \\
  \La\div + 4 \La\sh &\subset  \La\div, \\
  \La\div + 2\La\lg &\subset \La\div, \end{matrix}
  \qquad
  \begin{matrix}%\eqalign{
    (S\lg &\ne &\emptyset)\\
    (S  &= &\rmbc_1) \\
  (S & = &\rmbc_I , &|I|\ge 2).\end{matrix} %\eq{supplem123}
 $$

{\rm (b)} Conversely, given a subset $\La_0$ with $0\in \La_0$,
pointed reflection subspaces $\La\sh$ and $\La\lg$ (if $S\lg \ne
\emptyset$) and a symmetric reflection subspace $\La\div\subset \La$
(if $\{0\} \subsetneq S\div$) satisfying {\rm (i)} -- {\rm (iii)} of
{\rm (a)}, then (1) defines an extension datum of type
$(S,S\ind,Z)$. Moreover, the following hold where, we recall $k$ is
defined in (\ref{afrs}.1): \sm

\iitem{(i)$'$} $\La\sh$ is a subgroup of $Z$ if $S\ne \rmfa_1$,
$\rmb_I$ or $\rmbc_I$ for any $I$, \sm

\iitem{(ii)$'$} $\La\lg$ is a subgroup of $Z$ if $S=\rmb_I$ or\/
$\rmbc_I$ and $|I|\ge 3$, or if $S=\rmf$ or $\rmg$, \sm

\iitem{(iv)} $k\La\sh \subset \La\lg \subset \La\sh$, $4\La\sh
\subset \La\div\subset \La\sh$ and $\La\div \subset \La\lg \subset
\La\sh$. \sm

\iitem{(v)} The inclusions $\La\sh + \La\div \subset \La\sh$ and
$\La\div + 4 \La\sh \subset \La\div$ hold for all $S=\rmbc_I$. \ms

\noindent One therefore obtains the following description of
$\frL\ti:= (\La\sh, \La\lg, \La\div)$ for the various types of
irreducible root systems. Note that the only condition on $\La_0$ is
$0\in \La_0$ because of (ED2).
\smallskip

(I) $S$ is simply laced, i.e., $S=\rma_I$, $\rmd_I\;  (|I|\ge 4)$,
$\rme_6, \rme_7$ or $\rme_8$: $\frL\ti=(\La\sh)$, where $\La\sh$ is
a pointed reflection subspace for $S=\rmfa_1$ and a subgroup of $Z$
otherwise.  \smallskip

(II) $S=\rmb_I (|I|\ge 2), \rmc_I (|I|\ge 3)$ or  $\rmf$: $\frL\ti=
(\La\sh, \La\lg)$, where $\La\sh$ and $\La\lg$ are pointed
reflection subspaces satisfying $\La\sh + \La\lg \subset \La\sh $
and $\La\lg + 2 \La\sh \subset \La\lg$. Moreover, $\La\sh$ is a
subgroup of $Z$ if $S=\rmc_I$ or $\rmf$, while $\La\lg$ is a
subgroup if $S=\rmb_I$, $|I|\ge 3$ or $\rmf$.
\smallskip

(III)  $S=\rmg$: $\frL\ti=(\La\sh, \La\lg)$, where $\La\sh$ and
$\La\lg$ are subgroups of $Z$ satisfying $\La\sh + \La\lg \subset
\La\sh$ and $\La\lg + 3 \La\sh \subset \La\lg$. \smallskip

(IV) $S=\rmbc_1$: $\frL\ti=(\La\sh, \La\div)$, where $\La\sh$ is a
pointed reflection subspace, $\La\div$ is a symmetric reflection
subspace and $\La\sh + \La\div \subset \La\sh$ and $\La\div + 4
\La\sh \subset \La\div$. \smallskip

(V) $S=\rmbc_I, |I|\ge 2$: $\frL\ti=(\La\sh, \La\lg, \La\div)$ where
$\La\sh$ and $\La\lg$ are pointed reflection subspaces, $\La\div$ is
a symmetric reflection subspace and $\La\sh + \La\lg \subset
\La\sh$, $\La\lg + 2 \La\sh \subset \La\lg$, $\La\lg + \La\div
\subset \La\lg$, $\La\div + 2 \La\lg \subset \La\div$. Moreover, if
$|I|\ge 3$ we require that $\La\lg$ is a subgroup of $Z$.
 \es

\subsec[affo]{\bf Affine forms.} Our definition of affine reflection
systems follows the approach of \cite{brac} where root systems are
defined without reference to a bilinear form. In the literature, it
is customary to define affine root systems and their
generalizations, the extended affine root systems (EARS), in real
vector spaces using positive semidefinite forms. We will therefore
give a characterization of affine reflections systems in terms of
affine invariant forms where, by definition, an {\it affine
form}\index{affine form}\index{pre-reflection system!affine
form}\index{affine reflection system!affine form} of a
pre-reflection system $(R,X)$ over $\KK$ is an invariant form $b$
satisfying $R\im = R\cap \Rad b$, cf. (\ref{ibf}.2). In particular,
affine forms are strictly invariant in the sense of \ref{ibf}. \sm

As an example, the forms used in the theory of EARS are affine forms
in our sense. In particular, the form $\inpr$ in \ref{afrs} is an
affine form of the affine root system $R$.
\endsubsec

\subsec[affcar]{\bf Theorem.} {\it Let $(R,X)$ be a pre-reflection
system. Then $(R,X)$ is an affine reflection system if and and only
if it satisfies the following conditions: \sm

\iitem{\rm (i)} $(R,X)$ is integral,

\iitem{\rm (ii)} $(R,X)$ has an affine form, and

\iitem{\rm (iii)} $\lan R, \al\ch\ran$ is bounded for every $\al\in
R\re$. \smallskip

\noindent In this case: \sm

{\rm (a)}  Let $b$ be an affine form for $(R,X)$ and let $f\colon\,
X \to X/\Rad b$ be the canonical map. Then $(S,Y)= (f(R), X/\Rad b)$
is the quotient root system of $R$ and $f$ its canonical projection.
Moreover, ${\rm Re}(R)$ is connected if and only if $S$ is
irreducible. \sm

{\rm (b)} There exists a unique affine form\/ $\inpr_a$ on $(R,X)$
that is\/ {\rm normalized} in the sense of\/ {\rm \ref{lfrschar}},
i.e., for every connected component $C$ of\/ ${\rm Re}(R)$ we have
$$   2\in \{(\al|\al)_a : 0 \ne \al\in C\} \subset \{2,3,4,6,8\}.
$$
The form\/ $\inpr_a$ satisfies $ \{(\al|\al)_a : 0\ne \al\in C\} \in
\big\{ \{2\}, \{2,4\}, \{2,6\}, \{2,8\}, \{2,4,8\}\big\}$.
%\eq{affcar2}
 Its radical is\/ $\Rad \,\inpr_a=\Ker f$. If\/ $\KK=\RR$ then\/
$\inpr_a$ is positive semidefinite. } \endsubsec

\subsec[semi]{\Corollary}{\it  A pre-reflection system over the
reals is an affine reflection system if and only if it is integral
and has a positive semidefinite affine form. } \es

\subsec[affex]{\bf Special types of affine reflection systems.} As
usual, the {\it rank\/} of a reflection system $(R,X)$ is defined as
$\rank(R,X) = \dim X$. If $R$ has finite rank and $\KK=\RR$, we will
say that $R$ is {\it discrete\/}\index{discrete!reflection
system}\index{affine reflection system!discrete} if $R$ is a
discrete subset of $X$. \sm

(a) Let $(R,X)$ be an affine reflection system over $\KK=\RR$ with
the following properties: $R$ has finite rank, ${\rm Re}(R)$ is
connected, $R=-R$ and  $R$ is discrete. Then $R$ is called \sm

\iitem{(i)} an EARS,\index{EARS|see{extended affine root system}} an
abbreviation for {\it extended affine root system\/},\index{extended
affine root system}\index{root system!extended affine}\index{affine
root system!extended} if $R$ is reduced, tame (see
(\ref{affdef}.c)), and all root strings are unbroken; \sm

\iitem{(ii)} a SEARS,\index{SEARS|see{Saito's extended affine root
system}} an abbreviation of {\it Saito's extended affine root
system\/},\index{Saito's extended affine root system}\index{extended
affine root system!Saito's}\index{affine root system!Saito's
extended} if $R={\rm Re}(R)$. \sm

\noindent It was shown in \cite[Th.~18]{az:ext} that every reduced
SEARS can be uniquely extended to an EARS and, conversely, the
reflective roots of an EARS are the non-zero roots of a SEARS. This
is now immediate. Indeed, by the results mentioned in \ref{affdef},
an affine reflection system is tame and has unbroken root strings if
and only $R\im = \La_{\rm diff}$. The affine root systems $R(S,
t(S))$ of \ref{afrs} with $S$ finite are precisely the EARSs of
nullity $1$ (\cite{abgp}). \sm

(b) In \cite{morita-yoshii}, Morita and Yoshii define a
LEARS,\index{LEARS|see{locally extended affine root system}} an
abbreviation for {\it locally extended affine root
system\/}.\index{locally extended!affine root system}\index{extended
affine root system!locally}\index{affine root system!locally
extended} In our terminology, this is a symmetric affine reflection
system $R$ over $\KK=\RR$ such that $R={\rm Re}(R)$ is connected.
The equivalence of this definition with the one in
\cite{morita-yoshii} follows from \ref{semi}. \sm

(c) In \cite{az:grla}, Azam defines a
GRRS,\index{GRRS|see{generalized reductive root system}} an
abbreviation of a {\it generalized reductive root
system\/}.\index{generalized reductive!root system}\index{root
system!generalized reductive} In our terminology, this is a
symmetric real reduced, discrete affine reflection system $R$ which
has finite rank and unbroken root strings.
\es%\endsubsec

\subsec[peraffrsnotes]{\bf Notes.} With the exceptions mentioned in
the text and below,  all results in this section are proven in
\cite{prs}. \sm

If $S$ is an integral reflection system, our definition of an
extension datum in \ref{gexda} makes sense for any abelian group $Z$
instead of a $\KK$-vector space. This generality is however not
needed for extension data arising from Lie algebras with toral
subalgebras, \ref{inro}, like the generalizations of affine Lie
algebras to be discussed later. But we note that the results on
extension data stated in \ref{affdef} and \ref{affexd} are true for
an abelian group $Z$. In this generality, but for $S$ a finite
irreducible root system, extension data without $\La_0$ were defined
by Yoshii in \cite{y:ext} as ``root systems of type $S$ extended by
$Z$''. The results stated in \ref{affexd} extend
\cite[Th.~2.4]{y:ext} to the setting of locally finite root systems.
\sm

We point out that our definition of an EARS is equivalent to the one
given by Azam, Allison, Berman, Gao and Pianzola in \cite[II,
Def.~2.1]{aabgp}. That our definition of a SEARS is the same as
Saito's definition of an ``extended affine root system'' in
\cite{sa}, follows from \ref{semi}. \sm

Our characterization of extensions in terms of extension data
\ref{extchar} generalizes the Structure Theorem for extended affine
root systems \cite[II, Th.~2.37]{aabgp} as well as
\cite[Prop.~4.2]{morita-yoshii} and  -- modulo the limitations
mentioned above -- the description of root systems extended by an
abelian group in \cite[Th.~3.4]{y:ext}. The isomorphism criterion
(\ref{affdef}.1) is proven in \cite[]{y:LARS} for the case of LEARSs
(\ref{affex}(c)). A proof of the general case will be contained in
\cite{n:lect}. \sm

The extension data arising in the theory of extended affine root
systems are studied in detail in \cite[Ch.~II]{aabgp} and
\cite{Az1}. These references are only a small portion of what is
presently known on extended affine root systems and their Weyl
groups. Many of these results likely have generalizations to the
setting of affine reflection systems.
\endsubsec

\vfe

\section{Graded algebras} \label{pergralg}

In the previous sections \S\ref{perrootsys} and \S\ref{peraffrs} we
reviewed the ``combinatorics'' needed to describe extended affine
Lie algebras and their generalizations in \S\ref{pereala}. In this
and the following section we will introduce the necessary algebraic
concepts. All of them have to do with algebras graded by an
arbitrary abelian group, usually denoted $\La$ and written
additively. (The reader will notice that at least in the first few
subsections $\La$ need not be abelian.) \ms

\subsec[algperse]{\bf Algebras per se.} Throughout, we will consider
algebras $A$ over a unital commutative associative ring of scalars
$k$. Unless explicitly stated otherwise, we do not assume that $A$
belongs to some special variety of algebras, like associative
algebras or Lie algebras. Therefore, an {\it
algebra\/}\index{algebra} or a {\it
$k$-algebra\/}\index{algebra!$k$-algebra}\index{k-algebra@$k$-algebra}
if we want to be more precise, is simply a $k$-module $A$ together
with a $k$-bilinear map $A \times A \to A$, $(a,b) \mapsto ab$, the
{\it product of $A$}. As is customary, if $A$ happens to be a Lie
algebra, its product will be denoted $[a,b]$ or sometimes just
$[a\,b]$ instead of $ab$.\index{$[a\,b]$} \sm

Let $A$ be an algebra. The $k$-module of all its $k$-linear
derivations is always a Lie algebra, denoted
$\Der_kA$.\index{Der@$\Der$ (derivation
algebra}\index{algebra!derivation} A symmetric bilinear form $\inpr
: A \times A \to k$ is called {\it
invariant\/}\index{invariant!bilinear form}\index{algebra!invariant
bilinear form}\index{bilinear form!invariant} if $(ab | c) = (a |
bc)$ holds for all $a,b,c\in A$; it is {\it nondegenerate\/} if $(a
| A) = 0$ implies $a=0$. We put $AA=\Span_k \{ ab : a,b \in A\}$.
The algebra $A$ is {\it perfect}\index{perfect
algebra}\index{algebra!perfect} if $A=AA$, and is {\it
simple\/}\index{simple algebra}\index{algebra!simple} if $AA\ne 0$
and if $\{0\}$ and $A$ are the only ideals of $A$. Any unital
algebra is perfect, and any simple algebra is perfect since $AA$ is
an ideal of $A$.
\endsubsec

\subsec[centroid]{\bf The centroid of an algebra.} The {\it
centroid\/}\index{centroid}\index{algebra!centroid} of an algebra
$A$ is the subalgebra $\cent_k(A)$\index{cent@$\cent(A)$ (centroid
of $A$)} of the endomorphism algebra $\End_k(A)$ consisting of the
$k$-linear endomorphisms of $A$, which commute with the left and
right multiplications of $A$. One always has $k \Id_A \subset
\cent_k(A)$, but in general this is a proper inclusion. If the
canonical map $k \to \cent_k(A)$ is an isomorphism, $A$ is called
{\it central\/},\index{central}\index{algebra!central} and one says
that $A$ is {\it
central-simple\/}\index{central-simple}\index{algebra!central-simple}
if $A$ is central and simple, see \ref{grcesich} for a
characterization of central-simple algebras. \sm

The centroid is always a unital associative algebra, but not
necessarily commutative in general. For example, the centroid of the
null algebra (all products are $0$) is the full endomorphism
algebra. On the other hand, the centroid of a perfect algebra is
always commutative. For example, let $A$ be a unital associative
algebra and denote by $\rmZ(A) = \{ z\in A : za=az \hbox{ for all }
a\in A\}$\index{Z@$\rmZ(.)$ (the centre of an
algebra)}\index{algebra!centre (associative)} its {\it centre}, then
$$
   \rmZ(A) \to \cent_k(A), \quad z \mapsto L_z \qquad
      (\hbox{$A$ associative}) \eq{(1)}%{centroid1}
$$
is an isomorphism of algebras, where $L_z$ denotes the left
multiplication by $z$. However, even if the names are quite similar,
one should not confuse the centroid with the centre of $A$! A
non-zero Lie algebra $L$ may well be {\it
centreless}\index{centreless!Lie algebra} in the sense that
$\rmZ(L)= \{0\}$, but always has $0 \ne k\, \Id_L \in \cent(L)$, see
\ref{fiex}(a) for an example. \sm

Since $\cent_k(A)$ is a subalgebra of $\End_k(A)$ one can consider
$A$ as a left module over $\cent_k(A)$. This change of perspective
is particularly useful if $\cent_k(A)$ is commutative since then $A$
is an algebra over $\cent_k(A)$. In general, an algebra $A$ is
called {\it fgc}\index{fgc} (for $\underline{\rm f}$initely
$\underline{\rm g}$enerated over its $\underline{\rm c}$entroid) if
$A$ is a finitely generated $\cent_k(A)$-module. Fgc algebras are
more tractable than general algebras, and we will characterize
various classes of fgc algebras throughout this paper, see
\ref{grsico}, \ref{perweakv}, \ref{quatoex} and \ref{rocent}.
\es %\endsubsec

\subsec[pergradef]{\bf Graded algebras.} Recall that $\La$ is an
abelian group.  A {\it $\La$-graded algebra}\index{graded
algebra}\index{algebra!graded|see{graded algebra}} is an algebra $A$
together with a family $(A^\la : \la \in \La)$ of submodules $A^\la$
of $A$ satisfying $A = \bigoplus_{\la \in \La} A^\la$ and $A^\la
A^\mu \subset A^{\la + \mu}$ for all $\la,\mu \in \La$. We will
usually indicate this by saying ``Let $A= \bigoplus_\la A^\la$ be a
$\La$-graded algebra ...'' or simply ``Let $A$ be graded
algebra...''. \sm

We will encounter many example of graded algebras. But an immediate
and important example is $k[\La]$, the {\it group algebra of
$\La$\/}\index{group algebra}\index{kl@$k[\La]$ (group algebra of
$\La$)} over $k$. By definition, $k[\La]$ is a free $k$-module with
a $k$-basis in bijection with $\La$, say by $\la \mapsto z^\la$, and
product determined by $z^\la z^\mu = z^{\la + \mu}$ for $\la,\mu \in
\La$. It is graded by $(k[\La])^\la = k\, z^\la$ and is a unital
associative algebra. \ms

We will need some more terminology for a graded algebra $A$. The
submodules $A^\la$ of $A$ are referred to as {\it homogeneous
spaces\/}.\index{homogenous space}\index{graded algebra!homogeneous
space} We will occasionally use subscripts to describe them, in
particular whenever we consider an algebra with two gradings, say
$A= \bigoplus_{\la \in \La} A^\la$ and $A= \bigoplus_{q\in \scq}
A_q$. The two gradings are called {\it compatible\/}\index{graded
algebra!compatible grading}\index{compatible grading} if, putting
$A^\la_q := A^\la \cap A_q$, we have $A^\la = \bigoplus_{q \in \scq}
A^\la_q$ for all $\la\in \La$. We will indicate compatible gradings
by saying that {\it $A$ is $(\scQ,\La)$-graded.\/} \sm

The {\it support\/}\index{support}\index{graded algebra!support} of
a $\La$-graded algebra $A$ is the set $\supp_\La A = \{ \la \in \La
: A^\la\ne 0 \}$. In general $\supp_\La A$ is a proper subset of
$\La$. The subgroup of $\La$ generated by $\supp_\La A$ will be
denoted $\lan \supp_\La A \ran$. One says that {\it $A$ has full
support\/}\index{support!full}\index{graded algebra!full support} if
$\lan \supp_\La A \ran = \La$. Since the grading of $A$ only depends
on $\lan \supp_\La A\ran$, it is of course always possible, and
sometimes even useful, to assume that a graded algebra has full
support. \sm

Two graded algebras $A$ and $A'$ with grading groups $\La$ and
$\La'$ respectively are called {\it
isograded-isomorphic\/}\index{isograded-isomorphic}\index{graded
algebra!isograded-isomorphic} if there exists an isomorphism
$f\colon\, A \to A'$ of the underlying $k$-algebras and a group
isomorphism $\ph\colon\, \lan \supp_\La A \ran \to \lan \supp_{\La'}
A'\ran$ satisfying $f(A^\la ) = A'{}^{\ph(\la)}$ for all $\la \in
\supp_\La A$. Note that $\ph$ is uniquely determined by $f$. A more
restrictive concept is that of a {\it
graded-isomorphism}\index{graded-isomorphic}\index{graded
algebra!graded-isomorphic} between two $\La$-graded algebras $A$ and
$A'$, which by definition is an isomorphism $f\colon\, A \to A'$ of
the ungraded algebras satisfying $f(A^\la) = A'{}^\la$ for all $\la
\in \La$. \sm

A $\La$-graded algebra is called {\it graded-simple\/}\index{graded
algebra!graded-simple}\index{graded-simple} if $AA\ne 0$ and $\{0\}$
and $A$ are the only $\La$-graded ideals of $A$. For example,
$k[\La]$ is graded-simple if and only if $k$ is a field. But of
course in this case $k[\La]$ may have many non-trivial ungraded
ideals. \sm

A bilinear form $\inpr \:A \times A \to k$ on $A$ is {\it
$\La$-graded\/}\index{bilinear form!graded} if $(A^\la | A^\mu) = 0$
for $\la + \mu \ne 0$. For example, the Killing form of a
finite-dimensional $\La$-graded Lie algebra is $\La$-graded. The
radical of an invariant $\La$-graded symmetric bilinear form $\inpr$
is a graded ideal of $A$, namely $\Rad \inpr =\{ a\in A :
(a|A)=0\}$. A $\La$-graded bilinear form $\inpr$ is nondegenerate,
i.e. $\Rad \inpr =\{0\}$, if and only if for all $\la \in \La$ the
restriction of $\inpr$ to $A^\la \times A^{-\la}$ is a nondegenerate
pairing.
 \sm

We denote by $\GIF(A,\La)$\index{GIF@$\GIF$ (graded invariant
bilinear forms)} the $k$-module of all $\La$-graded invariant
symmetric bilinear forms on $A$. If $L$ is a perfect $\La$-graded
Lie algebra, then $\Rad \inpr \supset \rmZ(L) = \{z\in L : [z, \,L]
=0\}$,\index{Z@$\rmZ(.)$ (the centre of an algebra)} the {\it
centre\/}\index{algebra!centre (Lie)} of $L$. Since $\rmZ(L)$ is
also $\La$-graded, the quotient algebra inherits a canonical
$\La$-grading, and every $\inpr \in \GIF(L,\La)$ induces an
invariant symmetric bilinear form on $L/\rmZ(L)$ by $(\bar x | \bar
y) = (x|y)$, where $x \mapsto \bar x$ denotes the canonical map.
This gives rise to an isomorphism
$$ \GIF(L,\La) \cong \GIF(L\big/ \rmZ(L), \La)
 \qquad \hbox{(Lie algebras).} \eq{(1)}%{pergradef0}
$$
If $A$ is a unital algebra, say with identity element $1$, then
necessarily $1\in A^0$. If in addition $A$ is associative, then
$$\GIF(A,\La) \cong ( A^0\big/ [A,A]^0)^*
 \qquad \hbox{(associative algebras)}
 \eq{(2)}%{pergradef1}
$$ by assigning to any linear form $\ph$ of $A^0\big/
[A,A]^0$ the bilinear form $(a\mid b)_\ph = \ph( (ab)^0 + [A,A]^0)$.
Here $[A,A] = \Span_k \{ ab-ba : a,b\in A\}$,\index{$[A,A]$} and
$(ab)^0$ and $[A,A]^0$ denote the $0$-component of $ab$ and $[A,A]$.
For example, the group algebra $k[\La]$ has up to scalars just one
invariant symmetric bilinear form, which assigns to $(a,b)$ the
$z^0$-coefficient of $ab$ with respect to the canonical $k$-basis.
\sm

An endomorphism $f$ of the underlying $k$-module of $A$ is said to
have {\it degree $\la$\/}\index{endomorphism!of degree $\la$} if
$f(A^\mu) \subset A^{\la + \mu}$ for all $\mu \in \La$. The
submodules $(\End_kA)^\la$, consisting of all endomorphisms of
degree $\la$, are the homogeneous spaces of the $\La$-graded
subalgebra
$$
   \grEnd_k (A) = \ts \bigoplus_{\la\in \La}\, (\End_k A)^\la
     \index{gr@$\grEnd$ (graded endomorphism algebra)}
$$
of the associative algebra $\End_k (A)$. We put
 \begin{eqnarray*}%$$\eqalign{
   \grDer_k (A) &=&
    \grEnd_k(A) \cap \Der_k(A)
      = \ts \bigoplus_{\la\in \La} (\Der_k A)^\la,
           \quad\hbox{and} \index{gr@$\grDer$ (graded derivation algebra)}\\
  \grcent_k (A) &=& \grEnd_k(A) \cap \cent_k(A)
      = \ts \bigoplus_{\la \in \La} \, (\cent_k A)^\la.
        \index{gr@$\grcent$ (graded centroid)}
 \end{eqnarray*}%$$
Then $(\Der_k A)^\la$\index{Dera@$(\Der A)^\la$ (derivations of
degree $\la$)} and $(\cent_k A)^\la$\index{cente@$(\cent A)^\la$ (
centroidal transformations of degree $\la$)} are, respectively, the
derivations and centroidal transformations of $A$ of degree $\la$.
In general, $\grDer_k(A)$ or $\grcent_k(A)$ are proper subalgebras
of $\Der_k(A)$ and $\cent_k(A)$. However, $\grDer_k(A) = \Der_k(A)$
if $A$ is a finitely generated algebra (\cite[Prop.~1]{Fa1}).
Similarly, $\grcent_k(A) = \cent_k(A)$ if $A$ is finitely generated
as ideal. For example, any unital algebra is finitely generated as
ideal (namely by $1\in A$). We will discuss some examples of $({\rm
gr})\Der_k A$ and $({\rm gr})\cent_k A$ later, see \ref{fiex},
\ref{quatoex}, \ref{propa} and \ref{derit}. \sm

We call $A$ {\it graded-central\/}\index{graded
algebra!graded-central}\index{graded-central algebra} if the
canonical map $k \to (\cent_k A)^0$ is an isomorphism, and {\it
graded-central-simple\/}\index{graded
algebra!graded-central-simple}\index{graded-central-simple algebra}
if $A$ is graded-central and graded-simple. It is well-known that an
(ungraded) algebra $A$ over a field $k$ is central-simple if and
only if $A\ot_k K$ is simple for every extension field $K$ of $k$.
The following theorem extends this to the case of graded algebras.
\endsubsec

\subsec[grcesich]{\Theorem}{\it Let $A$ be a graded-simple
$k$-algebra, where $k$ is a field. Then the following are
equivalent: %\sm
\begin{description}
\item{\rm (i)} $A$ is graded-central-simple, \sm

\item{\rm (ii)} for any extension field $K$ of $k$ the base field
extension $A\ot_k K$ is a graded-central-simple $K$-algebra. \sm

\item{\rm (iii)} for any extension field $K$ of $k$ the algebra
$A\ot_k K$ is a graded-simple $K$-algebra. \sm

\item{\rm (iii)} $A\ot_k (\cent_k A)^0$ is graded-simple.
\end{description}  }
\ms

In this theorem, ``base field extension''\index{algebra!base field
extension} means that we consider $A\ot_k K$ as the $K$-algebra with
product $(a_1 \ot x_1)\, (a_2 \ot x_2) = (a_1a_2) \ot (x_1 x_1)$ for
$a_i \in A$ and $x_i \in K$ and the $\La$-grading determined by
$(A\ot_k K)^\la = A^\la \ot_k K$. \ms

We will describe the centroids of graded-simple algebras in
\ref{grsico}, after we have introduced special classes of
associative graded algebras.
%\endsubsec%
 \es

\subsec[assogr]{\bf Associative (pre)division-graded algebras and
tori.} A $\La$-graded unital associative algebra $A= \bigoplus_{\la
\in \La} A^\la$ is called \sm

\begin{description}
\item{$\bull$} {\it predivision-graded\/},\index{predivision-graded associative algebra}
\index{graded algebra!predivision-graded (associative)} if every
non-zero homogeneous space contains an invertible element; \sm

\item{$\bull$} {\it division-graded\/}\index{division-graded associative algebra}
\index{graded algebra!division-graded (associative)} if every
non-zero homogenous element is invertible; \sm

\item{$\bull$} an {\it associative torus\/}\index{torus!associative}
\index{graded algebra!torus} if $A$ is
predivision-graded, $k$ is a field and $\dim_k A^\la \le 1$.
\end{description}\sm

For example, any group algebra $k[\La]$ is predivision-graded.
Moreover, $k[\La]$ is division-graded if and only if $k$ is a field,
and in this case $k[\La]$ is a torus. \sm

Group algebras cover only a small part of all possible
predivision-graded algebras. If $A$ is such an algebra, it is easily
seen that the support of $A$ is a subgroup of the grading group
$\La$. To simplify the notation,
$$
\hbox{\it
   we assume $\supp_\La A = \La$ in this subsection.}
$$
Let $(u_\la : \la \in \La)$ be a family of invertible elements
$u_\la \in A^\la$ and put $B=A^0$. Then $(u_\la )$ is a $B$-basis of
the left $B$-module $A$, and the product of $A$ is uniquely
determined by the two rules
$$ u_\la \, u_\mu = \ta(\la, \mu ) \, u_{\la + \mu }
 \quad \hbox{and} \quad
 u_\la \, b = b^{\si(\la)}\, u_\la \quad(b\in B)  \eq{(1)}%{assogr1}
$$
where $\ta \colon\, \La \times \La \to B\ti$ and $\si \colon\, \La
\to \Aut_k (B)$ are functions. Associativity of $A$ leads to the two
identities
$$
 \begin{matrix}
      \ta(\la,\mu)^{\si(\nu)} \, \ta(\nu,\la+ \mu)
  &= & \ta(\nu,\la)\,  \ta(\nu\la,\mu),
    & \hbox{ and} \\
  b^{\si(\la)} )^{\si(\nu)}
   \,\ta(\nu,\la) &= & \ta(\nu,\la) \, b^{\si(\nu +\la)}
  \end{matrix} \eq{(2)}
$$
for $\nu,\la,\mu\in \La$ and $b\in B$. Conversely, given any unital
associative $k$-algebra $B$ and functions $\ta$, $\si$ as above
satisfying (2), one can define a $\La$-graded $k$-algebra by (1). It
turns out to be an associative  predivision-graded algebra. Algebras
arising in this way are called {\it crossed product
algebras\/}\index{crossed product algebra}\index{graded
algebra!crossed product algebra} and denoted $(B, (u_\la),
\ta,\si)$. To summarize,
$$\vcenter{\amen{0.8}{any associative predivision-graded algebra $A$ with support
$\La$ is graded-isomorphic to a crossed product algebra $(B,
(u_\la), \ta,\si)$. \\ $A$ is division-graded iff $B$ is a division
algebra, and $A$ is a torus iff $B=k$ is a field. }} \eq{(3)}%{assogr4}
$$

The description of predivision-graded algebras as crossed product
algebras may not be very illuminating, except that it (hopefully)
demonstrates how general this class is. The reader may be more
comfortable with the subclass of {\it twisted group
algebras\/}\index{group algebra!twisted}\index{graded algebra!group
algebra!twisted} defined by the condition that $\si$ is trivial,
i.e., $\si(\la) = \Id_B$ for all $\la \in \La$, and denoted
$B^t[\La]$\index{kla@$k^t[\La]$ (twisted group algebra of $\La$)}.
For example, any commutative crossed product algebra is a twisted
group algebra. Or, up to graded isomorphism, a $\La$-torus with full
support is the same as a twisted group algebra $k^t[\La]$. In the
torus  case, (2) says that $\ta$ is a {\it
$2$-cocycle\/}\index{co@$2$-cocyle!group} of the group $\La$ with
values in $k\ti$. We will describe the special case of associative
$\ZZ^n$-tori in \ref{quatoex}. \sm

Twisted group algebras arise naturally as centroids of graded-simple
algebras.%\endsubsec%
\es

\subsec[grsico]{\Proposition}{\it Let $A= \bigoplus_{\la \in \La}
A^\la$ be a graded-simple $k$-algebra. \sm

{\rm (a)} Then $\cent_k(A)$ is a commutative associative
division-graded $k$-algebra,
$$
   \cent_k (A) = \grcent_k(A) = \ts \bigoplus_{\ga \in \Ga}\,
      (\cent_k A)^\ga \quad\hbox{where} \eq{(1)}%{grsico0}
$$
$$\Ga = \supp_\La \cent_k A\index{G@$\Ga$ (centroidal grading
group)}$$ is a subgroup of $\La$, called the {\rm central grading
group of $A$\/}.\index{centroid!centroidal grading
group}\index{centroidal grading group} Hence $\cent_k(A)$ is a
twisted group algebra over the field $K=(\cent_k A)^0$. \sm

{\rm (b)} $A$ is a graded-central-simple $K$-algebra. \sm

{\rm (c)} $A$ is a free $\cent_k(A)$-module with $\Ga + \supp_\La A
\subset \supp_\La A$. Moreover, $A$ has finite rank as
$\cent_k(A)$-module, i.e. $A$ is fgc, if and only if
$$
   [\supp_\La A : \Ga] < \infty \quad \hbox{and} \quad
    \dim_K (A^\la) < \infty \hbox{ for all } \la \in \La.
 \eq{(2)}%{grsico1}
$$}%

The interest in algebras which are graded-central-simple and also
fgc comes from the following important realization theorem.
\endsubsec

\subsec[perweakv]{\Theorem}(\cite[Cor.~8.3.5]{abfp}) {\it Let $k$ be
an algebraically closed field of characteristic\/ $0$ and let $\La$
be free abelian of finite rank. Then the following are equivalent
for a $\La$-graded $k$-algebra $B$: \sm
\begin{description}
\item{\rm (i)} $B$ is graded-central-simple, fgc and $[\La
:\Ga]<\infty$ where $\Ga$ is the central grading group of $B$; \sm

\item{\rm (ii)} $B$ is isograded-isomorphic to a multiloop algebra
$\scM_\bfm(A,\boldsi)$ for some finite-dimen\-sio\-nal simple
$k$-algebra $A$, described below. \end{description}}\ms

The multiloop algebras referred to in (ii) provide an interesting
class of examples of graded algebras. They are defined as follows
($A$ can be arbitrary, but we assume that $k$ is a field with enough
roots of unity). We let $\bfm=(m_1, \ldots , m_n)\in \NN^n_+$ be an
$n$-tuple of positive integers and let $ \pi : \ZZ^n \to \ZZ/(m_1)
\,\oplus \, \cdots \, \oplus \, \ZZ/(m_n) =:\Xi $ be the canonical
map. Furthermore,  $\boldsi= (\si_1, \ldots, \si_n)$ is a family of
$n$ pairwise commuting automorphisms of the (ungraded) algebra $A$
such that $\si_i^{m_i} = \Id$ for $1\le i\le n$. For $\la =(\la_1,
\ldots, \la_n)\in \ZZ^n$ we define
$$   A^{\pi(\la)} = \{ a \in A : \si_i(a) = \zeta_{m_i}^{\la_i} a \hbox{
for } 1 \le i \le n \}.
$$
where for each $l\in \{m_1, \ldots, m_n\}$ we have chosen a
primitive $l^{\rm th}$-root of unity $\ze_l \in k$. Then $A=\ts
\bigoplus_{\xi \in \Xi} \, A^\xi$ is a $\Xi$-grading of $A$.
Denoting by $k[z_1^{\pm 1}, \ldots, z_n^{\pm 1}]$ the Laurent
polynomial algebra over $k$ in the variables $z_1, \ldots, z_n$,
$$ \scM_\bfm (A, \boldsi) = \ts \bigoplus_{(\la_1,
\ldots, \la_n) \in \ZZ^n} \, A^{ \pi(\la_1, \ldots , \la_n )} \ot k
\, z_1^{\la_1} \cdots z_n^{\la_n}
$$
is a $\ZZ^n$-graded algebra, called the {\it multiloop
algebra}\index{multiloop algebra}\index{graded
algebra!multiloop}\index{loop algebra!multiloop} of $\boldsi$ based
on $A$ and relative to $\bfm$. \ms

Thus, the theorem in particular says that $\scM_\bfm(A, \boldsi)$ is
graded-central-simple. In fact, the centroid can be determined
precisely. \endsubsec%\es

\subsec[]{\Proposition} (\cite[Cor.~6.6]{abp2.5}) {\it The centroid
of a multiloop algebra $\scM_\bfm(A, \boldsi)$ based on a
finite-dimensional central-simple algebra $A$ is graded-isomorphic
to the Laurent polynomial ring $k[z_1^{\pm m_1}, \, \ldots \, ,
z_n^{\pm m_n}]$. } \es

\subsec[centd]{\bf Centroidal derivations.} There is a close
connection between the centroid and the derivations of an algebra
$A$: $[d,\chi] \in \cent_k(A)$ for $d\in \Der_k(A)$ and $\chi \in
\cent_k(A)$. In other words, $\cent_k(A)$ is a submodule of the
$\Der_k(A)$-module $\End_k(A)$. In the graded case one can define
special derivations using the centroid. \sm

Let $A= \bigoplus_{\la \in \La} A^\la$ be a $\La$-graded algebra. We
abbreviate $\rmC= \cent_k(A)$, $\grrmC= \grcent_k(A)$ and $\rmC^\la
= (\cent_k A)^\la$, thus $\grrmC = \bigoplus_\la C^\la$. We denote
by $\Hom(\La, \rmC) \cong \Hom_{\rmC}(\La\ot_\ZZ \rmC, \,\rmC)$ the
(left) $\rmC$-module of all abelian group homomorphisms
$\vth\colon\, \La \to (\rmC,+)$. Every $\vth \in \Hom(\La, \rmC)$
gives rise to a so-called {\it centroidal
derivation\/},\index{centroidal
derivation}\index{derivation!centroidal derivation} defined by
$$
  \pa_\vth (a^\la) = \vth(\la)\, (a^\la) \qquad \hbox{for }
         a^\la \in A^\la. \eq{(1)}%{centd0}
$$
We put
$$
   \CDer (A,\La) = \{ \pa_\vth : \vth \in \Hom(\La, \rmC) \},
\index{CDer@$\CDer$ (centroidal derivations)}
$$ sometimes also denoted $\CDer(A)$ if $\La$ is
clear from the context. We emphasize that $\CDer(A,\La)$ not only
depends on $C$ but also on $\La$. For example, any algebra $A$ can
be graded by $\La=\{0\}$ in which case $\CDer (A,\La) = \{0\}$. \sm

We will also consider the $\La$-graded submodule $\grCDer (A,\La)$
of $\CDer A$, \begin{eqnarray*}%$$\eqalign{
   \grCDer (A,\La) &=&
       \ts \bigoplus_{\la \in \La}\, (\CDer A)^\la
 \qquad \hbox{where} \index{gr@$\grCDer$ (graded centroidal derivations)}
 \\
   (\CDer A)^\la &=& \CDer (A) \cap (\End_k A)^\la =
       \{\, \pa_\vth : \vth \in \Hom (\La, \rmC^\la) \}.
 \end{eqnarray*}% }$$
\ms

{\it Suppose now that $\cent_k(A)$ is commutative\/}, e.g. that $A$
is perfect. Then any $\pa_\vth$ is $\rmC^0$-linear and $\CDer(A)$ is
a $\rmC^0$-subalgebra of the Lie algebra $\Der_{{\rmC}^0} (A)$ of
all ${\rmC}^0$-linear derivations of $A$. Indeed, $[\pa_\vth,\,
\pa_\psi]= \pa_{\vth * \psi - \psi * \vth}$ for $\vth * \psi \in
\Hom(\La, \rmC)$ defined by
$$
    (\vth * \psi )\, (\la) = \tsum_{\mu \in \La}\,
                          \vth(\mu)\, \psi(\la)^\mu \eq{(2)}%{centd1}
$$
where $\psi(\la)^\mu$ denotes the $\mu$-component of $\psi(\la) \in
\rmC$. The sum in (2) converges in the finite topology. In
particular, for $\vth \in \Hom(\La, \rmC^\la)$ and $\psi \in
\Hom(\La, \rmC^\mu)$ one has
$$    [\pa_\vth,\, \pa_\psi] = \vth(\mu)\pa_\psi - \psi(\la) \pa_\vth
        \in (\CDer A)^{\la + \mu}, \eq{(3)}%{centd2}
$$
in particular $(\CDer A)^0$ is abelian.  \ms

{\it Suppose $A$ is perfect and $\inpr$ is a $\La$-graded invariant
nondegenerate symmetric bilinear form on $A$.\/} We denote by
$$ \grSCDer (A, \La) = \ts \bigoplus_{\la \in \La}\,
          (\SCDer A)^\la \index{gr@$\grSCDer$ (graded skew
centroidal derivations)}
$$
the graded skew centroidal derivations of $A$, i.e. the graded
subalgebra of $\grCDer (A,\La)$ consisting of those graded
centroidal derivations of $A$ which are skew-symmetric with respect
to $\inpr$. Since $A$ is perfect, any centroidal transformation is
symmetric with respect to $\inpr$. This implies that
$$
   (\SCDer A)^\la = \{ \pa_\vth \in (\CDer A)^\la: \vth(\la) = 0 \},
$$
in particular $(\SCDer A)^0 = (\CDer A)^0$. Because of (3) we have
$$
   [(\SCDer A)^\la, \, (\SCDer A)^{-\la}] = 0 %\eq{centd3}
$$
for any $\la \in \La$. Hence $\grSCDer (A,\La)$ is the semidirect
product of the ideal spanned by the non-zero homogeneous
spaces and the abelian subalgebra $(\CDer A)^0$.  \es%\endsubsec

\subsec[degd]{\bf Degree derivations.} Let again $A= \bigoplus_\la
A^\la$ be a $\La$-graded $k$-algebra. Recall that $k \Id_A \subset
(\cent_k A)^0$. We can therefore consider the submodule
$$
   \scD= \scD(A,\La) = \{\pa_\vth: \vth \in \Hom(\La, \, k\Id_A)\}
        \subset (\CDer A)^0 \index{D@$\scD$ (degree derivations)}
$$
of so-called {\it degree derivations of $(A,\La)$\/}.\index{degree
derivation}\index{derivation!degree derivation} For example, if
$\La=\ZZ^n$ we have $\scD= \Span_k \{ \pa_i : 1\le i \le n\}$ where
$\pa_i (a^\la) = \la_i a^\la$ for $\la = (\la_1, \ldots, \la_n) \in
\ZZ^n$. \sm

Observe that $[\pa_\vth ,\, \pa_\psi] = \vth(\mu) \, \pa_\psi$ for
$\pa_\vth \in \scD$ and $\psi \in \Hom(\La, \cent_k(A)^\mu)$. So the
action of $\scD$ on $\grCDer A$ is diagonalizable. In order to
actually get a toral action,
$$
\amen{0.55}{we will assume in this subsection and in\/ {\rm
\ref{degl}} that $k$ is a field and that $A$ has full support.}$$ In
this case, we can identify $\Hom(\La, k \Id_A) \equiv
\Hom_\ZZ(\La,k)$. Moreover, the canonical map $\Hom(\La,k) \to \scD=
\{ \pa_\vth : \vth \in \Hom(\La, k)\}$ is an isomorphism. We
therefore have an evaluation map $\ev\colon\, \La \to \scD^*$,
defined by
$$ \ev_\la (\pa_\vth) = \vth(\la). \index{ev@$\ev$ (evaluation map)}
$$
Observe that (\ref{centd}.1) says
$$
   A^\la \subset \{ a\in A : d(a) = \ev_\la(d) \, a
           \hbox{ for all } d\in \scD\}. \eq{(1)}%{degd1}
$$
The following lemma specifies a condition under which (1) becomes an
equality. \endsubsec

\subsec[degl]{\Lemma}{\it As in\/ {\rm \ref{degd}} we suppose that
$k$ is a field and $A$ is a $\La$-graded $k$-algebra with full
support. We let $T \subset \scD$ be a subspace satisfying the
condition that
$$
  \hbox{the restricted evaluation map $\La \to T^*$,
      $\la \mapsto \ev_\la | T$, is injective.} \eq{(1)}%{degl1}
$$
Then $A^\la = \{ a \in A : d(a) = \ev_\la(t) \hbox{ for all }t\in
T\}$, in particular equality holds in\/ {\rm (\ref{degd}.1)}.
Moreover, if\/ $\cent_k(A)$ is commutative, $T$ is a toral
subalgebra of $\grCDer_k(A)$, and if
 $A$ is perfect and has a $\La$-graded invariant nondegenerate
symmetric bilinear form, then $T$ is also a toral subalgebra of\/
$\grSCDer_k(A)$.} \ms

We mention that (1) always holds for $T=\scD$, $\La$ torsion-free
and $k$ a field of characteristic $0$. \es

\subsec[fiex]{\bf Examples.} (a) Let $\frg$ be a finite-dimensional
split simple Lie algebra over a field $\KK$ of characteristic $0$,
e.g. a finite-dimensional simple Lie algebra over an algebraically
closed $\KK$. Hence $\frg$ has a root space decomposition $\frg =
\bigoplus_{\al \in R} \frg_\al$ with respect to some splitting
Cartan subalgebra $\frh = \frg_0$ of $\frg$, where $R\subset \frh^*$
is the root system of $(\frg, \frh)$. We consider this decomposition
as a $\scq(R)$-grading. It is well-known that $\frg$ is central
(hence central-simple): Any centroidal transformation $\chi$ leaves
the root space decomposition of $\frg$ invariant and, by simplicity,
is uniquely determined by $\chi|\frg_\al$ for some $0\ne \al \in R$.
\sm

 For any $\vth \in \Hom
(\scq(R), \KK) \cong \Hom_\KK(\frh^*, \KK)$ there exists a unique
$h_\vth \in \frh$ such that $\al(h_\vth) = \vth(\al)$. Therefore
$\pa_\vth = \ad h_\vth$ and $$\CDer(\frg, \scq(R)) = \grCDer(\frg,
\scq(R)) = \grSCDer(\frg, \scq(R)) = \ad \frh.$$ \ms

(b) Let $\frg$ be as in (a) and let $C$ be a unital commutative
associative $\La$-graded $\KK$-algebra. We view $L=\frg \ot C$ as a
$\La$-graded Lie algebra over $\KK$ with product $[x\ot c, x'\ot c']
= [x, x'] \ot cc'$ for $x,x'\in \frg$ and $c,c'\in C$ and with
homogeneous spaces $L^\la = \frg \ot C^\la$ (here and in the
following all tensor products are over $\KK$). The centroid of $L$
is determined in \cite[Lemma~2.3]{abp2.5} or \cite[Cor.~2.23]{BN}:
$$
   \cent_\KK(L) = \KK\, \Id_\frg \ot C. \eq{(1)}%fiex2}
$$
(We note in passing that (2) also holds for the infinite-dimensional
versions of $\frg$ defined in \ref{afrreal}, as well as for the
version of $\frg$ over rings which we will introduce in
\ref{predivex}.) \ms

By \cite[Th.~7.1]{block} or \cite[Th.~1.1]{bemo}, the derivation
algebra of $L$ is
$$
 \Der_\KK (L) = \big( \Der (\frg) \ot C\big) \, \oplus \,
                 \big( \KK \,\Id_\frg \ot \Der_\KK (C) \big)
 = \IDer (L) \oplus \big( \KK \, \Id_\frg \ot \Der_\KK (C) \big)
 \eq{(2)}%fiex1}
$$
where $\IDer L$\index{I@$\IDer$ (inner derivations)} is the ideal of
inner derivations of $L$. It is then immediate that $\CDer (L)$ and
$\grCDer(L)$ are given by the formulas
$$
  {\rm (gr)CDer} (L, \La) = \KK \, \Id_\frg \ot {\rm (gr)CDer}(C,\La).
$$
where here and in the following (gr) indicates that the formula is
true for the graded as well as the ungraded case. \sm

Let $\ka$ be the Killing form of $\frg$ and suppose that $\eps$ is a
$\La$-graded invariant nondegenerate symmetric bilinear form on $C$.
The bilinear form $\ka\ot \eps$ of $L$, defined by $(\ka\ot
\eps)(x\ot c, x'\ot c') = \ka(x,x') \, \eps(c,c')$, is $\La$-graded,
invariant, nondegenerate and symmetric (if $\frg$ is simple, any
such form arises in this way by \cite[Th.~4.2]{be:derinv}). With
respect to $\ka \ot \eps$, the skew-symmetric derivations and
skew-symmetric centroidal derivations are
 \begin{eqnarray*} %$$\eqalign{
  {\rm (gr)SDer}_\KK (L) &=&
    \IDer L
         \oplus \big( \KK\, \Id_\frg \ot
            {\rm (gr)SDer}_\KK (C)\big), \\
  {\rm (gr)SCDer}_\KK (L)
       &=&   \phantom{IDer L \oplus} \, \phantom{+}
         \KK\,  \Id_\frg \ot  {\rm (gr)SCDer}_\KK (C).
 \end{eqnarray*} %}$$
\sm

In particular, let $C= \KK[t^{\pm 1}]$ be the $\KK$-algebra of
Laurent polynomials over $\KK$, which we view as $\ZZ$-graded
algebra with $C^i = \KK t^i$. Its derivation algebra is (isomorphic
to) the Witt algebra,
$$
   \Der_\KK \big( \KK[t^{\pm 1}]\big) =
            \ts \bigoplus_{i \in \ZZ}\, \KK d^{(i)},
$$
where $d^{(i)}$ acts on $\KK[t^{\pm 1}]$ by $d^{(i)}(t^j) = j
t^{i+j}$. Hence $d^{(i)}$ is the centroidal derivation given by the
homomorphism $\vth^{(i)} \colon\, \ZZ \to C^i$, $\vth^{(i)}(j) = j
t^i$. It follows that
$$
 \Der \big( \KK[t^{\pm 1}]\big) = \grDer\big( \KK[t^{\pm 1}]\big)
    = \CDer \big( \KK[t^{\pm 1}]\big)
  = \grCDer \big( \KK[t^{\pm 1}]\big).
$$
Up to scalars, the algebra $\KK[t^{\pm 1}]$ has a unique
$\ZZ$-graded invariant nondegenerate symmetric bilinear form $\eps$
given by $\eps(t^i,t^j) = \de_{i+j,0}$. With respect to this form $
\SCDer\big( \KK[t^{\pm 1}]\big) = \KK d^{(0)}$ and hence the skew
derivations of the untwisted loop algebra $L = \frg \ot \KK[t^{\pm
1}]$ are
$$
   \SCDer (\frg \ot \KK[t^{\pm 1}]) = \KK d^{(0)}
$$
It is important for later that $\SCDer (L)$ are precisely the
derivations needed to construct the untwisted affine Lie algebra in
terms of $L$. \ms

(c) Let $E$ be an affine Lie algebra (\ref{afrreal}). Then (b)
applies to $$L= [E,\, E]\big/ \rmZ([E,\,E]) \cong \frg \ot
\KK[t^{\pm 1}]$$ (the centreless core when $E$ is considered as an
extended affine Lie algebra, see \ref{ealadef}). In particular, $L$
has an infinite-dimensional centroid. But the centroid of the core
$K:= [E,\,E]$ is $1$-dimensional: $\cent(K) = \KK \, \Id_K$, while
the centroid of $E$ is $2$-dimensional, namely isomorphic to $\KK \,
\Id_E \oplus \Hom_\KK( \KK c, \KK  d)$ where $\KK c$ is the centre
of the derived algebra $K$
and $E = K\rtimes \KK d$ (\cite[Cor.~3.5]{BN}). \es% \endsubsec

\subsec[quatoex]{\bf Quantum tori.} In short, a {\it quantum
torus\/}\index{quantum torus}\index{torus!quantum}\index{graded
algebra!quantum torus} is an associative $\ZZ^n$-torus. As in
\ref{assogr} we may assume that a quantum torus has full support. In
this case it can be described as follows. \sm

Let $k=F$ be a field and let $q=(q_{ij}) \in \Mat_n(F)$ be a {\it
quantum matrix\/},\index{quantum matrix} i.e., $q_{ij}\,q_{ji}=
1=q_{ii} $ for all $1\le i,j\le n$. The {\it quantum torus
associated to $q$\/} is the associative $F$-algebra
$F_q$\index{fq@$F_q$ (the quantum torus associated to $q$)}
presented by generators $t_i, {t_i}^{-1}$ and relations
$$
 t_i {t_i}^{-1} = 1 = {t_i}^{-1} t_i
\quad\hbox{and} \quad t_i t_j = q_{ij} t_j t_i \quad\hbox{for $1\le
i,j\le n$.}
$$  Observe that if all $q_{ij}=1$ then $F_q = F[t_1^{\pm
1}, \cdots , t_n^{\pm 1}]$ is the Laurent polynomial ring in $n$
variables. Hence, for a general $q$, a quantum torus is a
non-commutative version (a quantization) of the Laurent polynomial
algebra $F[t_1^{\pm 1}, \cdots , t_n^{\pm 1}]$, which is the
coordinate ring of an $n$-dimensional algebraic torus. This explains
the name ``quantum torus'' for $F_q$.  \sm

Let $A= F_q$ be a quantum torus. It is immediate that
 $$
      F_q = \ts \bigoplus_{\la \in \ZZ^n} Ft^\la, \qquad
            \hbox{for}\quad  t^\la = t_1^{\la_1}\cdots t_n^{\la_n}.
$$
The multiplication of $F_q$ is determined by the $2$-cocycle
$\ta\colon\, \ZZ^n \times \ZZ^n \to F\ti$ where
$$
   t^\la t^\mu = \ta(\la,\mu) \, t^{\la + \mu}, \qquad
    \ta(\la,\mu) = \ts \prod_{1 \le j < i \le n}\,
q_{ij}^{\la_i\mu_j}. \eq{(1)}%quantoex1}
$$
It then clear that $F_q$ is an associative $\ZZ^n$-torus with
$(F_q)^\la = Ft^\la$. Conversely, every associative $\ZZ^n$-torus is
graded-isomorphic to some quantum torus. \sm

The structure of quantum tori is well understood. We mention some
facts which are relevant for this paper. Let $A=F_q$, $q\in
\Mat_n(F)$, be a quantum torus. \sm

(a) (\cite[Lemma~3.1]{OP}, \cite[Th.~2.4]{neeb:quantumtori}) Let
$\tl q \in \Mat_n(F)$ be another quantum matrix and denote by $\tl
\ta$ the $2$-cocycle associated in (1) to the quantum torus $F_{\tl
q}= \bigoplus_{\la \in \ZZ^n} F \tl t ^\la$. Then the following are
equivalent: \sm
\begin{description}

\item{(i)} $F_q$ and $F_{\tl q}$ are isomorphic as algebras;

\item{(ii)} $F_q$ and $F_{\tl q}$ are isograded-isomorphic as
$\ZZ^n$-graded algebras;

\item{(iii)} there exists $M\in \GL_n(\ZZ)$ such that the
cohomology classes of the $2$-cocycles $\ta$ and $\tl \ta \circ
(M\times M)$ agree, i.e., $\ta(\la, \mu)\big/ \tl \ta(M\la,M\mu) =
v_\la \, v_\mu \, v_{\la + \mu}^{-1}$ for some function $v\colon\,
\ZZ^n \to k\ti$.  \end{description} \sm

\noindent In this case, the map $f\colon\, F_q \to F_{\tl q}$, given
by $f(t^\la) \mapsto v_\la\, {\tl t}^{\,M\la}$, is an
isograded-isomorphism, and every isomorphism arises in this way for
suitable $(v_\la)$ and $M$. \sm

(b) (\cite[Lemma~1.1]{OP}, \cite[Prop.~2.44]{bgk}) The centre of
$F_q$ has the description $\rmZ(A) = \ts\bigoplus_{\ga \in \Ga}\, F
t^\ga$ where
 \begin{eqnarray*}%$$\eqalign{
\Ga &= \{ \ga \in \La : \ta(\ga,\mu) = \ta(\mu,\ga)
           \hbox{ for all } \mu \in \ZZ^n\} \cr
    &= \{ \ga \in \La : \ts \prod_{1\le j \le n} \, q_{ij}^{\ga_j} = 1
         \hbox{ for } 1\le i \le n\}
  \end{eqnarray*}%$$
is a subgroup of $\ZZ^n$. The centre is isomorphic to $F[\Ga]$,
hence to a Laurent polynomial ring in $m$ variables, $0\le m\le n$.
Recall from (\ref{centroid}.1) that $\rmZ(A) \cong \cent_F(A)$,
hence $\Ga$ is the centroidal grading group in the sense of
\ref{grsico}. It is also useful to know that $F_q = \rmZ(A) \oplus
[F_q,\, F_q]$ where $[F_q,\, F_q] = \Span_F \{ ab-ba : a,b\in
F_q\}$. \sm

(c) The following are equivalent, cf. \ref{grsico}(c):
 %\begin{description}

   (i) all $q_{ij}$ are roots of unity;

  (ii) $[\ZZ^n : \Ga] < \infty$;

  (iii) $F_q$ is fgc (\ref{centroid}). %\end{description}
\sm

(d) (\cite[Remark~2.45]{bgk}) Up to scalars, $F_q$ has a unique
$\ZZ^n$-graded invariant symmetric bilinear form $\inpr$, given by
$(t^\la | t^\mu) = \pi(t^\la t^\mu)$ where $\pi \colon\, F_q \to
Ft^0$ is the canonical projection. The form $\inpr$ is
nondegenerate. \sm

(e) (\cite[Cor.~2.3]{OP}, \cite[2.48--2.56]{bgk}) The derivation
algebra of $F_q$ is the semidirect product of the ideal $\IDer
(F_q)= \{ \ad a: a \in F_q\}$ of inner derivations, where $(\ad a)
(b) = ab-ba$, and the subalgebra $\CDer(F_q)$ of centroidal
derivations: $$\Der_F (F_q) = \IDer (F_q) \rtimes \CDer(F_q).$$ Both
$\IDer(F_q)$ and $\CDer(F_q)$ are graded subalgebras of $\Der_F
(F_q)$:
 \begin{eqnarray*} %$$  \eqalign{
   \IDer(F_q) &=& \ts \bigoplus_{\la \not\in \Ga} \,
                 (\Der_F (F_q))^\la \cong [F_q,\, F_q]
              = \bigoplus_{\la \not \in \Ga} (Ft^\la), \\
    \CDer (F_q) &=& \ts\bigoplus_{\la \in \Ga}\, (\Der (F_q))^\la
                   = \rmZ(F_q) \cdot \scD, \quad \hbox{where} \\
       \scD &=& \Span_F\{ \pa_1, \ldots, \pa_n\} \end{eqnarray*}%}$$
are the degree derivations of the $\ZZ^n$-graded algebra $F_q$ with
$\pa_i$ operating as $\pa_i (t^\la) = \la_i t^\la$.
\es%\endsubsec%

\subsec[pergralgno]{\bf Notes.} The notion of the centroid of an
algebra arose in the theory of forms of algebras, see e.g.
\cite[Ch.X, \S 1]{jake:lie} or \cite[Ch.II, \S  1.6]{mc:taste}. The
centroids of not necessarily split finite-dimensional simple Lie
algebras over a field of characteristic $0$ were determined in
\cite[V.1]{Se:rat}. Some recent references regarding the centroids
of Lie algebras are the papers \cite{abp2.5}, \cite{abfp} and
\cite{BN}. \sm

Th.~\ref{grcesich} is a standard result in the ungraded case, see
for example \cite[Ch.X, \S  1]{jake:lie} or \cite[Ch.II, \S  1.6,
Th.~1.7.1]{mc:taste}. The graded case can be proven by a graded
version of the Density Theorem, \cite{n:lect}. The implication (i)
$\Rightarrow$ (ii) is mentioned in \cite[Remark~4.3.2(ii)]{abfp}.
\sm

For more information about crossed product algebras see for example
the book \cite{pass}. Prop.~\ref{grsico}(a) and (b) are proven in
\cite[Prop.~2.16]{BN}. Part~(c) is shown in
\cite[Prop.~4.4.5]{abfp}; it was independently discovered by the
author and used in the proof of \cite[Th.~7(b)]{n:tori}. \sm

Theorem~\ref{perweakv} is only one of several realization theorems
proven in \cite{abfp}. The special case of Lie-$\ZZ^n$-tori is
developed in detail in \cite{abfp2}. Multiloop algebras are studied
in detail in \cite{abp1,abp2,abp2.5}. Their derivations are
determined in \cite{Azam:multi}. \sm

Quantum tori are the algebraic counterparts of the so-called
rotation algebras, a special class of $C^*$-algebras which has a
very well-developed theory. The terminology ``quantum torus'' seems
to go back to Manin \cite[Ch.~4, \S  6]{manin}. Much more is known
about quantum tori than what is mentioned in \ref{quatoex}. For
example, the structure of quantum tori $F_q$ for which all $q_{ij}$
roots of unity is elucidated in \cite[Prop.~9.3.1]{abfp},
\cite[Th.~4.8]{hartwig} and \cite[Th.~III.4]{neeb:quantumtori}. The
involutions of quantum tori are determined in \cite{y:inv}.
\endsubsec
\vfe

\section{Lie algebras graded by root systems} \label{perrsgr}

While the previous section \S\ref{pergralg} provided some background
material for arbitrary graded algebras, in this section we consider
Lie algebras with a special type of grading by the root lattice of a
root system, so-called root-graded Lie algebras. We will only
describe that part of the theory which is of importance for the
following sections. \sm

Throughout, and unless specified otherwise, all algebras are defined
over a unital commutative associative base ring $k$ in which $2$ and
$3$ are invertible. As  in \S  \ref{pergralg}, $\La$ denotes an
abelian group, written additively.

\subsec[rogr]{\bf Root-graded Lie algebras.} Let $R$ be a locally
finite root system. It is convenient to put
$$ R\ti = R\setminus \{0\}\quad\hbox{and} \quad
       R\ti\ind = R\ti \cap R\ind,
 \index{R@$R\ti$ (the non-zero roots of $R$)}
 \index{R@$R\ti\ind$ (the non-zero indivisible roots of $R$)}
$$
see \ref{resexama}. We denote the $\ZZ$-span of $R$ by $\scq(R)$.
\sm Let $L$ be a $(\scq(R), \La)$-graded Lie algebra. We will always
use subscripts to indicate the $\scq(R)$-grading and superscripts
for the $\La$-grading: $L= \bigoplus_{q\in \scq(R)} L_q$ and $L=
\bigoplus_{\la \in \La} L^\la$. By compatibility we therefore have
$$
  L= \ts\bigoplus_{q\in \scq(R), \, \la \in \La} L_q^\la
    \qquad \hbox{for } L_q^\la = L_q \cap L^\la.
$$

To stress the analogy between the definitions in \ref{assogr} and
the ones below, we will first define invertible elements of $L$. A
non-zero element $e\in L_\al^\la$ for $\al \in R\ti$ is called {\it
invertible\/}, more precisely, an {\it invertible
element\index{invertible element} of the $(\scq(R),\La)$-graded Lie
algebra $L$}, if there exists $f\in L_{-\al}^{-\la}$ such that
$h=[f,\,e]$ acts on $L_q$, $q\in \scq(R)$, by
$$  [h, \, x_q] = \lan q,\al\ch\ran \, x_q
   \qquad (x_q\in L_q). \eq{(1)}%iel1}
$$
In this case, $(e,h,f)$ is an ${\mathfrak s}{\mathfrak l}_2$-triple.
If $L_{3\al} = \{0\}$ (a condition which will always be fulfilled
later), then $f\in L_{-\al}^{-\la}$ is uniquely determined by (1).
We will therefore refer to $f$ as the {\it inverse\/}\index{inverse
of an invertible element} of $e$. \ms

A $(\scq(R),\La)$-graded Lie algebra $L= \bigoplus_{q\in \scq(R), \,
\la\in \La} L_q^\la$ is called {\it $(R,\La)$-graded\/}\index{graded
algebra!root-graded}\index{root-graded Lie
algebra}\index{R@$(R,\La)$-graded Lie algebra} if it satisfies the
axioms (RG1)--(RG3) below: \sm

(RG1) $\supp_{\scq(R)} L \subset R$, whence $L = \bigoplus_{\al \in
R} L_\al^\la$, \sm

(RG2) every $L_\al^0$, $\al\in R\ind\ti$,  contains an invertible
element of $L$, and \sm

(RG3) $L_0 = \sum_{0 \ne \al \in R}\, [L_\al, \, L_{-\al}] $.
 \ms

\noindent An $(R,\La)$-graded Lie algebra $L$ is called \sm

\begin{itemize}
\item {\it predivision-$(R,\La)$-graded\/}\index{graded
algebra!root-graded!predivision}\index{predivision-root-graded Lie
algebra}\index{predivision-$(R,\La)$-graded Lie algebra} if every
non-zero $L_\al^\la$, $\al \in R\ti$, contains an invertible
element;

\item {\it division-$(R,\La)$-graded\/}\index{graded
algebra!root-graded!division}\index{division-root-graded Lie
algebra}\index{division-$(R,\La)$-graded Lie algebra} if every
non-zero element of $L_\al^\la$, $\al \in R\ti$, is invertible; \sm

\item {\it a
Lie-$(R,\La)$-torus}\index{Lie-$(R,\La)$-torus}\index{torus!Lie}
\index{graded algebra!Lie torus} if $L$ is predivision-graded,
defined over a field $k$ and $\dim_k L_\al^\la \le 1$ for all
$\al\in R\ti$; \sm

\item {\it invariant\/}\index{invariant!root-graded Lie
algebra}\index{graded
algebra!root-graded!invariant}\index{root-graded Lie
algebra!invariant} if $L$ has a $\La$-graded invariant nondegenerate
symmetric bilinear form $\inpr$ such that for some choice of
invertible elements $e_\al\in L_\al^0$, $\al \in R\ind\ti$, with
inverses $f_\al$ the restriction of $\inpr$ to $\Span_k \{ [e_\al,
f_\al]: 0\ne \al \in R\ind\}$ is nondegenerate too. \end{itemize}\ms

\noindent Let $L$ be $(R,\La)$-graded. By (RG1) and (RG2) we have
$$R\ind \; \subset \; \supp_{\scq(R)} L = \{\al\in R : L_\al \ne 0 \}
\; \subset \; R,
$$
whence $L$ has full $\scq(R)$-support: $\lan \supp_{\scq(R)} L \ran
= \scq(R)$. One can show that $\supp_{\scq(R)} L$ is a subsystem of
$R$. We can therefore always assume that $\supp_{\scq(R)}L =R$, if
we so wish. Analogously, we can always assume that $L$ has full
$\La$-support in the sense that $\La= \lan \supp_\La L \ran$. \ms

If in the definitions above $\La=\{0\}$, we will speak of an {\it
$R$-graded\/}\index{R@$R$-graded Lie algebra}\index{root-graded Lie
algebra} or a {\it (pre)division-$R$-graded\/}\index{R@$R$-graded
Lie algebra!(pre)division}\index{predivision-$R$-graded Lie
algebra}\index{division-$R$-graded Lie algebra} Lie algebra. Lie
algebras that are $R$-graded for some locally finite root system $R$
are simply called {\it root-graded\/}.\index{graded
algebra!root-graded}\index{root-graded Lie algebra}\index{Lie
algebra!root-graded} Similarly, a {\it predivision-root-graded Lie
algebra\/}\index{graded
algebra!root-graded!predivision}\index{predivision-root-graded Lie
algebra}\index{graded algebra!predivision-graded (Lie)}\index{graded
algebra!division-graded (Lie)} is a Lie algebra which is
predivision-$(S,\La)$-graded for some $(S,\La)$, and a {\it Lie
torus\/}\index{torus!Lie} \index{graded algebra!torus!Lie} is a
Lie-$(R,\La)$-torus for a suitable pair $(R,\La)$. \ms

We will present examples of root-graded Lie algebras in
\ref{predivex} and in \S\ref{perexass}. But we mention here right
away the following examples: the derived algebra of an affine
Kac-Moody algebra (see \ref{afrreal}), toroidal Lie algebras,
Slodowy's intersection matrix algebras and finite-dimensional simple
Lie algebras over fields of characteristic $0$ which contain an
${\mathfrak s}{\mathfrak l}_2$-subalgebra. Moreover,
predivision-graded Lie algebras will turn out to be the basic
building blocks of the generalizations of affine Lie algebras that
we will be
describing in section \S\ref{pereala}. \endsubsec%\es

\subsec[isot]{\bf Isomorphisms and isotopy.} An $(R,\La)$-graded Lie
algebra $L$ is called {\it
isogra\-ded-isomorphic\/}\index{isograded-isomorphic}\index{graded
algebra!isograded-isomorphic}\index{R@$(R,\La)$-graded Lie
algebra!isograded-isomorphic} to an $(R',\La')$-graded Lie algebra
$L'$ if there exists an isomorphism $\ph_a \colon\, L \to L'$ of Lie
$k$-algebras, an isomorphism $\ph_r \colon\, R \to R'$ of locally
finite root systems and an isomorphism $\ph_e \colon\, \La \to \La'$
of abelian groups satisfying
$$
     \ph_a (L_\al^\la) = L_{\ph_r(\al)}^{\ph_e(\la)}
 $$
for all $\al\in R$ and $\la \in \La$. We say $L$ and $L'$ are {\it
graded-isomorphic\/}\index{graded-isomorphic}\index{graded
algebra!graded-isomorphic}\index{R@$(R,\La)$-graded Lie
algebra!graded-isomorphic} if there exists an isograded isomorphism
with $\ph_r= \Id$ and $\ph_e= \Id$. \sm

Let $\io\colon\, \scq(R) \to \La$ be a group homomorphism. We can
then define a new $(\scq(R),\La)$-graded Lie algebra $L^{(\io)}$ by
$$  \big( L^{(\io)} \big)^\la_\al = L_\al^{\la + \io(\al)}.
$$
The following are equivalent: \sm

  (i) $L^{(\io)}$ is $(R,\La)$-graded, \sm

 (ii) $L_\al^{\io(\al)}$ contains an invertible element for
all $\al \in R\ind\ti$. \sm

\noindent If $R$ has a root basis and, in addition to $\frac{1}{2}$,
$\frac{1}{3}$ also $\frac{1}{5} \in k$, then (i) and (ii) are
equivalent to \sm

 (iii) $L_\al^{\io(\al)}$ contains an invertible element for all
$\al\in B$. \sm

\noindent In this case, we say that $L^{(\io)}$ is an {\it
isotope\/}\index{isotope}\index{Rla@$(R,\La)$-graded Lie
algebra!isotope} of $L$. An $(R',\La')$-graded Lie algebra $L'$ is
called {\it isotopic\/}\index{isotopic}\index{Rla@$(R,\La)$-graded
Lie algebra!isotopic} to $L$ if some isotope of $L$ is
isograded-isomorphic to $L'$. Being isotopic is an equivalence
relation on the class of all root-graded Lie algebras. \sm

Isotopy preserves the classes of (pre)division-graded Lie algebras
and Lie tori. For these cases, condition (ii) can be replaced by
$L^{\io(\al)} \ne 0$ for all $\al \in R\ind\ti$, and similarly for
(iii).
\es%\endsubsec

\subsec[centex]{\bf Coverings and central quotients.} Let $\Xi$ be
an abelian group and let $L$ be a $\Xi$-graded Lie algebra. By
definition, a {\it $\Xi$-graded central extension\/}\index{central
extension!graded}\index{extension!central extension!graded} of $L$
is a $\Xi$-graded Lie algebra $K$ together with a graded epimorphism
$f\colon\, K \to L$ with $\Ker f$ contained in the centre of $K$. A
$\Xi$-graded central extension $f: K\to L$ is called a {\it
$\Xi$-covering\/}\index{covering}\index{central extension!graded
covering} if $K$ is perfect. We will use the terms {\it central
extension\/}\index{central extension}\index{extension!central
extension} and {\it covering\/}\index{covering}\index{central
extension!covering} in case $\Xi=\{0\}$. \sm

A central extension $\fru\colon\, \frL \to L$ is called a {\it
universal central extension}\index{universal central
extension}\index{extension!universal central extension} if for any
central extension $f\colon\, K \to L$ there exists a unique
homomorphism ${\mathfrak f}\colon\, \frL \to \frk $ such that $\fru
= f \circ \frf$:
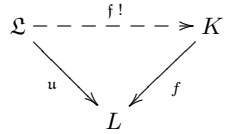
$$
\xymatrix{ \frL \ar@{-->}[rr]^{\frf\, !} \ar[dr]_\fru && K
\ar[dl]^f\cr
            & L }
$$
It is obvious from this universal property that {\it two universal
central extensions of $L$ are isomorphic\/} in the sense that one
has a commutative triangle as above with the horizontal map being an
isomorphism. A Lie algebra $L$ has a universal central extension,
say $\fru \colon\, \uce(L) \to L$,\index{U@$\uce$ (universal central
extension)} if and only if $L$ is perfect. In this case, $\uce(L)$
is perfect and is a universal central extension of all covering
algebras of $L$ as well as all central quotients. In particular, it
is a universal central extension of the centreless algebra
$L/\rmZ(L)$. If $L$ is $\Xi$-graded, then $\fru\colon\, \uce(L) \to
L$ is a $\Xi$-covering. One calls $L$ {\it simply
connected\/}\index{simply connected}\index{Lie algebra!simply
connected} if $L$ is a universal central extension of itself. For
example, a universal central extension is always simply connected.
\sm

As an example we mention that a universal central extension of the
untwisted and twisted loop algebra $L$ of \ref{afrreal} is the Lie
algebra $K$ defined there, \cite{gar,wilson}. \ms

Let now $L$ be an $(R,\La)$-graded Lie algebra. We will apply the
concepts above to the $\Xi$-grading of $L$ where $\Xi= \scq(R)
\oplus \La$. We observe first that the centre of $L$ is $\La$-graded
and contained in $L_0$,
$$   \rmZ(L) = \ts\bigoplus_{\la \in \La} \, (\rmZ(L) \cap L_0^\la).
\eq{(1)}%centex1}
$$ whence $\rmZ(L)$ is also $\Xi$-graded. \endsubsec

 \subsec[rogrcov]{\Proposition}{\it Let $f\colon\,
K\to L$ be a $\La$-covering of the $\La$-graded Lie algebra $L$.
Then $K$ is $(R,\La)$-graded if and only if $L$ is $(R,\La)$-graded.
In this case $\supp_{\scq(R)} K = \supp_{\scq(R)} L$ and $f|K_\al$
is bijective for all $\al \in R\ti$. Moreover, $K$ is predivision-
or division-graded or a Lie torus if and only if $L$ has the
corresponding property. } \ms

As a consequence of this proposition, the classification of
root-graded Lie algebras and their subclasses introduced above can
be done in three steps: \sm

\begin{enumerate}
\item Describe all centreless root-graded Lie algebras,

\item determine their universal central extensions and
\item characterize those that are (pre)division-graded or Lie
tori. \end{enumerate}

\noindent Note that any root-graded Lie algebra is perfect and
therefore has a universal central extension. Examples of this
approach are given in \ref{predivex}-\ref{jpclas} and \S
\ref{perexass}, see also \ref{perrsgrno} for a short survey.
 \es

\subsec[predivex]{\bf Example.} Let $\frg_\ZZ$ be a Chevalley order
of a finite-dimensio\-nal semisimple complex Lie algebra $\frg_\CC$
of type $R$, see \cite[Ch. VIII, \S  12]{bou:lie78}. Thus, we have a
Chevalley system $(e_\al : 0\ne \al \in R)$ such that, putting
$h_\al = [e_{-\al}, e_\al]$,
$$
   \frg_\ZZ = \big(\ts\bigoplus_{0\ne \al \in R} \ZZ e_\al \big)
   \oplus \big(\sum_{0\ne \al\in R} \ZZ h_\al\big)
$$
is a Lie algebra over $\ZZ$ which is a $\ZZ$-form of $\frg_\CC$,
i.e. $\frg_\ZZ \otimes _\ZZ \CC \cong \frg_\CC$.  Then $\frg_k =
\frg_\ZZ \otimes_\ZZ k$ is a $\scq(R)$-graded Lie $k$-algebra with
homogeneous spaces $\frg_{k,\,\al}= \frg_{\ZZ,\al} \otimes_\ZZ k$.
For every unital commutative associative $\La$-graded $k$-algebra
$C= \bigoplus_{\la \in \La}\, C^\la$ the Lie $k$-algebra
$$
   L =\frg_k \otimes_k C\cong \frg_\ZZ \otimes_\ZZ C \eq{(1)}% predivex1}
$$
is $(\scq(R),\La)$-graded with homogeneous spaces $L^\la_\al =
\frg_{k,\al} \otimes_k C^\la$. An element $e_\al \ot c$ is
invertible in $L$ if and only if $c\in C\ti$. Hence, $L$ is
$(R,\La)$-graded. It is (pre)division-$(R,\La)$-graded if and only
if $C$ is (pre)division-graded. If $k$ is a field (recall of
characteristic $\ne 2,3$), then $L$ is a Lie torus if and only if
$C$ is an associative torus. \sm

Let $\ka$ be the Killing form of $\frg_\CC$ and let $h\ch$ be the
dual Coxeter number. Then $\frac{1}{2 h\ch} \ka$ is an invariant
nondegenerate symmetric bilinear form on $\frg_\ZZ$
(\cite[I.4.8]{springer-steinberg}, \cite[\S  5]{gross-nebe}). The
canonical extension $\inpr$ of this form to $\frg_k$ is then an
invariant symmetric bilinear form on $\frg_k$. It is nondegenerate
if $R$ is irreducible and not of type $\rmfa$. In this case $L$ is
invariant with respect to $\inpr\ot \eps$, where $\eps$ is any
$\La$-graded invariant nondegenerate symmetric bilinear form on $C$,
see (\ref{pergradef}.2).  \sm

If $R$ is irreducible, the Lie algebra $\frg_k$ is simply connected
(\cite[6.1]{stein:gen}, \cite[Cor.~3.14]{vdK}, \cite{Kas}).  The
universal extension $\fru \colon\, \uce(L)\to L$ of an arbitrary $L$
as in (1) is determined in \cite[3.8]{Kas}. It turns out that the
kernel of $\fru$ is $\Om_{C/k}\big/ dC$ where $\Om_{C/k}$ is the
module of K\"ahler differentials of $C$ over $k$ and $d\colon\, C
\to \Om_{C/k}$ the universal $k$-linear derivation. In view of
(\ref{uce}.1) it is appropriate to point out that $\Om_{C/k}\big/
dC= \HC_1(C)$\index{H@$\HC_1(C)$ (first cyclic homology group of
$C$)}, the first cyclic homology group of $C$. \ms

For root systems of type $\rmd$ or $\rme$ the algebras of this
subsections are in fact all examples of $R$-graded Lie algebras.
\endsubsec

\subsec[declas]{\Theorem} (\cite{bm}, \cite{n:3g}) {\it Suppose $R$
is a finite root system of type $\rmd_l$, $l \ge 4$ or $\rme_l$,
$l=6$, $7$ or $8$, and that $k$ is a field of characteristic $0$ in
case of $\rme_8$. Then any $(R,\La)$-graded Lie algebra is
graded-isomorphic to a $\La$-covering of a Lie algebra $L$ as in\/
{\rm (\ref{predivex}.1)} for a suitable $C$.} \ms

If in case $\rmd$ one replaces $\frg$ by the corresponding infinite
rank Lie algebra (\ref{afrreal}), the result remains true for
$(\rmd_I,\La)$-graded Lie algebras, $|I|=\infty$. But it is no
longer true for the other types of root systems, see \S
\ref{perexass}, in particular \ref{grtypeAcl} and \ref{grtypeAco},
for $R$ of type $\rma_I$ and \ref{perrsgrno} for a short summary of
the classification results. The models used to describe $R$-graded
Lie algebras are different for every type of $R$ and involve all
important classes of nonassociative algebras (Jordan algebras,
alternative algebras, structurable algebras). On the other side, one
has the following type-free approach for most root systems, not only
finite ones.
\endsubsec

\subsec[jpclas]{\Theorem} (\cite{n:3g}) {\it Let $R$ be a locally
finite root system without an irreducible component of type
$\rme_8$, $\rmf$ or $\rmg$. Then a Lie algebra $L$ is $R$-graded if
and only if $L$ is a covering of the Tits-Kantor-Koecher algebra of
a Jordan pair covered by a grid whose associated root system is $R$.
}
%\es%
\endsubsec

\subsec[rocent]{\bf Centroids of root-graded Lie algebras.} Let $L=
\bigoplus_{\al \in R, \, \la\in \La} L_\al^\la$ be an
$(R,\La)$-root-graded Lie algebra over $k$. \sm

(a) Since $L$ is perfect, the centroid of $L$ is a commutative (and
of course also unital and associative) subalgebra of $\End_k(L)$.
One can show that a centroidal transformation of $L$ leaves every
$L_\al$, $\al \in R$, invariant. If $R$ has only finitely many
irreducible components, $\cent_k(L)$ is $\La$-graded. As an example,
we mention that the centroid of the Lie algebra $L$ of
(\ref{predivex}.1) is isomorphic to $C$. \sm

(b) Suppose now that the $\La$-graded algebra $L$ is graded-simple.
Then $R$ is irreducible and $L$ is also graded-simple with respect
to the $\scq(R)$-grading of $L$. Hence, the centroid of $L$ is a
twisted group algebra $K^t[\Ga]$ over the field $K=(\cent L)^0$, as
described in \ref{grsico}. The condition\/ {\rm (\ref{grsico}.2)}
characterizing fgc algebras now simplifies to the following, where
$\La_\al = \{ \la \in \La : L_\al^\la \ne 0\}$: $L$ is fgc if and
only if \sm

 (i) $R$ is finite,

  (ii) $[\La_\al : \Ga]<\infty$ for all short roots $\al \in
R$, and

  (iii) $\dim_K L_\al^\la < \infty$ for all short roots $\al$
and all $\la \in \La_\al$. \sm

\noindent If also $\frac{1}{5} \in k$, then for every $w\in W(R)$
there exists a $\cent (L)$-linear automorphism of $L$ which for all
$\al \in R$ and $\la \in \La$ maps $L_\al^\la$ to $L^\la_{w(\al)}$.
Therefore, in this case it is sufficient to require (ii) and (iii)
for one short root $\al$. \sm

If $L$ is in addition predivision-graded and $\La=\lan \supp_\La
L\ran$ is finitely generated, then the condition (ii) can be
replaced by

 (ii$'$)  $[\La :\Ga]<\infty$. \sm

%Finally, if $k$ is a field and $\dim_k L_\al^\la = 1$ for some
%$\al\in ^\la$ then $K=k$ and hence $L$ is graded-central-simple. \sm

(c) Let $L$ be division-graded with $\rmZ(L)=\{0\}$ and suppose that
$R$ is irreducible. Then $L$ is $\La$-graded-simple. So its centroid
is described in (b). Moreover, if $\La$ is torsion-free, $L$ is a
prime algebra, i.e., $[I,J]=0$ for ideals $I,J$ of $L$ implies $I=0$
or $J=0$. In this case, $\cent_k(L)$ is an integral domain, and $L$
embeds into the $\tl C$-algebra $\tl L = L \ot \tl C$, where $\tl C$
is the quotient field of $\cent_k(L)$. The algebra $\tl L$ is
central and prime (\cite{emo}). If the Lie algebra $L$ is fgc, then
$\tl L$ is finite-dimensional. Hence, if $C$ has characteristic $0$,
then $\tl L$ is a finite-dimensional central-simple $\tl C$-algebra.
 \sm

(d) Let $L$ be a centreless Lie-$(R,\La)$-torus and suppose that $R$
is irreducible. Then $L$ is graded-central-simple by (c) and (b). In
particular, if $\La=\lan \supp_\La L\ran$ is finitely generated,
then $L$ is fgc if and only if $R$ is finite and $[\La :
\Ga]<\infty$. See also \ref{lietor}. \ms

We will now consider root-graded Lie algebras in characteristic $0$.
It will turn out that the $\scQ(R)$-grading of a root-graded Lie
algebra is in fact the root space decomposition with respect to a
toral subalgebra, as we will explain now.
\endsubsec%\es

\subsec[perschn0]{\Proposition} {\it  Let $L$ be an $(R,\La)$-graded
Lie algebra, defined over a field $\KK$ of characteristic $0$. We
suppose $\supp_{\scq(R)} L = R$. For $\al \in R\ind\ti$ let $e_\al
\in L_\al^0$ be an invertible element with inverse $f_\al$. We put
$h_\al = [f_\al,\, e_\al]$, $h_0 = 0$, $h_{2\al} = \frac{1}{2}
h_\al$ in case $2\al \in R$, and
$$
    \frh = \Span_\KK \{ h_\al : \al \in R\}.
$$

{\rm (a)} For every $\al \in R$ there exists a unique $\tl\al \in
\frh^*$ defined by $\tl \al (h_\be) = \lan \al,\be\ch\ran$ for
$\be\in R$. The set $\tl R = \{ \tl\al \in \frh^* : \al \in R\}$ is
a locally finite root system, isomorphic to $R$ via the map $\al
\mapsto \tl\al$. \sm

{\rm (b)} The set $\{ h_\al : \al \in R\}$ is a locally finite root
system in $\frh$, canonically isomorphic to the coroot system $R\ch$
via the map $\al\ch \mapsto h_\al$. \sm

{\rm (c)} The subspace $\frh$ is a toral subalgebra of $L$ whose
root spaces are the homogeneous subspaces $L_\al$ with corresponding
linear form $\tl\al$:
$$
   L_\al = \{ l \in L : [h,\, l] = \tl\al(h) l \hbox{ for all }
    h\in \frh\}.
$$

{\rm (d)} Let $R$ be finite and let $B$ be a root basis of $R$. Then
$\{e_\be, f_\be : \be \in B\}$ generates a finite-dimensional split
semisimple subalgebra $\frg$ of $L$ with splitting Cartan subalgebra
$\frh$. The root system of $(\frg,\frh)$ is $\tl R\ind$.} \ms

It follows immediately from this result that for $\KK$ a field of
characteristic $0$, an $(R,\La)$-graded Lie $\KK$-algebra for $R$ an
irreducible root system and $\La$ an abelian group can be
equivalently defined as follows: $L= \bigoplus_{\la \in \La} L^\la$
is a $\La$-graded Lie algebra over $\KK$ satisfying (i)--(iv) below:
\sm

\begin{list}{}{\setlength{\leftmargin}{\parindent}}

 \item[(i)] $L^0$ contains as a subalgebra a finite-dimensional split
simple Lie algebra $\frg$, called the {\it grading
subalgebra\/},\index{grading subalgebra}\index{root-graded Lie
algebra!grading subalgebra} with splitting Cartan subalgebra $\frh$,
\sm

 \item[(ii)] either $R$ is reduced and equals the root system
$R_\frg\subset \frh^*$ of $(\frg, \frh)$, or $R= \rmbc_I$ and
$R_\frg=\rmfa_1$ in case $|I|=1$ or $R_\frg=\rmb_I$ for $|I|\ge
2$.\sm

\item[(iii)] $\frh$ is a toral subalgebra of $L$, and the weights of
$L$ relative $\frh$ are in $R$, whence $L = \bigoplus_{\al \in R}
L_\al$, \sm

\item[(iv)] $L_0 = \bigoplus_{0 \ne \al \in R} [L_\al, \, L_{-\al}]$;
\end{list} \sm

\noindent This is the original definition of a root-graded Lie
algebra, see \ref{perrsgrno}.
\es%\endsubsec

\subsec[perrsgrno]{\bf Notes.} Unless stated otherwise, the results
mentioned in \ref{rogr}--\ref{perschn0} are proven in \cite{n:lect}.
Some of them have also been proven by others, as indicated below.
\sm

 The definition of a root-graded Lie
algebra in \ref{rogr} is in the spirit of the definition in
\cite{n:3g}, which introduces root-graded Lie algebras over rings
for reduced root systems. The definition given at the end of
\ref{perschn0} is the original one of \cite{bm}, where $R$-graded
Lie algebras were introduced $R$ simply laced. It is mentioned in
\cite[Remark 2 of 2.1]{n:3g} that in characteristic $0$ the two
definitions are equivalent. $\rmfa_1$-graded Lie algebras had been
studied much earlier in \cite{tits:JA} and then in
\cite{Kan1,Kan2,Kan3} and \cite{Koe1,Koe2}. \sm

The paper \cite{bm} contains a classification of centreless
$R$-graded Lie algebras for $R$ simply laced. The classification for
reduced, non-simply laced $R$ is given in \cite{beze}. The theory
for a non-reduced $R$, i.e. $R=\rmbc_I$, is developed in \cite{abg2}
for $R= \rmbc_l$, $2\le l< \infty$) and \cite{BeSm} for $R=\rmbc_1$.
In this case, one can be more general by allowing the grading
subalgebra $\frg$ to have type $\rmC$ or $\rmd$. This too can be
done in the setting of \ref{rogr}, but we will not need this in the
following and have therefore refrained from describing the more
general setting. The papers \cite{bm}, \cite{beze}, \cite{abg2} and
\cite{BeSm} work with finite root systems and Lie algebras over
fields of characteristic $0$. The extensions to locally finite root
systems and Lie algebras over rings is given in \cite{n:3g} for $R$
reduced and $3$-graded.  As already mentioned in \cite[2.10]{gn2},
the classification results in \cite{n:3g}, viz. \ref{jpclas}, can
easily be extended to the $(R,\La)$-graded and
$(R,\La)$-(pre)division-graded case. \sm

At present, predivision-$(R,\La)$-graded Lie algebras have been
classified in the centreless case only for certain root systems and
base rings $k$: $R$ finite of type $\rmfa_l$, $(l\ge 3)$, $\rmd$ or
$\rme$, $\La=\ZZ^n$, $k=\KK$ a field of characteristic $0$ in
\cite{y2}; $R=\rmb_r$, $3\le r< \infty$, $L$ division-graded, $k=
\KK$ a field of characteristic $0$ in \cite{y:ext}; $R=\rmc_r$,
$2\le r < \infty$, $L$ division-graded over $k=\KK$ a field of
characteristic $0$ in \cite{BeYo}; $R=\rmb_2=\rmc_2$, $\La$
torsion-free, $L$ graded-simple (e.g. division-graded), $k$
arbitrary in \cite{NeTo}. Regarding the classification of Lie tori,
the situation the reader is referred to the survey in \cite[\S
12]{AF:isotopy}. \sm

The concept of isotopy was introduced in \cite{AF:isotopy} for Lie
tori over fields of characteristic $0$ and related to isotopy of
their coordinate algebras. In the context of Lie tori in
characteristic $0$, the equivalence of (i)--(iii) in \ref{isot} is
proven in \cite[Prop.~2.2.3]{abfp2}. \sm

Universal central extensions of Lie algebras over rings were
described in \cite[\S  1]{vdK} (this reference also describes the
universal central extensions of the Lie algebras $\frg_k$ of
\ref{predivex}, where $k$ not necessarily contains $\frac{1}{2}$ and
$\frac{1}{3}$. Later references on universal central extensions are
\cite[\S  1]{gar}, \cite[1.9]{mp}, \cite[7.9]{wei} or \cite[]{n:uce}
(the list is incomplete). In characteristic $0$, the universal
central extensions of root-graded Lie algebras are determined in
\cite{abg}, \cite{abg2} and \cite{BeSm}. \sm

Prop.~\ref{rogrcov} is standard, proven in \cite[Th.~5.36]{abg2},
\cite[\S  5]{BeSm},  \cite[Prop.~1.5--1.8]{beze} and
\cite[Prop.~1.29]{bm} for the root-graded Lie algebras studied in
these papers, viz. Lie algebras over fields of characteristic $0$,
where one can use that the grading subalgebra $\frg$ is simply
connected. The proof of \ref{rogrcov} in \cite{n:lect} is for
arbitrary $k$ and uses the model of the universal central extension
from \cite{vdK}.  \sm

\ref{predivex}: We give some more information on this example in
case $k=F$ is a field (recall of characteristic $\ne 2,3$). It
follows as in \ref{predivex} that $\frg_F$ is central. Moreover, the
description of $\Der_F L$ given there remains true, except that
$\Der_F (\frg_F)$ need not be equal to $\IDer \frg_F$ and hence
$\IDer L \ne \Der_F (\frg_F) \ot C$ in general. If $R\ne \rma$, $L$
is centreless, see \ref{ccc} for the case $R=\rma$. \sm

Th.~\ref{declas} is proven in \cite{bm} for $k$ a field of
characteristic $0$ and in \cite[5.4, 7.2, 7.3]{n:3g} for arbitrary
$k$ (of course containing $\frac{1}{2}$ and $\frac{1}{3}$) and $R\ne
\rme_8$. The methods of \cite{n:3g} do not allow to treat the case
$\rme_8$, but the result is probably also true in this case for any
$k$ containing $\frac{1}{30}$. Similarly, Th.~\ref{jpclas} can be
generalized to cover all locally finite root systems by using Kantor
pairs instead of Jordan pairs. With this method one can likely
extend the classification results to cover all root-graded Lie
algebras ($R$ not $3$-graded and Lie algebras over rings). \sm

\ref{rocent}: That a centreless division-$(R,\La)$-graded Lie
algebra with $R$ irreducible is graded-simple is proven in
\cite[Lemma~4.4]{y:ext} for finite $R$ (the proof easily
generalizes). The characterization of fgc Lie tori in
\ref{rocent}(d) is proven in \cite[Prop.~1.4.2]{abfp2}. For $R$ a
finite root system, $R\ne \rmbc_1$ or $\rmbc_2$, the centroid of an
$R$-graded Lie algebra over a field of characteristic $0$ is
described in \cite[\S  5]{BN} in terms of the centre of the
coordinate system associated to an $R$-graded Lie algebra. \sm

For $3$-graded $R$, Prop.~\ref{perschn0}(a) was shown in
\cite[3.2]{n:gen}. It is also proven in \cite[Prop.~1.2.2]{abfp2}
for Lie tori. \sm

Invariant forms and derivations of root-graded Lie algebras over
fields of characteristic $0$ are described in \cite{be:derinv},
\cite{abg2} and \cite{BeSm}. \sm
\endsubsec
 \vfe

\section{Extended affine Lie algebras and generalizations}
 \label{pereala}

In this section we will introduce Lie algebras whose set of roots
are affine reflection systems (\ref{ARLA}). We will then discuss
special types, most importantly invariant affine reflection algebras
(\ref{invax}) and extended affine Lie algebras (\ref{ealadef}). \sm

Unless specified otherwise, in this section all algebras are defined
over a field $\KK$ of characteristic $0$. Throughout, $(E,T)$ will
be a {\it toral pair\/},\index{toral pair} i.e., a non-zero Lie
algebra $E$ with a toral subalgebra $T$. We will denote the set of
roots of $(E,T)$ by $R$ and put $X= \Span_\KK(R) \subset T^*$, cf.
\ref{inro}. \bs

\subsec[impodef]{\bf Core and tameness.} Let $(E,T)$ be a toral
pair. Recall that $0\in R$ if $T\ne 0$ and that $E$ has a root space
decomposition
$$ E = \ts \bigoplus_{\al \in R} \, E_\al. \eq{(1)}%impodef0}
$$
If $(R,X)$ is a pre-reflection system with respect to some set of
real roots $R\re\subset R$, we define the {\it core of\/
$(E,T)$}\index{core}\index{toral pair!core of} as the subalgebra
$E_c$ of $E$ generated by $\{E_\al : \al \in R\re\}$. It is a graded
subalgebra with respect to the root space decomposition (1),
$$   E_c = \ts \big( \bigoplus_{\al \in R\re}\, E_\al\big) \oplus
   \big( \bigoplus_{\be \in R\im }\, (E_c \cap E_\be) \big),
   \eq{(2)}%impodef1}
$$
and also with respect to any $\La$-grading of $E$ which is
compatible with the root space decomposition (1). If $E_c$ is
perfect (which will be the case from \ref{ARLA} on), the quotient
algebra
$$E_{cc} = E_c/\rmZ(E_c)\eq{(3)}%impodef5}
$$ is centreless and is called the {\it centreless
core}.\index{centreless!core}\index{toral pair!centreless core of}
We call $(E,T)$ {\it tame\/}\index{tame toral pair}\index{toral
pair!tame} if $\rmC_E(E_c) \subset E_c$. Tameness of $(E,T)$ is not
the same as tameness of $R$ as defined in \ref{sbc}, cf.
\ref{afrsex}, \ref{invth} and \ref{invcop}.
\es%\endsubsec

\subsec[ARLA]{\bf Affine reflection Lie algebras.} Let $\al$ be an
integrable root of $(E,T)$. Hence, by definition, there exists an
integrable ${\mathfrak s}{\mathfrak l}_2$-triple $(e_\al, h_\al,
f_\al) \subset E_\al \times T \times E_{-\al}$ of $(E,T)$. It gives
rise to a reflection $s_\al$ of $(R,X)$, defined by
$$ s_\al (x) = x - x(h_\al) \al, \quad \hbox{whence } \lan
x,\al\ch\ran = x(h_\al). \eq{(1)}%ARLA0}
$$ We will call $(E,T,R\re)$
or simply $(E,T)$ an {\it affine reflection Lie
algebra\/}\index{affine reflection Lie algebra}\index{Lie
algebra!affine reflection} (abbreviated ARLA)\index{ARLA|see{affine
reflection Lie algebra}} if $(R,X)$ is an affine reflection system
with $R\re\subset R\int$ and the reflections $s_\al$, $\al \in
R\re$, defined by (1). Thus $R\re\subset R\int$, but not necessarily
$R\re=R\int$. \sm

Because of (\ref{affdef}.1), other choices of integrable ${\mathfrak
s}{\mathfrak l}_2$-triples for $\al \in R\re$ lead to the same
reflection $s_\al$. We will therefore  not consider the family
$(e_\al: \al \in R\re)$ as part of the structure of an affine
reflection Lie algebra. Following this point of view, two affine
reflection Lie algebras $(E,T, R\re)$ and $(E',T',\allowbreak
R'{}\re)$ are called {\it isomorphic\/}\index{affine reflection Lie
algebra!isomorphic}\index{isomorphic!affine reflection Lie algebras}
if there exists a Lie algebra isomorphism $\ph \colon\, E \to E'$
with $\ph(T) = T'$ and $\ph_r(R\re) = R'{}\re$ where $\ph_r(\al) =
\al \circ (\ph|T)^{-1}$.\ms

Recall from \ref{ext} and \ref{affdef} that any affine reflection
system is integral, coherent and has finite root strings. The affine
reflection systems occurring as roots in an affine reflection Lie
algebra have the following additional properties: \sm

\iitem{$\bull$} all root strings $\SS(\be,\al)$ for $\be \in R$,
$\al \in R\re$, are unbroken. \sm

\iitem{$\bull$} $R$ is reduced in case $T$ is a splitting Cartan
subalgebra. \sm

By definition, the {\it nullity of $(E,T,R\re)$\/}\index{affine
reflection Lie algebra!nullity}\index{nullity!of an affine
reflection Lie algebra} is the nullity of the affine reflection
system $R$, i.e., the rank of the subgroup $\La$ of $(T^*, +)$
generated by $R\im$, cf. \ref{affdef}. We always have
$$
 \hbox{nullity}\quad \ge \quad \dim_\KK\big(
\Span_\KK(R\im) \big).  \eq{(2)}%arlaine}
$$
In general this is a strict inequality, see \ref{LEALAdef}. \es

\subsec[arlacore]{\bf The core of an affine reflection Lie algebra.}
 To study the core of an affine reflection Lie algebra
$(E,T,R\re)$ we need some more notation:  \sm

\begin{list}{}{\setlength{\leftmargin}{\parindent}}

\item{$\bull$} $(S,Y)$ is the quotient root system of $(R,X)$
 (\ref{affdef}); \sm

\item{$\bull$} $f\colon\, (R,X) \to (S,Y)$ is the corresponding
extension (\ref{ext});\sm

\item{$\bull$} $Z=\Ker (f)= \Span_\KK(R\im)$; \sm

\item{$\bull$} $g\colon\, Y \to X$, $y\mapsto g(y)=:\dot y$, is a
partial section over $S\ind$ (\ref{parsecdef}); \sm

\item{$\bull$} $\frL = (\La_\xi)_{\xi\in S}$ is the extension datum
associated to $f$ and $g$, see \ref{gexda}--\ref{affexd}; thus $
f^{-1}(\xi) \cap R = g(\xi) + \La_\xi$ and the set $R$ can be
decomposed as
$$    R = \ts \bigcup_{\xi \in S} \, \big( \dot \xi \oplus \La_\xi \big)
   \, \subset \, X = g(Y) \oplus Z \quad\hbox{with}\quad
R\im = \dot 0 \oplus \La_0 = \La_0. \eq{(1)}%ARLA2}
$$ \sm

\item{$\bull$} $\La$ is the subgroup of $(Z,+)$ generated by
$\bigcup_{\xi \in S} \La_\xi$, \ms

\item{$\bull$} $K=E_c$ is the core of $(E,T,R\re)$, and $K_\xi^\la
\subset E_{g(\xi) \oplus \la}$, $K_\xi$ and $K^\la$ are defined for
$\xi \in S$ and $\la \in \La$ by
\begin{eqnarray*}
         K_\xi^\la &=& E_{\dot \xi \oplus \la}\,, \qquad\qquad(\xi \ne 0)\\
      K_0^\la  &= &\ts\sum_{\tau \in S\ti, \, \mu \in \La} \,
        [E_{\dot \tau \oplus \mu} , E_{-\dot \tau \oplus (\la -\mu)}]\,, \\
      K_\xi &= & \ts \bigoplus_{\la \in \La}\, K^\la_\xi,
                 \qquad(\xi \in S),    \\
      K^\la &=& \ts \bigoplus_{\xi \in S}\, K^\la_\xi,
              \qquad (\la \in \La).
\end{eqnarray*}

\item{$\bull$} To distinguish the $\La$-support sets of $K$ and the
extension datum $(\La_\xi)_{\xi \in S}$  we put
$$%\eqalign{
    \La^K_\xi =  \{ \la \in \La : K_\xi^\la \ne 0\}
     \quad \hbox{and} \quad %, \cr
    \La^K = \Span_\ZZ \{ \La_\xi^K : \xi \in S\}. $$ \end{list}
\ms

\noindent The data $S$, $f$, $Z$ and the Lie algebra $K$ are
invariants of $(E,T,R\re)$. Moreover, since $K_\xi=
\bigoplus_{f(\al) = \xi} E_\al $ for $\xi \in S\ti$, also the
subspaces $K_\xi$, $\xi \in S$, are uniquely determined by
$(E,T,R\re)$. But this is not so for the partial section $g$ and
hence for the family $(\La_\xi)_{\xi\in S}$. Part (d) of the
following Theorem~\ref{arlath} says that another choice of $g$ leads
to an isotope of $K$ as defined in
\ref{isot}. \es%\endsubsec

\subsec[arlath]{\Proposition}{\it We use the notation of\/ {\rm
\ref{arlacore}}. The core $K$ of an affine reflection Lie algebra
$(E,T,R\re)$ is a perfect ideal of $E$ and a
predivision-$(S,\La^K)$-graded Lie algebra with respect to the
decomposition $$K = \ts \bigoplus_{\xi \in S, \, \la \in \La}\,
K_\xi^\la.$$ Moreover, $K$ has the following properties. \sm

{\rm (a)} The root system $S$ embeds into $\frh_K^*$ where $\frh_K =
\Span_\KK \{ h_{g(\xi)} : \xi \in S\ind\ti\}$, such that
$K=\bigoplus_{\xi \in S} K_\xi$ is the root space decomposition of
$K$ with respect to the toral subalgebra $\frh_K$ of $K$. The set of
roots of $(K,\frh_K)$ is $S$, unless $K=\{0\}$ and thus $R= R\im$
and $S=\{0\}$. \sm

{\rm (b)} The family $\frL^K=(\La^K_\xi : \xi \in S)$ of
$\La$-supports of $K$ is an extension datum of type $(S,S\ind,
\Span_\KK(\La^K))$ with $\La_\xi = \La_\xi^K$ for $0\ne \xi \in S$
and
$$
   \La_0  \, \supset \, \La_0^K = \supp_\La K
   = \ts \bigcup_{\xi \in S\ti}\, (\La_\xi + \La_\xi),
        % \supset \La_\xi^K \quad \hbox{for all  } \xi \in S,
$$
hence $\lan \La_0 \ran = \La \,\supset \,\La^K = \lan \La_0^K \ran$.
\sm

{\rm (c)} The following are equivalent:  \sm

\iitem{\rm (i)} $K$ is a Lie torus,

\iitem{\rm (ii)}  $\dim_\KK E_\al \le 1 $ for all $\al \in R\re$,

\iitem{\rm (iii)}  $[E_\al, E_{-\al} ] \subset T$ for all $\al \in
R\re$. \sm

\noindent In particular, if $T$ is a splitting Cartan subalgebra,
then $K$ and the centreless core $E_{cc}$ are Lie tori. \sm

{\rm (d)} Suppose $(E,T)$ is isomorphic to the affine reflection Lie
algebra $(E',T')$. Denote by $R'$, $S'$, $\La'$ and $K'$ the data
for $(E',T')$ corresponding to the ones introduced in \ref{arlacore}
for $(E,T)$. Then $R$ and $R'$ are isomorphic as reflection systems,
and the $(S,\La)$-graded Lie algebra $K$ is isotopic to the
$(S',\La')$-graded Lie algebra $K'$.} \ms

We have excluded $K=\{0\}$ in (a) because, by definition in
\ref{inro}, the set of roots of $(K,\frh_K)$ consists of those $\al
\in \frh^*$ for which the corresponding root space $K_\al(\frh)\ne
0$. Of course, $K = \{0\}$ can still be viewed as a root space
decomposition. \es

\subsec[afrsex]{\bf A construction of an affine reflection Lie
algebra.} We have seen in \ref{arlath} that the core of an affine
reflection Lie algebra $(E,T)$ is a predivision-$(S,\La^K)$-graded
Lie algebra, where $S$ is the quotient root system of the affine
reflection system of $(E,T,R\re)$ and where $\La^K$ is a
torsion-free abelian group. Changing notation (replacing $\La^K$ by
$\La$) we will now show that, conversely,
$$\vcenter{\amen{0.8}{every
predivision-$(S,\La)$-graded Lie algebra, where $S$ is any locally
finite root system and where $\La$ is any torsion-free abel\-ian
group, is the core of an affine reflection Lie algebra.\/} }
\eq{(1)}%afrsex1}
$$

 Our construction uses $K$ and $D$, where $K$ is a
predivision-$(S,\La)$-graded Lie algebra over $\KK$ with $\La$ a
torsion-free abelian group and where $D$ is described below. Since
$\supp_{\scq(S)} K$ is a subsystem of $S$, we can assume $
\supp_{\scq(S)}K=S$. Also, the family $\frL = (\La_\xi : \xi \in
S)$, $\La_\xi = \{ \la \in \La : K_\xi^\la \ne 0 \}$, is an
extension datum of type $S$ satisfying $\supp_\La K = \La_0$. Again
it does no harm to assume that $\La$ is spanned by $\supp_\La K$.
\sm

Let $(e_\xi^\la \in K_\xi^\la: 0\ne \xi \in S, \la \in \La_\xi)$ be
a family of invertible elements of the $(\scQ(S), \La)$-graded Lie
algebra $K$. We denote by $f_\xi^\la\in L_{-\xi}^{-\la}$ the inverse
of $e_\xi^\la$, put $h_\xi^\la = [e_\xi^\la, f_\xi^\la]\in K_0^0$
and define
$$
   \frh = \Span_\KK \{ h_\xi^0 : \xi \in S\ti\ind\}
      \; \subset \;
    T_K = \Span_\KK \{h_\xi^\la : \xi \in S\ti, \la \in \La_\xi\}.
$$
Then $\frh$ is a toral subalgebra of $K$ and $S$ canonically embeds
into the dual space $\frh^*$ such that $\lan \xi, \ta\ch \ran =
\xi(h_\ta^0)$ for $\si, \ta \in S\ind$, $\ta\ne 0$, and the root
spaces of $(K, \frh)$ are the subspaces $K_\xi$, $\xi \in S$. Since
$T_K \subset \frh + \rmZ(K)$, $T_K$ is also a toral subalgebra of
$K$ with the same root spaces as $\frh$. We view $S$
(non-canonically) embedded into $T_K^*$. Let $\scD = \{ \pa_\vth :
\vth \in \Hom_\ZZ(\La, \KK) \}$ be the space of degree derivations
of the $\La$-graded Lie algebra $K$, \ref{degd}. Since $\La$ is
torsion-free, there exists an embedding of the abelian group $\La$
into $\scD^*$ mapping $\la$ to the linear form $\ev_\la$ given by
$\ev_\la(\pa_\vth) = \vth(\la)$, \ref{degl}. \sm

Let $D \subset \scD$ be a subspace such that $\la \mapsto \ev_\la|D$
is injective, e.g. $D= \scD$. For such a $D$, the subspaces
$K^\la\subset K$ are determined by $D$, namely $K^\la = \{ x\in K :
[d,x] = \ev_\la(d) \hbox{ for all } d\in D\}$. \sm

It is now clear that $T= T_K  \oplus D$ is a toral subalgebra of the
semidirect product $E = K \rtimes D$ whose root spaces are
$$
  E_0 = K_0^0 \oplus D \quad \hbox{and} \quad
  E_{\xi \oplus \la} = K_\xi^\la \hbox{ for } (\xi,\la) \ne (0,0).
$$
where for $\xi \in S$ and $\la \in \La$ the linear form $\xi \oplus
\la \in T^*$ is given by
$$
   (\xi \oplus \la)\, (t_K \oplus d) = \xi(t_K) + \ev_\la(d)
   \quad \hbox{for $t_K \in T_K$ and $d \in D$.}
$$
Thus the set of roots of $(E,T)$ is $R = \ts \bigcup_{\xi \in S}\,
(\xi \oplus \La_\xi) \subset T^*$. \sm

We let $R\re = \{ \xi \oplus \la : \xi \in S\ti, \la\in \La_\xi\}$
and observe that $R\re \subset R\int$. The corresponding reflection
$s_\al$, $\al = \xi \oplus \la$, maps $\be = \ta \oplus \mu \in R$,
$\mu \in \La_\ta$, to
$$
   s_\al (\be) = \be - \be(h_\xi^\la) \al = \be - \ta(h_\xi^\la) \al
   = \be - \lan \ta, \xi\ch\ran \al.
$$
Hence $R$ with the reflections $s_\al$, $\al \in R\re$, is the
affine reflection system associated to $S$ and the extension datum
$(\La_\xi : \xi \in S)$ in \ref{extchar}. Thus, $(E,T)$ is an affine
reflection Lie algebra. Moreover, since $K$ is $R$-graded, it equals
the core of $(E,T)$. This proves (1). \sm

Note that $\La_0= \La_0^K$ in the notation of \ref{ARLA}. This
implies that the affine reflection system $R$ of $(E,T,R\re)$ is
tame. However, $(E,T)$ is tame if and only if $\IDer K \cap D =
\{0\}$.  \ms

We will now give a sufficient criterion for a toral pair $(E,T)$ to
be an affine reflection Lie algebra. \es%\endsubsec

\subsec[affrefth]{\Theorem}{\it Suppose $(E,T)$ satisfies the
conditions  {\rm (AR1)}--{\rm (AR3)} below.\sm

\iitem{{\rm (AR1)}} There exists a symmetric bilinear form $\inpr$
on $X= \Span_\KK(R)$, \sm

\iitem{{\rm (AR2)}} for every $\al \in R\an = \{ \al \in R : (\al |
\al) \ne 0 \}$ there exists an integrable ${\mathfrak s}{\mathfrak
l}_2$-triple $(e_\al , h_\al , f_\al)$ of $(E,T)$ such that
$$
    (\be | \al) = \ts \frac{(\al|\al)}{2} \, \be(h_\al)
$$
holds for all $\be \in R$, and \sm

\iitem{{\rm (AR3)}} for every $\al \in R^0= \{ \al \in R : (\al|\al)
= 0\}$ there exists a triple $(e_\al , t_\al , f_\al) \in E_\al
\times T \times E_{-\al}$ such that
$$
     (\be | \al) = \be(t_\al)
$$
holds for all $\be \in R$. \sm

%\noindent
Then, with respect to the orthogonal reflections $s_\al$, $\al \in
R\an$, of\/ {\rm \ref{resexamd}}, $(R,X)$ is an affine reflection
system with $R\re=R\an \subset R\int$, $R\im = R^0$ and the bilinear
form $\inpr$ of {\rm (AR1)} as affine form. In particular,
$(E,T,R\an)$ is an affine reflection Lie algebra.} \ms

It follows that there exists an affine form which is
positive-semidefinite on the rational span of $R$ and hence on the
real span of $R$ if $\RR \subset \KK$ (Kac's conjecture). We also
note that the triple $(e_\al, t_\al, f_\al)$ in {\rm (AR3)} spans a
$2$-dimensional abelian or a $3$-dimensional Heisenberg algebra,
depending on $t_\al = 0$ or $t_\al \ne 0$.
\es%\endsubsec

\subsec[invax]{\bf Invariant affine reflection algebras.} Recall
that any bilinear form $\inpr$ on a vector space $V$ gives rise to a
linear map $\flat \colon\, V \to V^* : v \mapsto v^\flat$, where
$v^\flat$ is defined by $v^\flat(u) = (v|u)$  for $v,u\in V$. The
map $\flat$ is injective if and only if $\inpr$ is nondegenerate. In
this case, $\flat$ is an isomorphism if and only if $V$ is
finite-dimensional. \ms

%\noindent
We will call $(E,T)$ an {\it invariant affine reflection
algebra\/},\index{invariant!affine reflection algebra} or
IARA\index{IARA|see{invariant affine reflection algebra}} for short,
if the axioms {\rm (IA1)}--{\rm (IA3)} below hold. \sm
\begin{list}{}{\setlength{\leftmargin}{\parindent}}
\item[\rm (IA1)] $E$ has an invariant nondegenerate symmetric
bilinear form $\inpr$ such that

 \item[(i)] the restriction $\inpr_T$ of the form to $T \times T$
is nondegenerate and
 \item[(ii)] the set of roots $R$ of $(E,T)$ is contained in the
image $T^\flat$ of $T$ with respect to the form $\inpr_T$.
\end{list} \sm

Suppose $(E, T)$ satisfies {\rm (IA1)}. Then for every $\al \in R$
there exists a unique $t_\al \in T$ such that $(t_\al)^\flat = \al$,
i.e., $(t_\al \mid t)_T = \al(t)$ for all $\al \in R$ and $t\in T$.
We use $\flat\colon\, T \to T^*$ to transport the bilinear form
$\inpr_T$ to a symmetric bilinear form $\inpr_X$ on
$X=\Span_\KK(R)$. Thus, by definition $( \al \mid \beta)_X = (t_\al
\mid t_\beta)_T = \al(t_\be)$ for $\al,\beta \in R$. We define the
{\it null\/}\index{null roots}\index{root!null} and {\it
anisotropic\/}\index{anisotropic roots}\index{root!anisotropic}
roots as
$$
  R^0 = \{ \al \in R :
(\al |\al)_X = 0 \}\quad\hbox{and} \quad R\an = \{ \al \in R : (\al
|\al)_X \ne 0 \}.
   \index{R@$R^0$ (null roots)}\index{R@$R\an$ (anisotropic roots)}
  \eq{(1)}%invax3}
$$ We can now introduce the
remaining axioms {\rm (IA2)} and {\rm (IA3)}:  \sm

\iitem{{\rm (IA2)}} For any $0\ne \al \in R$ there exists $e_{\pm
\al} \in E_{\pm \al}$ such that $0\ne [e_\al,\, e_{-\al}] \in T$,
and \sm

\iitem{{\rm (IA3)}} $\ad x_\al $ is locally nilpotent for any $\al
\in R\an$ and $x_\al \in E_\al$. \ms

{\bf Remarks.} (a) Of particular interest are invariant affine
reflection algebras, whose toral subalgebras are finite-dimensional
splitting Cartan subalgebras. They have the following
characterization: A toral pair $(E,T)$ is an invariant affine
reflection algebra with $T$ a finite-dimensional splitting Cartan
subalgebra if and only if \sm

\iitem{$\bull$} $E$ has an invariant nondegenerate symmetric
bilinear form,

\iitem{$\bull$} $T=E_0$ is finite-dimensional, and

\iitem{$\bull$} {\rm (IA3)} holds. \sm

%\noindent
(b) We point out that an IARA may well have an empty set
 of anisotropic roots. For example, the extended Heisenberg algebra,
constructed in \cite[2.9]{kac} and denoted $\frg(0)$ there, is an
IARA with respect to the subalgebra $\frh$ of loc. cit. with $R\an =
\emptyset$.
\endsubsec

\subsec[invth]{\Theorem}{\it Suppose $(E,T)$ is an invariant affine
reflection algebra. Then the conditions {\rm (AR1)}--{\rm (AR3)}
of\/ {\rm \ref{affrefth}} are satisfied for the form $\inpr_X$
defined in\/ {\rm \ref{invax}}, and $(R,X)$ is a symmetric affine
reflection system with $R\re=R\an= R\int$ and $R^0=R\im$. Moreover,
the following hold. \sm

{\rm (a)} If $(E,T)$ is tame, then $R$ is tame.\sm

{\rm (b)} The centreless core $L= E_{cc} = K/ \rmZ(K)$ is a
predivision-$(S,\La)$-graded Lie algebra, which is invariant with
respect to the form induced from $\inpr\mid K\times K$ by the
canonical epimorphism $K \to L$.} \ms

This result and the examples below indicate that the bilinear form
$\inpr$ on an invariant affine reflection algebra $E$ plays only a
secondary role. We therefore define an isomorphism of IARAs as an
isomorphism of the underlying affine reflection algebras, i.e., two
IARAs $(E,T)$ and $(E',T')$ are {\it
isomorphic\/}\index{isomorphic!invariant reflection
algebras}\index{invariant!affine reflection algebra!isomorphic} if
there exists a Lie algebra isomorphism $\ph\colon\, E \to E'$ such
that $\ph(T)=T'$. It then follows that $\ph_r(R\an) = R'{}\an$. \es

\subsec[invcon]{\bf A construction of invariant affine reflection
algebras.} We have seen in \ref{invth} that the centreless core of
an invariant affine reflection algebra is an invariant
predivision-$(S,\La)$-graded Lie algebra with $\La$ a torsion-free
abelian group. We will now describe a construction of invariant
affine reflection algebras which starts from any invariant
predivision-root-graded Lie algebra. In particular, it will follow
from \ref{invcop} that
$$\vcenter{
   \amen{0.78}{every invariant
predivision-$(S,\La)$-graded Lie algebra $L$ with $S$ a locally
finite root system and $\La$ a tor\-sion-free abelian group arises
as the centreless core of an invariant affine reflection algebra.}
 } \eq{(1)}%invcon1}
$$

\noindent  Our construction uses data $(L, D, T_D, C, T_C, \ta)$ as
described below in \INV a -- \INV f.

\begin{description}
\iiitem{\INV a} $L= \bigoplus_{\la \in \La, \, \xi \in S} L_\xi^\la$
is an invariant predivision-$(S,\La)$-graded Lie algebra with $\La$
torsion-free. \end{description}

\noindent As explained in \ref{afrsex}, it is no harm to assume that
$\supp_{\scq(S)} L = S$ and that $\La$ is spanned by $\supp_\La L =
\bigcup_{\xi \in S} \La_\xi$ where $\La_\xi = \{ \la \in \La :
L_\xi^\la \ne 0 \}$. Since $L$ is predivision-graded, we can choose
families $(e_\xi^\la, \, h_\xi^\la, \, f_\xi^\la)\in L_\xi^\la
\times L_0^0 \times L_{-\xi}^{-\la}$, $\xi \in S$ and $\la \in
\La_\xi$, such that $e_\xi^\la$ for $\xi \ne 0$ is invertible with
inverse $f_\xi^\la$, in particular $(e_\xi^\la, h_\xi^\la,
f_\xi^\la)$ is an integrable ${\mathfrak s}{\mathfrak l}_2$-triple,
and such that $(e_0^\la | f_0^\la)_L \ne 0$ but $[f_0^\la, e_0^\la]=
0= h_0^\la$ . We put
$$
   T_L = \frh = \Span_\KK \{ h_\xi^0 : \xi \in S\ind\ti \}.
$$
We denote by $\grSCDer L= \bigoplus_{\la \in \La} (\SCDer L)^\la$
the $\La$-graded subalgebra of centroidal derivations which are
skew-symmetric with respect to $\inpr_L$, see \ref{centd}. One knows
that $\scD= \{ \pa_\vth : \vth \in \Hom_\ZZ(\La, \KK)\}$, the space
of $\La$-degree derivations of $L$ (\ref{degd}), is contained in
$(\SCDer L)^0$. As in \ref{afrsex} every $\la \in \La$ gives rise to
a well-defined linear form $\ev_\la \in \scD^*$. We can now describe
the remaining data $(D, T_D, C, T_C, \ta)$.

\begin{description}

 \item{\INV b} $D= \bigoplus_{\la \in \La}D^\la$ is a
$\La$-graded subalgebra of $\grSCDer L$. \sm

\item{\INV c} $T_D \subset D^0 \cap \scD$ is a subspace such that
the restricted evaluation map $\La \to T_D^*$, given by $\la \mapsto
\ev_\la | T_D$, is injective. \sm

\item{\INV d} $C= \bigoplus_{\la\in \La} C^\la$ is a $\La$-graded
subspace of $D\g*$, which is invariant under the contragredient
action of $D$ on $D\g*$ and contains $\Span_\KK \{ \si_D(l_1, l_2) :
l_i \in L\}$, where $$\si_D\colon\, L \times L \to D\g*$$ is the
central $2$-cocycle defined by $\si_D(l_1, l_2)(d) = \big( d(l_1)
\mid l_2\big)_L$. \sm

\item{\INV e} $T_C \subset C^0$ is a subspace such that the
restriction map $T_C \to T_D^*$, $t_C \mapsto t_C|T_D$, is injective
and such that $\si_D(e_\xi^\la, f_\xi^\la) \in T_C$ for all $\xi \in
S$ and $\la \in \La_\xi$.
 \sm

\item{\INV f} $\ta \colon\, D\times D \to C$ is an {\it invariant
toral\/ $2$-cocycle\/},\index{invariant!toral
$2$-cocycle}\index{co@$2$-cocyle!invariant toral} i.e., $\ta$ is a
bilinear map $\ta \colon\, D \times D \to C$ such that for $d,
d_1,d_2,d_3\in D$
\begin{eqnarray*} %$$\eqalignno{
 \tau (d,d)&=&0, \\
% \tau ([d_1\, d_2],d_3)
%    +\tau ([d_2\, d_3],d_1)+\tau([d_3\, d_1],d_2)
% &= d_1 \cdot \tau (d_2,d_3)
%   + d_2 \cdot \tau (d_3, d_1) + d_3 \cdot \ta(d_1, d_2), \cr
   \ts \sum_{\circlearrowright}\tau ([d_1\, d_2],d_3)
    &=& \ts\sum_{\circlearrowright} d_1 \cdot \tau (d_2,d_3), \\
 \ta(d_1, d_2)(d_3) &=& \ta(d_2, d_3)(d_1) \quad\hbox{and} \\
 \ta(T_D, D) &=& 0. \end{eqnarray*}%$$
Here $\sum_{\circlearrowright}$ indicates the sum over all cyclic
permutations of $(1,2,3)$.  \end{description} \ms

Assume now that $(L,T_L, D, T_D, C, T_C, \ta)$ satisfy \INV a --
\INV f. Then $E=C \oplus L \oplus D$ becomes a Lie algebra with
respect to the product \begin{eqnarray*} %$$\eqalign{
  [c_1 \oplus l_1 \oplus d_1 \, , \, c_2\oplus l_2 \oplus d_2]
        &=&
    \big( \, \si_D(l_1, l_2) + d_1 \cdot c_2 - d_2 \cdot c_1 +
                    \ta(d_1, d_2) \, \big) \\
   &\quad & \oplus \big(\,  [l_1, \, l_2]_L + d_1 \cdot l_2 - d_2 \cdot l_1\, \big)
       \oplus [d_1, \, d_2]_D \end{eqnarray*}
 %} \eq{invcon1} $$
where $c_i \in C$, $l_i \in L$, $d_i \in D$, and
$[.,.]_L$ and $[.,.]_D$ denote the Lie algebra product of $L$ and
$D$ respectively. Moreover, $T= T_C \oplus \frh \oplus T_D $ is a
toral subalgebra of $E$ such that \INV{1.i} holds for $(E,T)$ with
respect to the bilinear form $\inpr_E$ given by
 $$ (\,  c_1 \oplus l_1 \oplus d_1 \mid c_2 \oplus l_2 \oplus d_2 \, )_E =
            c_1(d_2) + c_2 (d_1)  + (l_1 \mid l_2)_L.
             \eq{(2)}%invcon3}
$$
To describe the roots of $(E,T)$, recall that $S$ uniquely embeds
into $\frh^*$. For $\xi \in S$ and $\la \in \La_\xi$ we define a
linear form $\xi \oplus \la \in T^*$ by
$$
   (\xi \oplus \la)\, (t_C \oplus h \oplus t_D) = \xi(h) + \ev_\la (t_D)
         %\eq{invcon2}
$$
($t_C\in T_C$, $h\in \frh$ and $t_D \in T_D$). The root spaces of
$(E,T)$ then are
$$
   E_{\xi \oplus \la}
    = \begin{cases}
      C^\la \oplus L_0^\la \oplus D^\la, & \xi=0, \la \in \La_0, \\
       L_\xi^\la , &\xi \ne 0, \la \in \La_\xi.\end{cases}
           %\eq{invcon5}
$$
Hence the set of roots of $(E,T)$ is $R = \{\xi \oplus \la  \in T^*
: \xi \in S, \la \in \La_\xi\}$. \ms

{\bf Examples.} Suppose $L$ satisfies \INV a. Of course, $D=
\grSCDer L$ or $D=\scD$ are always possible choices for $D$. By
\ref{centd}, $\grSCDer L$ is the semidirect product of the abelian
subalgebra $(\CDer L)^0$ and the ideal $\bigoplus_{\la\ne 0} (\SCDer
L)^\la$, which indicates that there are lots of possibilities for
$D$. An example for $C$ is $C_{\min} = \Span_\KK \{ \si_D(l_1, l_2)
: l_i \in L\}$. In general, any $D$-submodule between $C_{\min}$ and
$C_{\rm max} = (\grSCDer L)\g*$ is a possible choice for $C$. \sm

The conditions on $T_D$ and $T_C$ are interrelated. Since $\KK$ has
characteristic $0$ and $\La$ is torsion-free, $T_D=\scD$ fulfills
the conditions in \INV{c}, see \ref{degl}. Then $T_{C, {\min}} =
\Span_\KK \{ \si_D(e_\xi^\la, f_\xi^\la) : \xi \in S, \la \in
\La_\xi\}$ satisfies \INV{e}, and hence any bigger space does too.
\sm

For $\ta$, we can of course always take $\ta=0$. But there are in
general many more interesting and non-trivial cocycles. For example,
this already happens in the case where $D=\grSCDer L$ is a
generalized Witt algebra and $C=D\g*$, see \cite{Rao-Moody}.\sm
% or\cite{LaN}. \sm

In particular, the choices  $(D,T_D, C, T_C) = (\scD, \scD,
C_{\min}, C_{\min})$ and $(D,T_D, C, T_C) \allowbreak = (\grSCDer L,
\scD, \allowbreak (\grSCDer L)\g*, \allowbreak \scD\g*)$ fulfill all
conditions \INV{b}--\INV{e}. \endsubsec

\subsec[invcop]{\Theorem}{\it The pair $(E,T)$ constructed in\/ {\rm
\ref{invcon}} is an invariant affine reflection algebra with respect
to the form\/ {\rm (\ref{invcon}.2)}. Its centreless core is $L$.
Moreover: \sm

{\rm (a)} The core $E_c$ of $E$ is always contained in $C_{\min}
\oplus L$, and $(E,T)$ is tame if and only if $E_c = C \oplus L$,
which is in turn equivalent to $C\oplus L$ being perfect. \sm

{\rm (b)} If $D=T_D=\scD$ and $C=T_C=C_{\min}$ then $(E,T)$ is
tame.\sm

{\rm (c)} The root system $R$ of $(E,T)$ is always tame.\sm

{\rm (d)} $E_0=T$ if and only if $L$ is an invariant Lie torus,
$T_D=D^0$ and $T_C = C^0$.} %\endsubsec%
\es

\subsec[ealadef]{\bf Extended affine Lie algebras.} These are
special types of invariant affine reflection algebras and are
defined in the introduction. However, since the goal here is a
description of all extended affine Lie algebras in terms of their
centreless cores (\ref{ealamt}), we will only consider tame extended
affine Lie algebras in this section. To simplify the presentation,
we will include tameness in the definition of an EALA. Hence, {\it
an EALA here is the same as a tame EALA in the introduction.} \sm

For the convenience of the reader, who is solely interested in
extended affine Lie algebras, and also for easier comparison with
\cite{n:eala} and other papers on extended affine Lie algebras, we
list the complete set of axioms, following \cite{n:eala}. Moreover,
following the tradition in EALA theory, we will denote the toral
subalgebra of $E$ by $H$ (and not by $T$). This is ``justified''
since (1) shows that $H$ is a splitting Cartan subalgebra of $E$.
\sm

And now, without any further ado, an {\it extended affine Lie
algebra\/}\index{extended affine Lie algebra}\index{Lie
algebra!extended affine}\index{affine Lie algebra!extended} or
EALA\index{EALA|see{extended affine Lie algebra}} for short, is a
pair $(E,H)$ consisting of a Lie algebra $E$ over a field $\KK$ of
characteristic $0$ and a toral subalgebra $H$ satisfying the
following axioms \EA1--\EA6. We denote by $R\subset H^*$ the set of
roots of $(E,H)$, hence $E$ has a root space decomposition
$E=\bigoplus_{\al\in R} E_\al$ with respect to $H$.

\begin{description}
\item{\EA1} $E$ has an invariant nondegenerate symmetric bilinear
form $\inpr$. \sm

\iiitem{\EA2} $H$ is nontrivial, finite-dimensional and
self-centralizing subalgebra
$$E_0 = H.\eq{(1)}%ealadef1}
$$
 \end{description}

Before we can state the other axioms, we need to draw some
consequences of \EA1 and \EA2. The condition (1) implies that the
restriction of the form $\inpr$ to $H\times H$ is nondegenerate.
Since $H$ is finite-dimensional, \INV1 is satisfied. As in
\ref{invax} we can therefore define $t_\al \in H$ for all $\al\in
H^*$, transport the restricted form $\inpr \mid H \times H$ to a
symmetric bilinear form $\inpr_X$ on $X=H^*$ and define isotropic
roots $R^0 = \{ \al \in R : (\al|\al)_X = 0 \}$ and anisotropic
roots $R\an = \{ \al \in R : (\al|\al)_X \ne 0\}$. We let $K=E_c$ be
the core of $(E,H)$, i.e., the subalgebra of $E$ generated by
$\bigcup \{ E_\al : \al \in R\an\}$. We can now state the remaining
axioms.

\begin{description}
\item{\EA3} For every $\al \in R\an$ and $x_\al \in E_\al$, the
operator $\ad x_\al$ is locally nilpotent on $E$. \sm

\item{\EA4} $R\an$ is connected in the sense that for any
decomposition $R\an = R_1 \cup R_2$ with $(R_1 \mid R_2)= 0$ we have
$R_1 = \emptyset$ or $R_2 = \emptyset$. \sm

\item{\EA5} $(E,H)$ is tame.\sm

\item{\EA6} The subgroup $\La = \lan R^0\ran \subset H^*$ generated
by $R^0$ in $(H^*,+)$ is a free abelian group of finite rank.
\end{description} \sm

Suppose the pair $(E,H)$ satisfies \EA1--\EA3. Then \INV1 holds, and
because of \EA3 also \INV2 and \INV3 hold. It therefore follows that
$$
 \amen{0.75}{an EALA is an IARA $(E,H)$ with the following $4$
additional properties: $H$ is a finite-dimensional splitting Cartan
subalgebra and \EA4--\EA6 hold.}$$ In particular, we will define an
{\it isomorphism}\index{extended affine Lie
algebra!isomorphic}\index{isomorphic!extended affine Lie algebras}
of EALAs as an isomorphism of their underlying IARAs. \sm

Since any pair satisfying \EA1--\EA3 is an IARA, Th.~\ref{invth}
implies that $(R,X)$ is an affine reflection system. It therefore
follows from \ref{ext} that \EA4 is equivalent to
\begin{description}

\item{\EA{4$'$}} The quotient root system $(S,Y)$ of the affine
reflection system $(R,X)$ is an irreducible finite root system.
\end{description}

\noindent  Without axiom \EA6,  $\La$ is always a torsion-free
abelian group. Because a torsion-free finitely generated abelian
group is free of finite rank, the axiom \EA6 is equivalent to $\La$
being finitely generated. As we mentioned earlier, $\La$ is an
invariant of an EALA whose rank is called the nullity of $E$. \ms

Since an EALA $(E,H)$ is in particular an affine reflection algebra,
the root system of $(E,H)$ and its core are described in
\ref{arlacore} and \ref{arlath}. Because of (1), the core and the
centreless core are invariant Lie-$(S,\La)$-tori. \ms

We will next describe some properties of this class of Lie tori,
without assuming that they are the centreless core of an EALA. This
is important since in \ref{ealacon} we will specialize the
construction \ref{invcon} to obtain EALAs starting with any
invariant Lie-$(S,\La)$-torus. Because of the special properties of
Lie tori we will in fact obtain a greatly simplified version of
\ref{invcon}. %
\endsubsec%\es

\subsec[lietor]{\bf Some properties of invariant Lie tori.} Let $L$
be an invariant Lie-$(S,\La)$-torus where, as in \ref{ealadef}, $S$
is a finite irreducible root system and $\La$ is a free abelian
group of finite rank. Thus $L$ has a $(\scq(S), \La)$-grading
$$L = \ts \bigoplus_{\xi
\in  S , \, \la \in \La}\, L_\xi^\la. \eq{(1)}%ealasch6}
$$

(i) (\ref{perschn0}) The decomposition $L= \bigoplus_{\xi \in
S}L_\xi$, where $L_\xi = \bigoplus_{\la \in \La} L_\xi^\la$, is the
root space decomposition with respect to the toral subalgebra
$\frh_L = \Span_\KK \{ h_\al : \al \in S\ind\ti\}$. We have $L_0^0 =
\frh_L$. \ms

(ii) $L$ is finitely generated as Lie algebra and has bounded
dimension with respect to the $(\scq(S), \La)$-grading (1), i.e.,
there exists a constant $C$ such that $\dim_\KK L_\xi^\la \le C$ for
all $\xi \in S$ and  $\la \in \La$ (this is of course only of
interest for $\xi=0$).

Let $t$ be the maximum cardinality of minimal generating sets of the
finitely generated reflection spaces $\La_\xi$, $\xi \in S\ind$.
Then $L$ can be generated by $\le 2t\rank S$ elements and $\dim_\KK
L_0^\la \le 4 (t \rank S)^2$. \ms

(iii) Because of (ii), $\Der_\KK L$ is $(\scq(S),\La)$-graded:
$\Der_\KK L = \bigoplus_{\xi \in S, \, \la \in \La} \, (\Der
L)_\xi^\la$, where $(\Der_\KK L)^\la_\xi$ consists of those
derivations mapping $L_\ta^\mu$ to $L_{\ta + \xi}^{\mu + \la}$.
Moreover, the dimension of the homogeneous subspaces $(\Der_\KK
L)_\xi^\la$ is bounded. The analogous statements also hold for
$\SCDer L$. \ms

(iv) (\ref{rocent}(c)) The $\La$-graded algebra $L$ is graded-simple
and prime. The centroid of $L$ is an integral domain, and $L$ embeds
into the $\tl C$-algebra $L \ot \tl C= \tl L$ (tensor product over
$\cent_\KK(L)$), where $\tl C$ be the quotient field of
$\cent_\KK(L)$. \ms

(v) (\ref{rocent}(b)) $\Ga = \supp_\La \cent_\KK(L)$ is a subgroup
of $\La$, and $\cent_\KK(L)$ is isomorphic to the group algebra
$\KK[\Ga]$ over $\KK$, hence to a Laurent polynomial ring in $m$
variables, $0\le m \le \rank \La$ (all possibilities occur). Because
of $(\cent_\KK L)_0\cong \KK \cdot \Id$, it follows that $(\SCDer
L)^0 = \scD$.  \ms

(vi) (\ref{rocent}(b)) $L$ is free as $\cent_\KK L$-module; the
following are equivalent:

\begin{description}
\item{(a)} $L$ is a finite rank $\cent_k(L)$-module,

\item{(b)} $L$ is fgc,

\item{(c)} the number $[\La_\al : \Ga]$ of cosets of $\La_\al$ in
$\Ga$ is finite for some (and hence for every) short root $\al$,

\item{(d)} $[\lan \supp_\La L \ran  : \Ga] < \infty$. \end{description}

\noindent In this case, the Lie algebra $\tl L$ is
finite-dimensional and central-simple over $\tl C$. \ms

(vii) $L$ is fgc, if $S$ is not of type $\rmfa$. \sm

(viii) $\Der L$ is the semidirect product of the ideal $\IDer L$ of
all inner derivations and the subalgebra $\CDer L$ of centroidal
derivations: $\Der L = \IDer L \rtimes \CDer L$. Consequently, ${\rm
SDer} L = \IDer L \rtimes \SCDer L$, where ${\rm SDer} L$ denotes
the subalgebra of skew-symmetric derivations of $L$.
%\endsubsec%
\es

\subsec[ealacon]{\bf A construction of EALAs.} We specialize the
construction and the data $(L,D,T_D, C, T_C)$ of \ref{invcon} with
the aim of obtaining an EALA $(E,H)$. As explained in \ref{invcop},
$L$ must therefore be an invariant Lie-$(S,\La)$-torus where $S$ is
a finite irreducible root system and $\La$ is a free abelian group
of finite rank. Also, since $H$ is a splitting Cartan subalgebra, we
must take $T_D = D^0$ and $T_C = C^0$, see \ref{invcop}. By
\ref{lietor}, we have $\grSCDer L = \SCDer L$ with $(\SCDer L)^0 =
\scD$. Moreover $\SCDer L$ and hence any $\La$-graded subalgebra
have finite homogenous dimension. Since $C$ and $D$ are
nondegenerately paired, it follows that $C=D\g*$. Therefore,
specializing \ref{invcon}, we obtain the following construction,
based on data $(L, D, \ta)$, where \begin{description}

\iiitem{\EA{a}} $L$ is an invariant Lie-$(S,\La)$-torus with
$(S,\La)$ as above, \sm

\item{\EA{b}} $D$ is a graded subalgebra of $\SCDer L$ such that
the canonical map $\La \to (D^0)\g*$, $\la \mapsto \ev_\la|D^0$, is
injective, and \sm

\iiitem{\EA{c}} $\ta\colon\, D \times D \to D\g*$ is an invariant
toral $2$-cocycle, as in \INV{f}. \end{description}\sm

\noindent We denote the resulting invariant affine reflection
algebra by $E(L,D,\ta)$. By construction, $H=(D^0)^* \oplus \frh
\oplus D^0$ is a splitting Cartan subalgebra and \EA4 and \EA6 hold.
One can also show that $D\g* \oplus L$ is perfect, so that $(E,H)$
is tame. Thus we have the first part of the following theorem,
describing the structure of EALAs in terms of Lie tori.
\endsubsec

\subsec[ealamt]{\Theorem}{\it The invariant affine reflection
algebra $E(L,D,\ta)$ constructed in\/ {\rm \ref{ealacon}} is an
EALA. Conversely, let $E$ be an EALA  and let $L$ be its centreless
core. Then there exist unique data $(D,\ta)$ as in\/ {\rm
\ref{ealacon}} such that $E=E(L,D,\ta)$. } \ms

Since $D=\scD$ or $D=\SCDer L$ and $\ta=0$ always satisfy the
conditions \EA{b} and \EA{c}, any invariant Lie torus arises as the
centreless core of an EALA. \sm

{\bf Example.} If $L=\frg \ot_\KK \KK[t^{\pm 1}]$ is the untwisted
loop algebra as in \ref{afrreal} and \ref{fiex}(b), then $\grSCDer L
= \SCDer L = \KK\, d^{(0)}$. Hence $D=\KK\,d^{(0)}$ is the only
possible choice for $D$ and then necessarily $\ta=0$. Thus, in this
case $E(L,D,\ta)$ is the Lie algebra $E$ constructed in
\ref{afrreal}(III), and this is the only extended affine Lie algebra
based on $L=\frg\ot \KK[t^{\pm 1}]$. \sm

In the remainder of this section we will discuss other classes of
invariant affine reflection algebras. \endsubsec

\subsec[diseala]{\bf Discrete extended affine Lie algebras.} Let $E$
be a Lie algebra over $\KK=\CC$ and let $H$ be a toral subalgebra of
$E$. We call $(E,H)$ a {\it discrete extended affine Lie
algebra\/}\index{extended affine Lie
algebra!discrete}\index{discrete!extended affine Lie
algebra}\index{Lie algebra!discrete extended affine} if $(E,H)$
satisfies the axioms \EA1--\EA5 above and in addition \sm

\begin{description}
 \iiitem{(DE)} $R$ is a discrete subset of $H^*$. \end{description}\sm

\noindent In (DE) ``discrete'' refers to the natural topology of
$H^*$. Suppose $(E,H)$ satisfies \EA1--\EA5. Then tameness of
$(E,H)$, which is axiom (EA5), implies tameness of $R$. Also, the
quotient root system $S$ is irreducible by \EA4. One can easily show
that (DE) holds iff $R^0$ is discrete iff $\La$ is discrete. Since a
discrete subgroup of $H^*$ is always free of finite rank, (DE)
implies \EA6. Hence, a discrete EALA is in particular an EALA, but
not every EALA over $\KK=\CC$ or $\RR$ is discrete, see
\ref{LEALAdef}. \ms

In the construction \ref{ealacon}, $E(L,D,\ta)$ is a discrete EALA
if the evaluation map $\La \to (D^0)^*$ has a discrete image.
Conversely, any discrete EALA arises in this way. For example, the
embedding $\La \to \scD^*$ always has a discrete image. \es

\subsec[stinva]{\bf Reductive Lie algebras.} Let $E$ be a
finite-dimensional reductive Lie algebra over $\KK$, and let $T$ be
a maximal toral subalgebra of the semisimple subalgebra $[E,E]$.
Then $(E,T)$ is an invariant affine reflection algebra of nullity
$0$. \sm

Observe that $E=\rmZ(E) \oplus [E,\, E]$ with $[E,\,E]$ being the
core and the centreless core of the IARA $(E,T)$. By \cite[\S
I.2]{Se:rat}, the set of roots $R$ of $(E,T)$ is a finite, possibly
non-reduced root system. (Indeed, by \ref{invth}, $R$ is a finite
affine reflection system, hence an extension of a finite root system
$S$ by an extension datum $\frL= (\La_\xi)$ where each $\La_\xi$ is
a finite subset of torsion-free group, whence $\La_\xi = \{0\}$ and
$S=R$.) The maximality of $T$ implies that $E_c$ is not only
predivision-graded, but actually division-graded,
\cite[Lemma~I.3]{Se:rat}, but note that $E_c$ is in general not a
Lie torus. Also, $(E,T)$ is tame if and only if $E$ is semisimple.
\sm

Of course, all of this is not very helpful if $T=\{0\}$, i.e., the
semisimple part $[E,E]$ of $E$ is anisotropic, which is allowed. On
the other extreme, if $E$ is simple and $T$ is a splitting Cartan
subalgebra, e.g. $\KK$ is algebraically closed, then $E$ is an EALA
of nullity $0$. Conversely, every EALA of nullity $0$ is a
finite-dimensional split simple Lie algebra.
\endsubsec

\subsec[LEALAdef]{\bf Locally extended affine Lie algebras.} A toral
pair $(E,H)$ with set of roots $R$ is called a {\it locally extended
affine Lie algebra\/}\index{locally extended!affine Lie
algebra}\index{extended affine Lie algebra!locally}\index{Lie
algebra!locally extended affine}\index{affine Lie algebra!locally
extended} or LEALA\index{LEALA|see{locally extended affine Lie
algebra}} for short (\cite{morita-yoshii}), if $(E,H)$ satisfies the
following axioms: \begin{description}

\item[(LEA1)] $H=E_0$, \sm

\item[(LEA2)] $E$ has an invariant nondegenerate symmetric bilinear
form $\inpr$, \sm

\item[(LEA3)] $R\subset H^\flat$, where $\flat \colon\, H \to H^*$ is
the canonical map induced by $\inpr \mid H\times H$, \sm

\item[(LEA4)] $\ad_E x$ is locally nilpotent for $x\in E_\al$ and
$\al \in R\an = \{ \al \in R : (\al | \al) \ne 0\}$.
\end{description}
\sm

%\noindent
Let $(E,H)$ be such a Lie algebra. We point out that $R\an =
\emptyset$ is allowed. Since $\inpr \mid H \times H$ is
nondegenerate by invariance of $\inpr$ and (LEA1), the axiom
\INV{1.i} is fulfilled. Observe \INV{1.ii} = (LEA3) and \INV3 =
(LEA4). Because also \INV2 holds, we have
$$
\amen{0.83}{a LEALA is the same as an invariant affine reflection
algebra $(E,H)$ with a splitting Cartan subalgebra $H$ and a
connected set of anisotropic roots.}
$$
In this case, the core $K$ of $(E,H)$ and hence also the centreless
core $L$ are Lie tori, \ref{arlath}. Moreover, by \ref{invth}, $L$
is an invariant Lie torus. \sm

We note that LEALAs generalize EALAs: A LEALA $(E,H)$ is an EALA if
and only if $H$ is finite-dimensional, $(E,H)$ is tame and $\lan R^0
\ran$ is finitely generated. \ms

{\bf Examples.} (a) By \ref{invcop}(e), the invariant affine
reflection algebras constructed in \ref{invcon} are LEALAs if and
only if $T_D=D^0$, $T_C=C^0$ and $L$ is an invariant Lie torus, say
of type $(S,\La)$ where $S$ is an irreducible locally finite root
system. For example, the IARAs $(L, \scD, \scD, C_{\min}(\scD),
\allowbreak C_{\min}(\scD), \ta)$ are LEALAs. \sm

(b) Let $\frg$ be a finite-dimensional split simple Lie
$\KK$-algebra or a centreless infinite rank affine algebra, see
\ref{afrs}. Let $\La$ be a subgroup of $(\KK, +)$ generated by a set
$B$ of $\QQ$-linearly independent elements, and let $\KK[\La]$ be
the group algebra of $\La$. Then $L= \frg \ot_\KK \KK[\La]$ is a
Lie-$(S,\La)$-torus  which is invariant with respect to the bilinear
form $\ka \ot \eps$, where $\ka$ is a invariant nondegenerate
symmetric bilinear form on $\frg$ and $\eps$ is the canonical trace
form of $\KK[\La]$. Let $\vth \colon\, \La \to \KK$ be given by
$\vth(\la) = \la$ and denote by $\pa_\vth$ the corresponding degree
derivation. Then $D=\KK\pa_\vth = T_D$ and $C=D^*=T_C$ satisfy
\INV{b}--\INV{e}, and hence the construction \ref{invcon} produces a
LEALA $E$, e.g. by taking $\ta=0$. It is tame.

Note that for this example the affine root system is $R=\{ \xi
\oplus \la \in \frh^* \oplus (\KK \pa_\vth)^* : \xi \in S, \la \in
\La\}$ with isotropic roots $R^0 \cong \La \subset (\KK
\pa_\vth)^*$. Hence, the nullity of $(E,H)$ is the cardinality
$|B|$, while $\dim_\KK \Span_\KK (R\im)= 1$. In particular, this
shows that \ref{ARLA}(2) is in general a strict inequality.

It also follows that this LEALA is an EALA if and only  if $S$ and
$B$ are finite. Assume this to be the case and let $\KK=\CC$. Then
$\La=R^0$ is a discrete subset of $D \cong \KK$, i.e., $E$ is a
discrete EALA, if and only if $B$ is linearly independent over
$\RR$. For example, the EALA corresponding to $B=\{1, \sqrt{2}\}$ is
not discrete.
\endsubsec

\subsec[graladef]{\bf Generalized reductive Lie
algebras.}\index{generalized reductive!Lie algebra}\index{Lie
algebra!generalized reductive} By definition \cite{az:grla}, these
are the invariant affine reflection algebras $(E,H)$ over $\KK=\CC$
with the following properties: \sm

(GRLA1) $H$ is a non-trivial finite-dimensional splitting Cartan
subalgebra, \sm

(GRLA2) the set of roots $R$ of $(E,H)$ is a discrete subset of
$H^*$. \sm

\noindent We point out that $R\an$ is not necessarily connected,
i.e., the quotient root system of the corresponding affine
reflection system is not necessarily irreducible. It follows from
Th.~\ref{arlath}(c), \ref{invth} and the axioms (GRLA1),
(GRLA2) that the affine reflection system of such an algebra has
finite rank, unbroken root strings and is symmetric, discrete and
reduced. In other words, $R$ is a generalized reductive root system
in the sense of \ref{affex}. \sm

As in \ref{LEALAdef}, the core and centreless core of such an
algebra are invariant Lie tori, say of type $(S,\La)$. Here $S$ is a
finite, not necessarily irreducible root system and $\La$ is a free
abelian group of finite rank. \sm

The construction \ref{invcon} provides examples of this type of
IARAs, assuming in \INV{a} that $L$ is an invariant
Lie-$(S,\La)$-torus  with $S$ a finite root system and $\La\cong
\ZZ^n$ and in \INV{c} and \INV{e} that $T_D=D^0$ and $T_C= C^0$.
\es%\endsubsec

\subsec[stinvd]{\bf Toral type extended affine Lie
algebras.}\index{toral type extended affine Lie algebra}\index{Lie
algebra!toral type extended}\index{extended affine Lie algebra!toral
type}\index{affine Lie algebra!toral type extended} By definition
\cite{AKY}, these are invariant affine reflection algebras $(E,T)$
which also satisfy the conditions (TEA1)--(TEA4): \sm

(TEA1) $T$ is non-trivial and finite-dimensional, \sm

(TEA2) the affine root system $R$ of $(E,T)$ is tame and $R\an$ is
connected, \sm

(TEA3) $(E,T)$ has finite nullity, and \sm

(TEA4) $R$ is a discrete subset of $\Span_\QQ(R) \ot \RR$. \sm

\noindent It follows from Th.~\ref{arlath}(c), \ref{invth} and the
axioms (TEA2)--(TEA4) that the affine reflection system of such an
algebra has finite rank, a connected ${\rm Re}(R)$, unbroken root
strings and is symmetric, discrete and tame. In other words, $R$ is
an extended affine root system as defined in \ref{affex} with the
exception that $R$ need not be reduced.
\endsubsec%\es

\subsec[KMex]{\bf Kac-Moody algebras.} Concerning Kac-Moody algebras
we will follow the definitions and terminology of \cite{kac}. We
have the following characterization of invariant affine reflection
algebras and extended affine Lie algebras among the Kac-Moody
algebras: \ms

Let $\frg(A)$ be a Kac-Moody algebra for some generalized Cartan
matrix $A$ and let $\frh$ be its standard Cartan subalgebra. Then
the following are equivalent: \sm \begin{description}

\item{(i)}  $(\frg(A), \frh)$ is an invariant affine reflection algebra, \sm

\item{(ii)} $A$ is symmetrizable and every imaginary root is a null root,
\sm

\item{(iii)}  every connected component of the Dynkin diagram of $A$ is
either of finite or of affine type. \end{description}

\noindent In this case, the anisotropic roots of the associated
affine reflection system equals the set of real roots in the
terminology of \cite{kac}, and the core and the centreless core of
$\frg(A)$ are Lie tori. In particular, $(\frg(A),\frh)$ is an EALA
if and only if $A$ is of finite or affine type, in which case the
nullity of $\frg(A)$ is, respectively, $0$ or $1$, and $\frg(A)$ is
discrete. In fact, one has a much stronger result:
\endsubsec

\subsec[null1char]{\Theorem} (\cite{abgp}) {\it A complex Lie
algebra $E$ is (isomorphic to) a discrete extended affine Lie
algebra of nullity $1$ if and only if $E$ is (isomorphic to) an
affine Kac-Moody algebra.} \es

\subsec[ealanotes]{\bf Notes.} The results in this section are due
to the author, unless indicated otherwise in the text or below.
Details will be contained in \cite{n:lect}. Some results have been
independently proven by others, as indicated below. \sm

The importance of isotopy for the theory of extended affine Lie
algebras is explained in \cite{AF:isotopy}. Assuming the Structure
Theorem for extended affine Lie algebras, Th.~\ref{arlath}(d) is in
fact proven in \cite[Th.~6.1]{AF:isotopy} for the case of extended
affine Lie algebras. \sm

The notion of ``nullity'' goes back to the early papers \cite{bgk}
and \cite{bgkn} on EALAs, where it was implicitly assumed that
\ref{ARLA}(2) is an equality. This was later corrected in
\cite{gao:deg}. Beware that in the paper \cite{morita-yoshii} the
term ``null rank'' is used for the nullity as defined here, while
nullity has another meaning.  \sm

Kac's conjecture (\ref{affrefth}) for EALAs was proven in \cite[I,
\S 2]{aabgp}. It was shown for LEALAs in
\cite[Th.~3.10]{morita-yoshii} and for GRLAs in \cite[\S
1]{az:grla}. The proof of \ref{affrefth} generalizes techniques
found in these papers. Th.~\ref{invth}(a) is proven in
\cite[Lemma~3.62]{abp1} for EALAs. \sm

The definition of an extended affine Lie algebra in \ref{ealadef} is
due to the author (\cite{n:eala}). It extends the previous
definition from \cite{aabgp}, which only made sense for Lie algebras
over $\CC$ or $\RR$. Tame extended affine Lie algebras in the sense
of \cite{aabgp} are here called discrete extended affine Lie
algebras (see \ref{ealadef}). To describe the structure of the core
and the centreless core of an EALA $E$, it would be enough to
require that $R$ is tame. We have included tameness of $(E,H)$ in
the definition of an EALA since there seems to be little hope to
classify non-tame EALAs, as the examples in \ref{invcop} hopefully
demonstrate. The reader can find another example at the beginning of
\S  3 of \cite{bgk}. That a tame extended affine Lie algebra has a
tame root system, see \ref{invth}(b), is proven in
\cite[3.62]{abp1}. The structure of the root system and the core of
an EALA is determined in \cite[\S  1]{AG2}. \sm

Some of the properties of Lie tori summarized in \ref{lietor}, the
construction \ref{ealacon} and Th.~\ref{ealamt} were announced in
\cite{n:eala}. The characterization of fgc Lie tori in
\ref{lietor}(vi) is also proven in \cite[Prop.~1.4.2]{abfp2}. The
interest in fgc Lie tori comes from the multiloop realization
described \ref{perweakv}. This description is made more precise in
\cite{abfp2}. In this context it is of importance to know
\ref{lietor}(vii). Part (viii) is crucial for the Structure Theorem
\ref{ealamt}. \sm

Yoshii has shown in \cite{y:lie} that for $S$ a finite root system
every centreless complex Lie-$(S,\ZZ^n)$-torus automatically has a
graded invariant nondegenerate symmetric bilinear form $\inpr_L$
with respect to which $(L, \frh)$ is invariant. Such a form is
unique, up to scalars. Hence, in this case one could leave out the
quantifier ``invariant''. Yoshii's proof uses the classification of
Lie tori of type $(\rmfa_1,\ZZ^n)$ and hence ultimately Zelmanov's
Structure Theorem for strongly prime Jordan algebras
(\cite{zel:primeII}, \cite{mczel}). \sm

The construction \ref{ealacon} with $D=\scD$ and $\ta=0$ appears in
\cite[III]{aabgp} and a modified (twisted) version of it in
\cite{Az5}. For $S$ of type $\rma_l$, $(l\ge 3)$, $\rmd$ or $\rme$,
the construction with $D=\SCDer L$ and arbitrary $\ta$ can be found
in \cite{bgk}. The fact that the Moody-Rao cocycle provides a
non-trivial example of an invariant toral $2$-cocycle is also
mentioned there. Finally, for $S=\rmfa_2$ the general cosntruction
is given in \cite{bgkn}, modulo the mistake mentioned above. The two
papers \cite{bgk} and \cite{bgkn} rely heavily on a concrete
realization of the corresponding Lie torus $L$, and do not present
the construction in a way suitable for generalization to the case of
an arbitrary $S$. \sm

LEALAs (\ref{LEALAdef}) were introduced in \cite{morita-yoshii}. The
Example(b) of \ref{LEALAdef} appears there. This paper also contains
examples of LEALAs with $R\an=\emptyset$ and a classification of
LEALAs of nullity $0$. \sm

It is shown in \cite{az:grla} that the set of roots of a generalized
reductive Lie algebra is a generalized root system and that its core
and centreless core are Lie tori. A proof that the set of roots of a
toral type extended affine Lie algebra is an extended affine root
system is given in \cite{AKY}. This paper also gives an example of a
toral type EALA of type $\rmbc$ with a non-reduced set of roots.
\es%\endsubsec

\vfe

\section{Example: $\frslI (A)$ for $A$ associative} \label{perexass}

We describe the various concepts introduced in the previous sections
for the special linear Lie algebra $\frslI (A)$ and its central
quotient $\frpsl_I(A)$, where $A$ is an associative algebra.\sm

Unless specified otherwise, $k$ is a commutative associative ring
with unit element $1_k$ such that $2\cdot 1_k$ and $3\cdot 1_k$ are
invertible in $k$. All algebraic structures are defined over $k$.
Throughout, $I$ is a not necessarily finite set of cardinality
$|I|\ge 3$ and $\La$ denotes an abelian group, written additively.
\ms

\subsec[assba]{\bf Basic definitions.} Let $A$ be a unital
associative $k$-algebra with identity element $1_A$. We recall that
$[a,b]:= ab-ba$ for $a,b\in A$ denotes the {\it commutator\/} in
$A$. We put $[A,A]= \Span_k\{ [a,b] : a,b\in A\}$ and denote by
$\rmZ(A) = \{ z\in A : [z,a]=0\}$ the centre of $A$. The underlying
$k$-module of $A$ becomes a Lie algebra with respect to the
commutator product $[.,.]$, denoted $A^-$. \sm

We denote by $\Mat_I(A)$\index{M@$\Mat$ (finitary matrices)} the set
of all finitary $I\times I$-matrices with entries from $A$. By
definition, an element of $\Mat_I(A)$ is a square matrix
$x=(x_{ij})_{i,j \in I}$ with all $x_{ij}\in A$ and $x_{ij} \ne 0$
for only finitely many indices $(ij)\in I\times I$. With respect to
the usual addition and multiplication of matrices, $\Mat_I(A)$ is an
associative $k$-algebra. For $i,j\in I$ we let $E_{ij}\in
\Mat_I(A)$\index{E@$E_{ij}$ (matrix unit)} be the matrix with entry
$1_A$ at the position $(ij)$ and $0$ elsewhere. Then any $x\in
\Mat_I(A)$ can be uniquely written as $x=\sum_{i,j} \, x_{ij}
E_{ij}$, $x_{ij} \in A$, with almost all $x_{ij} = 0$. \sm

We denote by $\frgl_I(A)$\index{G@$\frgl$ (general linear Lie
algebra)} the Lie algebra $\Mat_I(A)^-$ and by $\frslI(A)=
[\frgl_I(A), \, \frgl_I(A)]$\index{S@$\frslI(A)$ (special linear Lie
Lie algebra)} its derived algebra. One easily verifies that
$\frslI(A)$ is the subalgebra of $\frgl_I(A)$ generated by all
off-diagonal matrices $AE_{ij}$, $i\ne j$ and that $\frslI(A)$ can
also be described as
$$
  \frslI(A) = \{ x\in \frgl_I(A) : \tr(x)\in [A,A]\} \eq{(1)}%assba1}
$$
where $\tr(x) = \sum_i x_{ii}$\index{T@$\tr$ (trace)} is the {\it
trace\/} of $x$. Since $\frslI(A) \cong {\mathfrak s}{\mathfrak
l}_J(A)$ if $|I|=|J|$ we put
$$
   {\mathfrak s}{\mathfrak l}_n(A) = \frslI(A) \quad \hbox{if $|I| =n\in \NN$.}
$$
It is important that the algebra $A$ can be recovered from the Lie
algebra $\frslI(A)$ by the product formula:
$$
  ab\, E_{ij} = \big[ \big[ [ aE_{ij}, \, E_{jl}], \, E_{li}\big],
    \, bE_{ij}\big]
   \eq{(2)}%assba4}
$$
where $a,b\in A$ and $i,j,l\in I$ are distinct. \sm

What makes $\frslI(A)$ somewhat delicate, are the diagonals of its
elements. For a commutative $A$ the description (1) of
$\frslI(A)$ shows that $\frslI(A)$ is the ``usual'' Lie algebra:
$$
 \frslI(A) = \{ x\in \frgl_I(A) : \tr(x) = 0 \} \cong
   \frslI(k) \ot_k A \qquad \hbox{($A$ commutative)}. \eq{(3)}%assba1a}
$$
But (3) is no longer true for a non-commutative $A$. Rather, if $n
\cdot 1\in k\ti$ then any $x\in \frsln(A)$ can be uniquely written
in the form $x=aE_n + x'$ with $a\in [A,A]$ and $x'\in \Mat_n(A)$
with $\tr(x') = 0$, namely $x=\frac{1}{n} \tr(x) E_n + (x-
\frac{1}{n}\tr(x)E_n)$. Hence, as $k$-module,
$$
   \frsln(A) = [A,A]E_n \oplus (\frsln(k) \ot_k A)
   \qquad \hbox{($n\cdot 1_k \in k\ti$)} \eq{(4)}%assba5}
$$
where $E_n$ is the $n\times n$-identity matrix. One can describe the
Lie algebra product with respect to this decomposition. To do so,
let us abbreviate $a\circ b = ab+ba$ for $a,b$ in an associative
algebra $B$, like $B=A$ or $B=\Mat_I(A)$ (the experts will recognize
$a\circ b$ as the ``circle product'' of the Jordan algebra
associated to $B$). Then, for $x,y\in \frsln(k)$, $a,b\in A$ and
$c,d\in [A,A]$, \begin{eqnarray*}
  [x\ot a, \, y\ot b] &= & (\ts \frac{1}{n} \tr(xy)\, [a,b] \,E_n )
             \oplus \big( \,[x,y] \ot \frac{1}{2} (a\circ b )
      \nonumber \\
  & & \qquad \ts
     + ( \frac{1}{2} (x \circ y) - \frac{1}{n} \tr(xy) E_n
                )  \ot [a,b] \,\big), \\
  {[ cE_m,\, dE_n]}  &=& [c,d]\, E_n, \\
  { [cE_n, \, x\ot a]}  &= & x \ot [ca].
 \end{eqnarray*} %assba3}

%$$\eqalign{
%  [x\ot a, \, y\ot b] &= (\ts \frac{1}{n} \tr(xy)\, [a,b] \,E_n )
%            \oplus \big( \,[x,y] \ot \frac{1}{2} (a\circ b )
%    \cr &\qquad \ts
%      + ( \frac{1}{2} (x \circ y) - \frac{1}{n} \tr(xy) E_n
%                 )  \ot [a,b]\,\big), \cr
%  [cE_n, \, dE_n] &= [c,d]\, E_n , \cr
%    [cE_n, \, x\ot a] &= x \ot [ca].
%} \eq{assba3}
%$$
\endsubsec%\es

\subsec[sligr]{\bf $\frslI(A)$ as $\rma_I$-graded Lie algebra.}
Recall the irreducible root system $ R=\rma_I = \{ \eps_i-\eps_j :
i,j\in I \}$ defined in \ref{lfrsclas}. For $\al= \eps_i-\eps_j\ne
0$ and $\be=\eps_m-\eps_n\in R$ the integers $\lan \be,\al\ch\ran$
are $ \lan \eps_m - \eps_n , (\eps_i - \eps_j)\ch\ran
  = \de_{mi} - \de_{mj} - \de_{ni} + \de_{nj}$, where $\de_{**}$ is the Kronecker delta. \sm

The Lie algebra $\frslI(A)$ has a natural $\scq(\rma_I)$-grading
with support $\rma_I$, whose homogeneous spaces are given by
$$
 \frslI(A)_\al = \begin{cases} AE_{ij}, & \al= \eps_i - \eps_j \ne 0,\\
     \{ x\in \frslI(A) :  x \hbox{ diagonal} \} , & \al= 0.\end{cases}
 \eq{(1)}%sligr1}
$$
Because of our assumption that $2$ and $3$ are invertible in $k$,
the space $\frslI(A)_\al$, $\al = \eps_i - \eps_j \ne 0$, has the
description
$$
  \frslI(A)_\al = \{ x\in \frslI(A) : [E_{ii}-E_{jj},\, x] = 2 x\}
         \eq{(2)}%sligr2}
$$

An element $aE_{ij} \in \frslI(A)_\al$, $\al= \eps_i-\eps_j \ne 0$,
is an invertible element of the $\scq(\rma_I)$-graded Lie algebra
$\frslI(A)$, cf. \ref{rogr}, if and only if $a$ is invertible in
$A$. In this case, $-a^{-1} E_{ji}$ is the inverse of $aE_{ij}$ and
$$
  (a E_{ij}, \, E_{ii}- E_{jj},\, - a^{-1}E_{ji})
   \hbox{ {\it  is an ${\mathfrak s}{\mathfrak l}_2$-triple.}}
$$
In particular, every matrix unit $E_{ij}$, $i\ne j$, is invertible
in $\frslI(A)$.  Since by construction
$$
  \frslI(A)_0 = \tsum_{0\ne \al\in R} \, [\frslI(A)_\al, \,
          \frslI(A)_{-\al}]
$$
it follows that {\it $\frslI(A)$ is an $\rma_I$-graded Lie algebra}
in the sense of \ref{rogr}. Hence, either by \ref{rogrcov} or by an
immediate verification,
$$  \frpsl_I(A) = \frslI(A)/ \rmZ(\frslI(A))
     \index{P@$\frpsl$ (projective special linear Lie algebra)}
$$
and any covering of $\frpsl_I(A)$ are $\rma_I$-graded Lie algebras.
For $\rank \rma_I \ge 3$ these are in fact all $\rma_I$-graded Lie
algebras:
%\endsubsec%

\es

\subsec[typeAcl]{\Theorem}{\it Let $|I|\ge 4$. Then a Lie algebra
$L$ is $\rma_I$-graded if and only if $L$ is a covering of
$\frpsl_I(A)$ for a unital associative algebra $A$. }\ms

In the light of this theorem and the definition of $\frpsl_I(A)$ it
is of interest to identify the centre of $\frslI(A)$, which will be
done in \ref{ccc}, and the coverings of $\frpsl_I(A)$. Since any
covering of a perfect Lie algebra $P$ is a quotient of its universal
central extension $\uce(P)$ by a central ideal, it is enough to
describe the universal central extension $\uce(\frpsl_I(A)) \to
\frpsl_I(A)$. Because $\frslI(A)$ is a covering of $\frpsl_I(A)$,
the Lie algebras $\uce(\frpsl_I(A))$ and $\uce(\frslI(A))$ are the
same (in any case, in many interesting cases $\frslI(A) =
\frpsl_I(A)$, see \ref{ccc} and \ref{grtypeAco}). It is therefore
sufficient to identify $\uce(\frslI(A))$, which we will do in
\ref{uce} after we introduced gradings in \ref{slila}.
\endsubsec

\subsec[ccc]{\bf Centre of $\frslI(A)$.} The centre of $\frslI(A)$
has the following description, cf. (\ref{centex}.1):
$$
 \rmZ\big( \frslI(A)\big) =
   \begin{cases} \{0\}, & \text{if $|I| = \infty$,} \\
    \{ z E_n \in \frslI(A) : z\in \rmZ(A)\}, &
          \text{if $|I|=n\in \NN$.} \end{cases}
 \eq{(1)}%assbac}
$$
In particular $\frslI(A) = \frpsl_I(A)$ for $|I|=\infty$, while
$$
 \rmZ\big( \frsln(A)\big) = \{ z E_n : z\in \rmZ(A), nz\in
    [A,A]\} \quad \hbox{for $n \in \NN$,}
 \eq{(2)}%assbac1}
$$
and therefore also $\frsln(A) \cong  \frpsl_n(A)$ if, for example,
$k$ is a field of characteristic $0$ and $A$ is commutative. \sm

Since $\rmZ(\frslI(A))$ is diagonal, the homogeneous spaces
$\frpsl_I(A)_\al$ are given by the obviously modified formula
(\ref{sligr}.1). In particular, we can identify $AE_{ij} =
\frpsl_I(A)_\al$ for $0\ne \al=\eps_i-\eps_j$ and then also have
(\ref{sligr}.2) and the product formula (\ref{assba}.2) available in
$\frpsl_I(A)$.
%\es%
\endsubsec

\subsec[slila]{\bf Graded coefficient algebras.} Let $A=
\bigoplus_{\la\in \La} A^\la$ be a $\La$-grading of $A$ and define
$\frslI(A)^\la$ as the submodule of matrices in $\frslI(A)$ whose
entries all lie in $A^\la$. Then
 $$\frslI(A) = \ts \bigoplus_{\la\in \La} \frslI(A)^\la \eq{(1)}%slila1}
$$
is a $\La$-grading of $\frslI(A)$ which is compatible with the
$\scq(\rma_I)$-grading (\ref{sligr}.1). The homogeneous spaces
$$\frslI(A)_\al^\la = \frslI(A)_\al \cap \frslI(A)^\la$$ are given by
$\frslI(A)_\al^\la=  A^\la E_{ij}$, $\al=\eps_i-\eps_j \ne 0$, while
$\frslI(A)_0^\la$ consists of the diagonal matrices whose entries
all lie in $A^\la$. Since $1_A\in A^0$, it follows that
$\frslI(A)_\al^0$, $\al\ne 0$, contains an invertible element. Thus
$$\hbox{\it
$\frslI(A) = \bigoplus_{\al\in \rma_I, \, \la\in \La} \,
\frslI(A)^\la_\al$ is an $(\rma_I, \La)$-graded Lie
algebra.}\eq{(2)}%slila2}
$$
It is immediate from the product formula (\ref{assba}.2) that,
conversely, any $\La$-grading of $\frslI(A)$ which is compatible
with the $\rma_I$-grading is obtained from a $\La$-grading of $A$ as
described above. \sm

Since the centre of $A$ and $[A,A]$ are graded submodules, it
follows from (\ref{assba}.2) that $\rmZ(\frslI(A))$ is $\La$-graded.
Hence  $\frpsl_I(A)$ is $(\scq(\rma_I), \La)$-graded with homogenous
spaces \begin{eqnarray*}
    \frpsl_I(A)_\al^\la & \cong& \frslI(A)_\al^\la \quad
         \text{for $0 \ne \al\in \rma_I$ and} \\
 \frpsl_I(A)_0^\la & =& \frslI(A)_0^\la \big/ \rmZ(\frslI(A))^\la.
\end{eqnarray*}
 This is of course a special case of \ref{rogrcov}. The graded version of \ref{typeAcl} is the following Theorem.
%\es%
\endsubsec

\subsec[grtypeAcl]{\Theorem}{\it Let $|I|\ge 4$. Then a Lie algebra
$L$ is $(\rma_I,\La)$-graded if and only if $L$ is a $\La$-covering
of\/ $\frpsl_I(A)$ for some unital associative $\La$-graded algebra
$A$. } \ms

Combining this with the description of invertible elements in
$\frslI(A)$ immediately shows that $\frslI(A)$ and $\frpsl_I(A)$ are
predivision- or division-$(\rma_I,\La)$-graded Lie algebras if and
only if $A$ is, respectively, a predivision- or a
division-$\La$-graded algebra. Similarly, $\frslI(A)$ is a Lie torus
of type $(\rma_I, \La)$ if and only if $A$ is an associative
$\La$-torus. Hence %, by \no{cenrogra},
we get the following corollary (a special case of \ref{rogrcov}).
\endsubsec

\subsec[grtypeAco]{\Corollary}{\it Let $|I|\ge 4$. Then $L$ is a
predivision- or a division-$(\rma_I,\La)$-graded Lie algebra if and
only if $L$ is a $\La$-covering of\/ $\frpsl_I(A)$ where $A$ is,
respectively, a predivision- or a division-$\La$-graded algebra. If
$k$ is a field then $L$ is a Lie torus of type $(\rma_I,\La)$ if and
only if $L$ is a $\La$-covering of $\frpsl_I(A)$ where $A$ is a
$\La$-torus. } \ms

Refer to \ref{assogr} for a discussion of (pre)division-graded
associative algebras. Concerning the structure of Lie tori, we note
that
$$
   A = \rmZ(A) \oplus [A,\, A] \eq{(1)}%grtypeAco1}  \quad
$$
for every associative $\La$-torus $A$. Hence, by (\ref{ccc}.1) and
(\ref{ccc}.2), we have $\frslI(A) = \frpsl_I(A)$ for a Lie torus if
$|I|=\infty$ or if $A$ does not have $|I|$-torsion. In particular,
the corollary says that for $\KK$ a field of characteristic $0$, a
Lie algebra $L$ over $\KK$ is a Lie-$(\rma_I, \ZZ^n)$-torus if and
only if $L$ is a $\ZZ^n$-covering of $\frslI(\KK_q)$ where $\KK_q$
is a quantum torus (see \ref{quatoex}). \es

\subsec[uce]{\bf The universal central extension of $\frslI(A)$ and
$\frpsl_I(A)$.} Let $\lan A,A \ran = (A \wedge A)\big / \scB$ where
$\scB = \Span_\KK\{ab \ot c + bc \ot a + ca \ot b : a,b,c\in A\}$.
It follows that $\lan A,A\ran$ is spanned by elements $\lan a,b\ran
= a \wedge b + \scB$, $a,b\in A$, and that $\lan a,b\ran \mapsto
[a,b]$ is a well-defined map $\lan A, A \ran  \to [A,A]$.  The
kernel of this map is $\HC_1(A)$, the first cyclic homology group of
$A$. \sm

We also need to recall the {\it Steinberg Lie algebra
$\frst_I(A)$\/}\index{S@$\frst$ (Steinberg Lie algebra)}\index{Lie
algebra!Steinberg}\index{Steinberg Lie algebra} associated to $A$.
It is the Lie $k$-algebra presented by generators $X_{ij}(a)$, where
$i,j\in I$, $i \ne j$, and $a\in A$, and relations (st1)--(st3): \sm

(st1) The map $a \mapsto X_{ij}(a)$ is $k$-linear, %\sm

(st2) $[X_{ij}(a), \, X_{jl}(b)] = X_{il}(ab)$ for distinct
$i,j,l\in I$, and %\sm

(st3) $[X_{ij}(a), \, X_{lm}(b)] = 0$ for $i\ne m$ and $j\ne
l$. \sm%\ms

\noindent Under our assumptions on $k$ and $I$, {\it the universal
central extension of $\frslI(A)$ is the map}
$$
   \fru \colon\, \frst_I(A) \to \frslI(A), \quad X_{ij}(a) \mapsto
aE_{ij}. \quad \hbox{\it Its kernel is $\HC_1(A)$,}  \eq{(1)}%uce0}
$$
cf. \ref{predivex}. The $\rma_I$-grading of $\frst_I(A)$ is
determined by $\frst_I(A)_\al = X_{ij}(A)$ for $0 \ne \al =
\eps_i-\eps_j$. Any invertible element $aE_{ij} \in \frslI(A)$ lifts
to an invertible element $X_{ij}(a) \in \frst_I(A)$. It then follows
that $\frst_I(A)$ is $\rma_I$-graded, as stated in
Th.~\ref{typeAcl}. Similarly, if $A$ is $\La$-graded, then
$\frst_I(A)$ is $(\scq(\rma_I) \oplus \La)$-graded with homogeneous
spaces $\frst_I(A)_\al^\la = X_{ij}(A^\la)$ for
$\al=\eps_i-\eps_j\ne 0$. This explains Th.~\ref{grtypeAcl} and
Cor.~\ref{grtypeAco}, at least for $\frst_I(A)$, but then the
general case easily follows. \sm

If the extension $\fru$ is split on the level of $k$-modules, e.g.
if $\frslI(A)$ is projective as $k$-module (for example $k$ is a
field), we can identify $\frslI(A)$ with a $k$-submodule of
$\frst_I(A)$ such that $$\frst_I(A) = \HC_1(A) \oplus \frslI(A)$$ as
$k$-modules and the Lie algebra product $[.,.]_{\frst}$ of
$\frst_I(A)$ is given in terms of the Lie algebra product $[.,.]$ of
$\frslI(A)$ by $[\HC_1(A) ,\, \frst_I(A)]_\frst = 0$ and
$$
    [x,y]_{\frst} = \tsum_{i,j} \,\lan x_{ij}, y_{ji}\ran \oplus
 [x,y] \eq{(2)}%uce1}
$$
for $x=(x_{ij})$ and $y=(y_{ij})\in \frslI(A)$. \sm

In case $n \cdot 1_k \in k\ti$, one can give a description of
$\uce(\frslI(A))$ in the spirit of \ref{assba}: We define a Lie
algebra on
$$ \uce\big( \frsln(A) \big) = \lan A, A \ran  \oplus ( \frsln(k) \ot_k A)
 \eq{(3)}%uce2}
$$
by modifying the product formulas in \ref{assba} as follows:
\begin{eqnarray*} %$$\eqalign{
  [x\ot a, \, y\ot b] &=& (\ts \frac{1}{n} \tr(xy)\, \lan a,b\ran )
            \oplus (\; \ldots \text{ as in \ref{assba} } \ldots )
     \\
  {[ \lan c_1,c_2 \ran  , \, \lan d_1, d_2\ran  ]}
     &=& \lan \, [c_1,c_2] , \, [d_1, d_2]\, \ran , \\
 {[\lan c_1, c_2 \ran , \, x\ot a]} &=& x \ot \big[ [c_1, c_2], \, a\big]
 \end{eqnarray*}
for $x,y\in \frsln(k)$ and $a,b,c_i, d_i \in A$. The obvious map
$\uce\big( \frslI(A) \big) \to \frslI(A)$ is then a universal
central extension. \sm

For example, let $A=\KK[t^{\pm 1}]$ be the Laurent polynomial ring
in one variable $t$, where $\KK$ is a field of characteristic $0$.
Since $\HC_1(A)$ is $1$-dimensional, the Lie algebra $K=
\frslI(\KK[t^{\pm 1}]) \oplus \KK c$ constructed in \ref{afrreal} is
a universal central extension of $\frslI(A)$.
 \es

\subsec[propa]{\bf The centroid of $\frslI(A)$.} It is convenient
here and in the description of derivations of $\frslI(A)$ to use the
following notation. Any endomorphism $f\in \End_k(A)$, which leaves
$[A,A]$ invariant, induces an endomorphism of $\frslI(A)$, denoted
$\frslI(f)$,  by applying $f$ to every entry of a matrix in
$\frslI(A)$:
$$
   \frslI(f) \, (x) = (f(x_{ij})) \quad \hbox{for } x=(x_{ij}) \in
                       \frslI(A). \eq{(1)}%propa1}
$$
The map $$\frslI(.) \colon\, \{ f\in \End_k (A) : f([A,A]) \subset
[A,A] \} \to \End_k (\frslI(A))$$ is a monomorphism of $k$-algebras.
If $f\in \End_k(A)$ also leaves the centre $\rmZ(A)$ of $A$
invariant, the map $\frslI(f)$ descends to an endomorphism of
$\frpsl_I(A)$, naturally denoted $\frpsl_I(f)$. \ms

Recall (\ref{centroid}) that $\rmZ(A) \to \cent(A)$, $z \mapsto
L_z$, is an isomorphism of $k$-algebras. For $z\in \rmZ(A) \cong
\cent_k(A)$ the endomorphism $\chi_z = \frslI(L_z)$ of $\frslI(A)$,
mapping $x=(x_{ij})$ to $\chi_z(x) = (zx_{ij})$, belongs to the
centroid of $\frslI(A)$. The map
$$
  \rm\rmZ(A) \to \cent\big( \frslI(A)\big), \quad z \mapsto \chi_z$$
is an isomorphism of $k$-algebras. Moreover, $z \mapsto
\frpsl_I(L_z)$ is also an isomorphism between $\rmZ(A)$ and
$\cent(\frpsl_I(A))$. If $A$ is $\La$-graded, then so is the centre
$\rmZ(A)$ and hence also $\cent( \frslI(A))$. \sm

For the description of $\frslI(A)$ as a multiloop algebra
(\cite{abfp}) it is important to know when is $\frslI(A)$ finitely
generated over its centroid. This is the case if and only if $I$ is
finite and $A$ is finitely generated over $\rmZ(A)$. For example, if
$A=k_q$, $q=(q_{ij})\in \Mat_n(k)$, is a quantum torus over a field
$k$, $A$ is a finitely generated module over its centre if and only
all $q_{ij}\in k$ are roots of unity, cf. \ref{quatoex}.
\es%\endsubsec

\subsec[assIBF]{\bf Invariant bilinear forms.} For a $\Xi$-graded
algebra $X$, $\Xi$ an abelian group, we denote by $\GIF(X,\Xi)$, the
$k$-module of $\Xi$-graded invariant symmetric bilinear forms on
$X$, see \ref{pergradef}. In this subsection we will determine
$\GIF(\frslI(A), \scq(\rma_I)\oplus \La)$ for a $\La$-graded algebra
$A$. For simpler notation we put $\Xi= \scq(\rma_I) \oplus \La$. \sm

To $\inpr_A \in \GIF(A,\La)$ we associate the bilinear form
$\inpr_{{\mathfrak s}{\mathfrak l}}$ defined on $\frslI(A)$ by
$$
 \ts \big( \, \tsum_{i,j} \, x_{ij} E_{ij}  \mid
    \tsum_{p,q}\, y_{pq}E_{pq}\big)_{{\mathfrak s}{\mathfrak l}}
    = \tsum_{i,j}\, (x_{ij} \mid y_{ji })_A. \eq{(1)}%assIBF0}
$$
It is an easy exercise to show that $\inpr_{{\mathfrak s}{\mathfrak
l}}$ is invariant. Since it is obviously also $\Xi$-graded, we have
$\inpr_{{\mathfrak s}{\mathfrak l}} \in \GIF(\frslI(A), \Xi)$.
Observe
$$
   (aE_{ij} \mid bE_{ji})_{{\mathfrak s}{\mathfrak l}}
          = (a\mid b)_A \eq{(2)}%assIBF1}
$$
for $a,b\in A$ and distinct $i,j\in I$. Conversely, this formula
assigns to every $\inpr_{{\mathfrak s}{\mathfrak l}} \in
\GIF(\frslI(A), \Xi)$ a graded bilinear form $\inpr_A$ on $A$. It is
immediate, using (\ref{assba}.2), that $\inpr_A$ is
invariant. Since %by \no{inforsga}
every invariant form on $\frslI(A)$ is uniquely determined by its
restriction to $AE_{ij} \times AE_{ji}$, $i\ne j$, we have an
isomorphism of $k$-modules$$
  \GIF(\frslI(A), \Xi)\; {\buildrel \cong \over
   \rightarrow} \; \GIF(A,\La),
 \quad \inpr_{{\mathfrak s}{\mathfrak l}} \mapsto \inpr_A.
     \eq{(3)}%assIBF1a}
$$
Moreover, by (\ref{pergradef}.2), $\GIF(A,\La) \cong ( A^0\big/
[A,A]^0)^*$. Combining these two isomorphisms, it follows that {\it
any $\La$-graded invariant bilinear form $b$ of $\frslI(A)$ can be
uniquely written in the form
$$  b(x,y) = \ph (\tr(xy) + [A,A]) \quad
 \hbox{for $\ph \in (A^0/[A,A]^0)^*$.} \eq{(4)}%assIBF2}
$$}
Since $\frslI(A)$ is perfect, $\GIF({{\mathfrak s}{\mathfrak
l}}_I(A)) \cong \GIF(\frpsl_I(A))$ by (\ref{pergradef}.1). Combining
(3) with this isomorphism implies
$$ \GIF(\frpsl_I(A)) \cong \GIF(A,\La)
     \cong  (A^0 / [A,A]^0)^*  . \eq{(5)}%assIBF4}
$$ \sm

We now turn to the description of nondegenerate invariant forms.
Since the centre of a perfect Lie algebra is contained in the
radical of any invariant form, $\rmZ(\frslI(A)) = \{0\}$ is a
necessary condition for the existence of such a form on $\frslI(A)$.
But rather than assuming $\rmZ(\frslI(A)) = \{0\}$ we will work with
the centreless Lie algebra $\frpsl_I(A)$ and the isomorphism (5),
written as $\inpr_\frpsl \mapsto \inpr_A$. It is easily seen that
$$
 \hbox{\it $\inpr_\frpsl$ is nondegenerate}\quad \Iff \quad
 \hbox{\it $\inpr_A$ is nondegenerate.} \eq{(6)}%assIBF5}
$$

{\bf Example.} Let $A$ be a $\La$-torus over a field $F$. Since
$[A,A]^0 = [A,A]\cap A^0= \{0\}$ by (\ref{grtypeAco}.1),
$$
           \GIF(\frpsl_I(A)) \cong (A^0)^* %\eq{assIBF6}
$$ {\it is a
$1$-dimensional space. Moreover, any non-zero graded invariant
bilinear form on $\frpsl_I(A)$ is nondegenerate}. Indeed, any
non-zero bilinear form $\inpr$ on $A$ is of type $(a_1 \mid
a_2)=\psi(a_1a_2)$ for a linear form $\psi\in A^*$ with $\Ker \psi =
\bigoplus_{0\ne \la} A^\la$. Since $\Rad \inpr$ is graded and $1 \in
z A$ for every non-zero homogeneous $z\in A$, it follows that $\Rad
\inpr =\{ z\in A : zA \subset \Ker \psi \} =\{0\}$.
 \es

\subsec[deri]{\bf Derivations.} We denote by
$\wi\Mat_I(A)$\index{M@$\wi\Mat$ (row- and column-finite matrices)}
the space of all row- and column-finite matrices with entries from
$A$. This is an associative $k$-algebra with respect to the usual
matrix multiplication. Hence it becomes a Lie algebra, denoted
$\wi{\frgl}_I(A)$,\index{G@$\wi\frgl$ (general Lie algebra of row-
and column-finite matrices)} with respect to the commutator. It is
immediate that $[\wi {\frgl}_I(A),\, \frslI(A)] \subset \frslI(A)$.
We obtain a Lie algebra homomorphism $
   \wi {\frgl}_I(A) \to \Der_k \frslI(A), \, %\quad
     \hat x \mapsto \ad \hat x |\frslI(A)
$ with kernel $\{ z E : z\in \rmZ(A)\}$ where $E$ is the identity
matrix of $\wi{\frgl}_I(A)$. Its image will be denoted
$$
 \wi \IDer \,\frslI(A) = \{ \ad \hat x | \frslI(A) : \hat x \in
  \hat \frgl_I(A)\} \cong \wi {\frgl}_I(A) \big/
            \{ z E: z \in \rmZ(A) \}. \index{I@$\wi\IDer$ (completed
inner derivations)}
$$
Of course, $\wi \IDer\, \frslI(A) = \IDer \frslI(A)$ if $I$ is
finite. Since $[\hat x, \rmZ(\frslI(A))] = 0$ for $\hat x\in {\frgl}_I(A)$, $\ad \hat x$ descends to a derivation of $\frpsl_I(A)$. %We denote by $\wi\IDer\, {\fr psl}_I(A)$ the derivations
%obtained in this way.
\sm

Another class of derivations of $\frslI(A)$ and $\frpsl_I(A)$ are
given by the maps $\frslI(d)$ and $\frpsl_I(d)$ for $d\in
\Der_k(A)$, cf. (\ref{propa}.1). The map $\frslI \colon\, \Der_k(A)
\to \Der_k \frslI(A)$, $d \mapsto \frslI(d)$ is clearly a
monomorphism of Lie algebras. We denote by $\frslI( \Der_kA)$ its
image. One checks that
$$
   [ \frslI(d),\, \ad \hat x] = \ad d(\hat x) \quad
     \hbox{for $\hat x \in \hat {\frgl}_I(A)$.}
$$
%The derivations $\frslI(d)$ descend to derivations $\frpsl_I(d)$
%of $\frpsl_I(A)$.
\endsubsec

\subsec[derit]{\Theorem}{\it {\rm (a)} The ungraded derivation
algebras are \begin{eqnarray*}%$$\eqalignno{
 \Der_k \big(\frslI(A)\big) &=& \wi \IDer \, \frslI(A) +
\frslI(\Der_kA) \quad \hbox{with} \\ %eq{derit1}
        \frslI( \IDer A) &=& \wi \IDer \, \frslI(A) \cap
         \frslI(\Der_kA) \cong \IDer A, \\
 \CDer \big( \frslI(A) \big) &=& \frslI( \CDer A) \cong \CDer A, \\
 {\rm SDer} \big( \frslI(A) \big) &=& \wi \IDer \, \frslI(A) +
        \frslI({\rm SDer} A), \end{eqnarray*}%}$$
where in the last formula the derivations in ${\rm SDer}$ are
skew-symmetric with respect to a bilinear form $\inpr_{{\mathfrak
s}{\mathfrak l}}$ and the corresponding bilinear form $\inpr_A$,
cf.\/ {\rm (\ref{assIBF}.2)} \sm

{\rm (b)} Suppose $A$ is $\La$-graded. Then the homogeneous
subspaces $(\Der \,\frslI(A))^\la_\al$ of derivations of degree $\al
\oplus \la$ are given by
\begin{eqnarray*}%$$\eqalign{
 \big( \Der\, \frslI(A) \big)^\la_\al
 &=& \ad \big( \frslI(A)_\al^\la)
    = \{ \ad x : x \in A^\la E_{ij} \}, \qquad \al= \eps_i-\eps_j\ne 0,
 \\
\big( \Der\, \frslI(A) \big)^\la_0 &=&
 \{ \ad \hat x: \hat x \in \wi {\frgl}_I(A)
         \hbox{ diagonal, all }  x_{ii}\in A^\la\}
        \oplus \frslI\big( (\Der\, A)^\la\big)
\end{eqnarray*}%$$
In particular,
$$
   \grSCDer\, (\frslI(A)) = \frslI\big( \grSCDer( A)\big)
       \cong \grSCDer A.
$$

{\rm (c)} All statements in {\rm (a)} and {\rm (b)} hold, mutatis
mutandis, for $\frpsl_I(A)$. }\es

\subsec[asstor]{\bf The standard toral subalgebra.} In this
subsection we suppose that $k=\KK$ is a field of characteristic $0$.
We put
$$
   \frh = \Span_\KK \{ E_{ii} - E_{jj} : i,j\in I\}
$$
which is the space of diagonal trace-$0$-matrices over $\KK$. We
write an element $h\in \frh$ as $h={\rm diag}(h_i)_{i\in I}$ where
$h_i\in \KK$ is the entry of the matrix $h$ at the position $(ii)$.
We embed the root system $\rma_I$ into $\frh^*$ by
$$
    (\eps_i - \eps_j)(h) = h_i-h_j \quad \hbox{for } h = {\rm diag}(h_i).
$$ Then $\frh$ is a toral subalgebra of $\frslI(A)$
with set of roots equal to $\rma_I$.  \sm

Suppose $A= \bigoplus_{\la \in \La}A^\la$ is $\La$-graded and let
$\psi \in (A^0)^*$ with $[A\, A]^0 \subset \Ker \psi$. Hence $\psi$
gives rise to a $\La$-graded invariant symmetric bilinear form
$\inpr_{A, \, \psi}$ on $A$, given by $(a_1|a_2)_{A, \, \psi} =
\psi(a_1a_2)$. By (\ref{assIBF}.1) we can extend $\inpr_{A,\, \psi}$
to an invariant $(\scQ(\rma_I) \oplus \La)$-graded symmetric
bilinear form $\inpr_{{{\mathfrak s}{\mathfrak l}},\, \psi}$ on
$\frslI(A)$. Note that by (\ref{assIBF}.4) for $h={\rm diag}(h_i)$
and $h'={\rm diag}(h_i')$ we get
$$
   ( h | h')_{{{\mathfrak s}{\mathfrak l}},\, \psi} = \tsum_{i\in I} \psi(h_i h'_i) =
      \big(\tsum_{i\in I}\, h_i h'_i\big)\,  \psi(1).
$$
This implies
$$
   \inpr_{{{\mathfrak s}{\mathfrak l}},\, \psi} \mid \frh \times \frh \hbox{\it\  is
nondegenerate if and only if } \psi(1) \ne 0. %\eq{asstor1}
$$
Indeed, $\psi(1) \ne 0 $ is obviously necessary for nondegeneracy.
Conversely, if $(h | \frh)_{{{\mathfrak s}{\mathfrak l}},\, \psi} =
0$ for $h={\rm diag}(h_i)$, then $h_i=h_j$ for all $i,j\in I$ whence
$h=0$ in case $I$ is infinite, but this also holds in case of a
finite $I$, say $|I|=n$, since then  $h=h_1 E_n$ has trace $nh_1 =
0$. \sm

As we have seen in (\ref{assIBF}.6), the form $\inpr_{{{\mathfrak
s}{\mathfrak l}},\, \psi}$ descends to a $(\scq(\rma_I)\oplus
\La)$-graded invariant symmetric bilinear form $\inpr_{\frpsl,\,
\psi}$ on $\frpsl_I(A)$. By (\ref{assIBF}.6),
$$\amen{0.70}{$\frpsl_I(A)$ is an invariant $(\rma_I,\La)$-graded
Lie algebra
%as defined in \no{invarex}
if and only if $A$ has a $\La$-graded invariant nondegenerate
symmetric bilinear form, say $\inpr_{A,\, \psi}$, and $\psi(1) \ne
0$.}$$

\noindent In particular, if $A$ is a $\La$-torus the Lie torus
$\frpsl_I(A)$ is invariant with respect to any form $\inpr_{\frpsl,\, \psi}$ with $\psi \ne 0$. %\endsubsec%
\es

%{\petit
\subsec[exassecnotes]{\bf Notes.} The definition of $\frslI(A)$ is
standard. Some authors call it the {\it elementary Lie algebra\/}
and denote it by ${\mathfrak e}_I(A)$. One can of course also define
$\frslI(A)$ for $|I|=2$ (recall our convention $|I|\ge 3$, but then
the claims in this section, e.g. (\ref{assba}.1), are no longer
true. Some results on ${{\mathfrak s}{\mathfrak l}}_n(A)$, $n\in
\NN$, are proven in \cite{bgk}, e.g. the formula in \ref{assba} for
the Lie algebra product of $\frsln(A)$ can be found in
\cite[(2.27)--(2.29)]{bgk}. \sm

In the generality stated, Th.~\ref{typeAcl} is proven in
\cite[(3.5)]{n:3g}. For $k$ a field of characteristic $0$ and $I$
finite, i.e., $\rma_I$ finite, it was proven earlier by Berman-Moody
\cite[0.7]{bm}. Their proof shows at the same time that $\frst_n(A)$
is a universal central extension of any $\rmfa_{n-1}$-graded Lie
algebra. \sm

As stated, Th.~\ref{grtypeAcl} is an immediate consequence of
\cite[2.11 and 3.4]{gn2} and \cite[4.7]{n:s}. Cor.~\ref{grtypeAco}
is proven in \cite[Prop.~2.13]{y2} for $I$ finite and $k$ a field of
characteristic $0$. This paper also contains a construction of
crossed products $B* \ZZ^n$ (\cite[\S  3]{y2}). The description of
Lie tori of type $({\rm A}_l,\ZZ^n)$ is given in \cite{bgk}. The
formula (\ref{grtypeAco}.1) is proven in \cite[Prop.~3.3]{NY}. The
case $\La=\ZZ^n$ is already contained in \cite[Prop.~2.44]{bgk}. \sm

That (\ref{uce}.1) is a universal central extension is proven in
\cite{Bloch} for $A$ commutative and $|I|\ge 5$, for $A$ arbitrary
and again $|I|\ge 5$ in \cite{KasLod} and for $3\le |I|\le 4$ in
\cite[2.63]{gao:unitary} (recall our assumption that $\frac{1}{6}
\in k$). The universal central $2$-cocycle (\ref{uce}.2) is also
given in these papers. The results in \ref{uce} are no longer true
if $2$ or $3$ are not invertible in $k$, see \cite{GaoShang} for the
case of an algebra $A$ that is free as $k$-module.  If $\frslI(A)$
is a free $k$-module or a finitely generated and projective
$k$-module, the kernel of the universal central extension is the
second homology with trivial coefficients, i.e., $H_2\big(
{{\mathfrak s}{\mathfrak l}}_I(A)\big) = \HC_1(A)$. That the Lie
algebra (\ref{uce}.3) is the universal central extension of
$\frsln(A)$ is shown in \cite[page 359--360]{bgk}. A lot is known
about the cyclic homology groups $\HC_1(A)$, see \cite{loday}. For
example, if $A$ is commutative, then $\HC_1(A) \cong
\Omega_{A|k}^1\big/ dA$ \cite[2.1.14]{loday}. For a complex quantum
torus $\HC_1(\CC_q)$ is described in \cite[Prop.~3.19]{bgk} (the
result is true for any field). In particular, it follows that $\dim
\HC_1(\KK[t^{\pm 1} ]) =1$. The formula (\ref{assIBF}.4) is proven
in \cite[Lemma~2.8]{bgk}. \sm

The description of $\Der \frslI(A)$ in \ref{derit}(a) is proven for
finite $I$ and $k$ a field (but $A$ arbitrary non-associative) in
\cite[Th.~4.8 and Cor.~4.9]{BenOsb:der}. The proof can easily be
adapted to yield the result as stated. For $k=A$ a field of
characteristic $0$ one has $\Der_k (k) = \{0\}$, $\frslI(k) =
\frpsl_I(k)$ and hence $\Der_k \big(\frslI(k)\big) = \wi \IDer \,
\frslI(k)$, which is proven in \cite[Th.~I.3]{neeb:der} with a
different method. The derivation algebra of $\CC_q$ is determined in
\cite[Lemma~2.48]{bgk}. \ms

%old A2-section \subsec[a2rem]{Remark on the ${\rm A}_2$-case.}
The astute reader will have noticed that we have excluded the ${\rm
A}_2$-case in all statements referring to the description of
$(\rma_I, \La)$-graded Lie algebras (\ref{typeAcl}, \ref{grtypeAcl}
and \ref{grtypeAco}). The reason for this is that not all $(\rmfa_2,
\La)$-graded Lie algebra are $\La$-coverings of $\frpsl_3(A)$ for
$A$ a $\La$-graded associative algebra. Rather, one has to allow
alternative algebras $A$. Various models of $\rmfa_2$-graded Lie
algebras have been given, see e.g. \cite[3.3]{bm}, \cite[\S
2]{bgkn}, \cite[3.3]{n:3g} or \cite[\S 9]{AF:isotopy}. Once a
replacement of $\frpsl_3(A)$ has been defined, the Theorems
\ref{typeAcl}, \ref{grtypeAcl} and \ref{grtypeAco} hold, mutatis
mutandis, by replacing the associative algebra $A$ by an alternative
algebra $A$. For example, for $\La= \ZZ^n$ and Lie algebras over
fields of characteristic $0$ this has been worked out in \cite[\S
6]{y:alt}, while the classification of the corresponding coordinate
algebras, i.e., division-$\ZZ^n$-graded alternative algebras, is
given in \cite[Th.~5.7]{y:alt} and is valid in characteristic $\ne
2$. For arbitrary $\La$, special types of alternative division
$\La$-graded algebras, so-called quasi-algebras, are described in
\cite{AlElPe}. An alternative nonassociative $\ZZ^n$-torus over a
field $F$ of characteristic $\ne 2$ is graded-isomorphic to an
octonion torus \cite[Th.~1.25]{bgkn}, \cite[Cor.~5.13]{y:alt}, see
\cite[Example~9.2]{AF:isotopy} for a presentation.
\es%\endsubsec

\vfe

\def\cprime{$'$} \def\cprime{$'$}
\def\bysame{\leavevmode\hbox to3em{\hrulefill}\thinspace}

\printindex

\begin{thebibliography}{99999.}

\bibitem[AEP]{AlElPe}
H.~Albuquerque, A.~Elduque, and J.~M. P{\'e}rez-Izquierdo, {\it
Alternative
  quasialgebras}, Bull. Austral. Math. Soc. {\bf 63} (2001), no.~2, 257--268.

\bibitem[AABGP]{aabgp}
B.~Allison, S.~Azam, S.~Berman, Y.~Gao, and A.~Pianzola, {\it
Extended affine
  {L}ie algebras and their root systems}, Mem. Amer. Math. Soc. {\bf 126}
  (1997), no.~603, x+122.

\bibitem[ABG1]{abg}
B.~Allison, G.~Benkart, and Y.~Gao, {\it Central extensions of {L}ie
algebras
  graded by finite root systems}, Math. Ann. {\bf 316} (2000), no.~3, 499--527.

\bibitem[ABG2]{abg2}
\bysame, {\it Lie algebras graded by the root systems {${\rm BC}_r,\
r\ge2$}},
  Mem. Amer. Math. Soc. {\bf 158} (2002), no.~751, x+158.

\bibitem[ABFP1]{abfp}
B.~Allison, S.~Berman, J.~Faulkner, and A.~Pianzola, {
 {\it Realization of graded-simple algebras as loop
algebras},
 arXiv:math.RA/0511723, to appear in Forum Mathematicum.

\bibitem[ABFP2]{abfp2}
\bysame,  \it Muliloop realization
  of extended affine {L}ie algebras and {L}ie tori},
arXiv:0709.0975 [math.RA].

\bibitem[ABGP]{abgp}
B.~Allison, S.~Berman, Y.~Gao, and A.~Pianzola, {\it A
characterization of
  affine {K}ac-{M}oody {L}ie algebras}, Comm. Math. Phys. {\bf 185} (1997),
  no.~3, 671--688.

\bibitem[ABP1]{abp1}
B.~Allison, S.~Berman, and A.~Pianzola, {\it Covering algebras. {I}.
{E}xtended
  affine {L}ie algebras}, J. Algebra {\bf 250} (2002), no.~2, 485--516.

\bibitem[ABP2]{abp2}
\bysame, {\it Covering algebras. {II}. {I}somorphism of loop
algebras}, J.
  Reine Angew. Math. {\bf 571} (2004), 39--71.

\bibitem[ABP3]{abp2.5}
\bysame, {\it Iterated loop algebras}, Pacific J. Math. {\bf 227}
(2006),
  no.~1, 1--41.

\bibitem[AF]{AF:isotopy}
B.~Allison and J.~Faulkner, {\it Isotopy for extended affine {L}ie
algebras and {L}ie tori}, this volume and arXiv:0709.1181 [math.RA].

\bibitem[AG]{AG2}
B.~N. Allison and Y.~Gao, {\it The root system and the core of an
extended
  affine {L}ie algebra}, Selecta Math. (N.S.) {\bf 7} (2001), no.~2, 149--212.

\bibitem[Az1]{Az1}
S.~Azam, {\it Nonreduced extended affine root systems of nullity
{$3$}}, Comm.
  Algebra {\bf 25} (1997), no.~11, 3617--3654.

\bibitem[Az2]{Az5}
\bysame, {\it Construction of extended affine {L}ie algebras by the
twisting
  process}, Comm. Algebra {\bf 28} (2000), no.~6, 2753--2781.

\bibitem[Az3]{az:ext}
\bysame, {\it Extended affine root systems}, J. Lie Theory {\bf 12}
(2002), 515--527.

\bibitem[Az4]{az:grla}
\bysame, {\it Generalized reductive {L}ie algebras: connections with
extended
  affine {L}ie algebras and {L}ie tori}, Canad. J. Math. {\bf 58} (2006),
  no.~2, 225--248.

\bibitem[Az5]{Azam:multi}
\bysame, {\it Derivations of multi-loop algebras}, Forum Math. {\bf
19} (2007),
  no.~6, 1029--1045.

%\bibitem[Az6]{Azam:tensor}
%\bysame, {\it Derivations of tensor product of algebras},
%arXiv:math.QA/0504368.

\bibitem[ABY]{ABY}
S.~Azam, S.~Berman, and M.~Yousofzadeh, {\it Fixed point subalgebras
of
  extended affine {L}ie algebras}, J. Algebra {\bf 287} (2005), no.~2,
  351--380.

\bibitem[AKY]{AKY}
S.~Azam, V.~Khalili, and M.~Yousofzadeh, {\it Extended affine root
systems of
  type {BC}}, J. Lie Theory {\bf 15} (2005), no.~1, 145--181.

\bibitem[Be]{be:derinv}
G.~Benkart, {\it Derivations and invariant forms of {L}ie algebras
graded by
  finite root systems}, Canad. J. Math. {\bf 50} (1998), no.~2, 225--241.

\bibitem[BM]{bemo}
G.~Benkart and R.~Moody, {\it Derivations, central extensions, and
affine {L}ie
  algebras}, Algebras Groups Geom. {\bf 3} (1986), no.~4, 456--492.

\bibitem[BN]{BN}
G.~Benkart and E.~Neher, {\it The centroid of extended affine and
root graded
  {L}ie algebras}, J. Pure Appl. Algebra {\bf 205} (2006), no.~1, 117--145.

\bibitem[BO]{BenOsb:der}
G.~Benkart and J.~M. Osborn, {\it Derivations and automorphisms of
  nonassociative matrix algebras}, Trans. Amer. Math. Soc. {\bf 263} (1981),
  no.~2, 411--430.

\bibitem[BS]{BeSm}
G.~Benkart and O.~Smirnov, {\it Lie algebras graded by the root
system {${\rm
  BC}_1$}}, J. Lie Theory {\bf 13} (2003), no.~1, 91--132.

\bibitem[BY]{BeYo}
G.~Benkart and Y.~Yoshii, {\it Lie {$G$}-tori of symplectic type},
Q. J. Math.
  {\bf 57} (2006), no.~4, 425--448.

\bibitem[BZ]{beze}
G.~Benkart and E.~Zelmanov, {\it Lie algebras graded by finite root
systems and
  intersection matrix algebras}, Invent. Math. {\bf 126} (1996), 1--45.

\bibitem[BGK]{bgk}
S.~Berman, Y.~Gao, and Y.~Krylyuk, {\it Quantum tori and the
structure of
  elliptic quasi-simple {L}ie algebras}, J. Funct. Anal. {\bf 135} (1996),
  339--389.

\bibitem[BGKN]{bgkn}
S.~Berman, Y.~Gao, Y.~Krylyuk, and E.~Neher, {\it The alternative
torus and the
  structure of elliptic quasi-simple {L}ie algebras of type {A}$_2$}, Trans.
  Amer. Math. Soc. {\bf 347} (1995), 4315--4363.

\bibitem[BeMo]{bm}
S.~Berman and R.~Moody, {\it Lie algebras graded by finite root
systems and the
  intersection matrix algebras of {S}lodowy}, Invent. Math. {\bf 108} (1992),
  323--347.

\bibitem[Bl1]{Bloch}
S.~Bloch, {\it The dilogarithm and extensions of {L}ie algebras},
Algebraic
  $K$-theory, Evanston 1980 (Proc. Conf., Northwestern Univ., Evanston, Ill.,
  1980), Lecture Notes in Math., vol. 854, Springer, Berlin, 1981, pp.~1--23.

\bibitem[Bl2]{block}
R.~Block, {\it Determination of the differentiably simple rings with
a minimal
  ideal.}, Ann. of Math. (2) {\bf 90} (1969), 433--459.

\bibitem[Bo1]{bou:lie78}
N.~Bourbaki, {\it Groupes et alg{\`e}bres de {L}ie, chapitres 7--8},
Hermann,
  Paris, 1975.

\bibitem[Bo2]{brac}
\bysame, {\it Groupes et alg{\`e}bres de {L}ie, chapitres 4--6},
Masson, Paris,
  1981.

\bibitem[EMO]{emo}
T.~S. Erickson, W.~S. Martindale, and J.~M. Osborn, {\it Prime
nonassociative
  algebras}, Pacific J. Math. {\bf 60} (1975), no.~1, 49--63.

\bibitem[ERM]{Rao-Moody}
S.~Eswara~Rao and R.~V. Moody, {\it Vertex representations for
{$n$}-toroidal
  {L}ie algebras and a generalization of the {V}irasoro algebra}, Comm. Math.
  Phys. {\bf 159} (1994), no.~2, 239--264.

\bibitem[Fa]{Fa1}
R.~Farnsteiner, {\it Derivations and central extensions of finitely
generated
  graded {L}ie algebras}, J. Algebra {\bf 118} (1988), 33--45.

\bibitem[Ga1]{gao:unitary}
Y.~Gao, {\it Steinberg unitary {L}ie algebras and skew-dihedral
homology}, J.
  Algebra {\bf 179} (1996), no.~1, 261--304.

\bibitem[Ga2]{gao:deg}
\bysame, {\it The degeneracy of extended affine {L}ie algebras},
Manuscripta
  Math. {\bf 97} (1998), no.~2, 233--249.

\bibitem[GS]{GaoShang}
Y.~Gao and S.~Shang, {\it Universal coverings of {S}teinberg {L}ie
algebras of
  small characteristic}, J. Algebra {\bf 311} (2007), no.~1, 216--230.

\bibitem[GaN]{gn2}
E.~Garc{\'{\i}}a and E.~Neher, {\it Tits-{K}antor-{K}oecher
superalgebras of
  {J}ordan superpairs covered by grids}, Comm. Algebra {\bf 31} (2003), no.~7,
  3335--3375.

\bibitem[Ga]{gar}
H.~Garland, {\it The arithmetic theory of loop groups}, Inst. Hautes
\'Etudes
  Sci. Publ. Math. (1980), no.~52, 5--136.

\bibitem[GrN]{gross-nebe}
B.~H. Gross and G.~Nebe, {\it Globally maximal arithmetic groups},
J. Algebra
  {\bf 272} (2004), no.~2, 625--642.

\bibitem[Ha]{hartwig}
J.~T. Hartwig, {\it Locally finite simple weight modules over
twisted
  generalized {W}eyl algebras}, J. Algebra {\bf 303} (2006), no.~1, 42--76.

\bibitem[He]{hee}
J.-Y. H{\'e}e, {\it Syst\`emes de racines sur un anneau commutatif
totalement
  ordonn\'e}, Geom. Dedicata {\bf 37} (1991), no.~1, 65--102.

\bibitem[HT]{HKT}
R.~H{\o}egh-Krohn and B.~Torr{\'e}sani, {\it Classification and
construction of
  quasisimple {L}ie algebras}, J. Funct. Anal. {\bf 89} (1990), no.~1,
  106--136.

\bibitem[Ho]{HofG}
G.~W. Hofmann, {\it Symmetric systems and their applications to root
systems
  extended by abelian groups}, arXiv:0712.0104 [math.GR].

\bibitem[Hu]{hum:cox}
J.~Humphreys, {\it Reflection groups and {C}oxeter groups},
Cambridge Studies
  in Advanced Mathematics, vol.~29, Cambridge University Press, Cambridge,
  1990.

\bibitem[Ja]{jake:lie}
N.~Jacobson, {\it Lie algebras}, Interscience Tracts in Pure and
Applied
  Mathematics, No. 10, Interscience Publishers (a division of John Wiley \&
  Sons), New York-London, 1962.

\bibitem[Kac1]{kac:sup1}
V.~G. Kac, {\it Lie superalgebras}, Advances in Math. {\bf 26}
(1977), no.~1,
  8--96.

\bibitem[Kac2]{kac}
\bysame, {\it Infinite dimensional {L}ie algebras}, third ed.,
Cambridge
  University Press, 1990.

\bibitem[Kac3]{kac:locality}
\bysame, {\it The idea of locality}, Physical applications and
mathematical
  aspects of geometry, groups and algebras (Singapure), World Sci., 1997,
  pp.~16--32.

\bibitem[Kan1]{Kan1}
I.~L. Kantor, {\it Classification of irreducible transitive
differential
  groups}, Dokl. Akad. Nauk SSSR {\bf 158} (1964), 1271--1274.

\bibitem[Kan2]{Kan2}
\bysame, {\it Non-linear groups of transformations defined by
general norms of
  {J}ordan algebras}, Dokl. Akad. Nauk SSSR {\bf 172} (1967), 779--782.

\bibitem[Kan3]{Kan3}
\bysame, {\it Certain generalizations of {J}ordan algebras}, Trudy
Sem. Vektor.
  Tenzor. Anal. {\bf 16} (1972), 407--499.

\bibitem[Kas]{Kas}
C.~Kassel, {\it K\"ahler differentials and coverings of complex
simple {L}ie
  algebras extended over a commutative algebra}, Proceedings of the Luminy
  conference on algebraic $K$-theory (Luminy, 1983), vol.~34, 1984,
  pp.~265--275.

\bibitem[KL]{KasLod}
C.~Kassel and J.-L. Loday, {\it Extensions centrales d'alg\`ebres de
{L}ie},
  Ann. Inst. Fourier (Grenoble) {\bf 32} (1982), no.~4, 119--142 (1983).

\bibitem[Ko1]{Koe1}
M.~Koecher, {\it Imbedding of {J}ordan algebras into {L}ie algebras
I}, Amer.
  J. Math. {\bf 89} (1967), 787--816.

\bibitem[Ko2]{Koe2}
\bysame, {\it Imbedding of {J}ordan algebras into {L}ie algebras
II}, Amer. J.
  Math. {\bf 90} (1968), 476--510.

\bibitem[Lod]{loday}
J.-L. Loday, {\it Cyclic homology}, Grundlehren, vol. 301,
Springer-Verlag,
  Berlin Heidelberg, 1992.

\bibitem[Loo]{l:sp}
O.~Loos, {\it {S}piegelungsr\"aume und homogene symmetrische
{R}\"aume}, Math.
  Z. {\bf 99} (1967), 141--170.

\bibitem[LN1]{lfrs}
O.~Loos and E.~Neher, {\it Locally finite root systems}, Mem. Amer.
Math. Soc. {\bf 171}
  (2004), no.~811, x+214.

\bibitem[LN2]{prs}
\bysame,{\it Reflection systems and partial root systems}, Jordan
Theory Preprint Archives, paper \#182.

\bibitem[Ma]{manin}
Y.~I. Manin, {\it Topics in noncommutative geometry}, M. B. Porter
Lectures,
  Princeton University Press, Princeton, NJ, 1991.

\bibitem[Mc]{mc:taste}
K.~McCrimmon, {\it A taste of {J}ordan algebras}, Universitext,
  Springer-Verlag, New York, 2004.

\bibitem[MZ]{mczel}
K.~McCrimmon and E.~Zel{\cprime}manov, {\it The structure of
strongly prime
  quadratic {J}ordan algebras}, Adv. in Math. {\bf 69} (1988), no.~2, 133--222.

\bibitem[MP]{mp}
R.~V. Moody and A.~Pianzola, {\it Lie algebras with triangular
decompositions},
  Can. Math. Soc. series of monographs and advanced texts, John Wiley, 1995.

\bibitem[MY]{morita-yoshii}
J.~Morita and Y.~Yoshii, {\it Locally extended affine {L}ie
algebras}, J.
  Algebra {\bf 301} (2006), no.~1, 59--81.

\bibitem[Nee1]{neeb:split}
K.-H.~Neeb, {\it Integrable roots in split graded {L}ie algebras},
J. Algebra {\bf
  225} (2000), no.~2, 534--580.

\bibitem[Nee2]{neeb:der}
\bysame, {\it Derivations of locally simple {L}ie algebras}, J. Lie
Theory {\bf
  15} (2005), no.~2, 589--594.

\bibitem[Nee3]{neeb:quantumtori}
\bysame, {\it On the classification of rational quantum tori and
structure of their automorphism groups}, arXiv:RA/0511263.

\bibitem[NS]{neeb-stu}
K.-H.~Neeb and N.~Stumme, {\it The classification of locally finite
split
  simple {L}ie algebras}, J. Reine Angew. Math. {\bf 533} (2001), 25--53.


\bibitem[Neh1]{n:cr}
E.~Neher, {\it Syst{\`e}mes de racines 3-gradu{\'e}s}, C. R. Acad.
Sci. Paris
  S{\'e}r. I {\bf 310} (1990), 687--690.

\bibitem[Neh2]{n:gen}
\bysame, {\it Generators and relations for {$3$}-graded {L}ie
algebras}, J.
  Algebra {\bf 155} (1993), no.~1, 1--35.

\bibitem[Neh3]{n:3g}
\bysame, {\it Lie algebras graded by 3-graded root systems and
{J}ordan pairs
  covered by a grid}, Amer. J. Math {\bf 118} (1996), 439--491.

\bibitem[Neh4]{n:uce}
\bysame, {\it An introduction to universal central extensions of
{L}ie
  superalgebras}, Groups, rings, Lie and Hopf algebras (St. John's, NF, 2001),
  Math. Appl., vol. 555, Kluwer Acad. Publ., Dordrecht, 2003, pp.~141--166.

\bibitem[Neh5]{n:s}
\bysame, {\it Quadratic {J}ordan superpairs covered by grids}, J.
Algebra {\bf
  269} (2003), no.~1, 28--73.

\bibitem[Neh6]{n:tori}
\bysame, {\it Lie tori}, C. R. Math. Acad. Sci. Soc. R. Can. {\bf
26} (2004),
  no.~3, 84--89.

\bibitem[Neh7]{n:eala}
\bysame, {\it Extended affine {L}ie algebras}, C. R. Math. Acad.
Sci. Soc. R.
  Can. {\bf 26} (2004), no.~3, 90--96.


\bibitem[Neh8]{n:lect}
\bysame, {\it Lectures on extended affine {L}ie algebras, {L}ie
tori, and
  beyond}, in preparation.

\bibitem[NT]{NeTo}
E.~Neher and M.~Toc\'on, {\it Graded-simple {L}ie algebras of type
$\rmb_2$ and
  graded-simple {J}ordan pairs covered by a triangle}, in preparation.

\bibitem[NY]{NY}
E.~Neher and Y.~Yoshii, {\it Derivations and invariant forms of
{J}ordan and
  alternative tori}, Trans. Amer. Math. Soc. {\bf 355} (2003), no.~3,
  1079--1108.

\bibitem[OP]{OP}
J.~M. Osborn and D.~S. Passman, {\it Derivations of skew polynomial
rings}, J.
  Algebra {\bf 176} (1995), no.~2, 417--448.

\bibitem[Pa]{pass}
D.~S. Passman, {\it Infinite crossed products}, Pure and Applied
Mathematics,
  vol. 135, Academic Press Inc., Boston, MA, 1989.

\bibitem[Sa]{sa}
K.~Saito, {\it Extended affine root systems. {I}. {C}oxeter
transformations},
  Publ. Res. Inst. Math. Sci. {\bf 21} (1985), no.~1, 75--179.

\bibitem[Se]{Se:rat}
G.~B. Seligman, {\it Rational methods in {L}ie algebras}, Marcel
Dekker Inc.,
  New York, 1976, Lecture Notes in Pure and Applied Mathematics, Vol. 17.

\bibitem[SS]{springer-steinberg}
T.~A. Springer and R.~Steinberg, {\it Conjugacy classes}, Seminar on
Algebraic
  Groups and Related Finite Groups (The Institute for Advanced Study,
  Princeton, N.J., 1968/69), Lecture Notes in Mathematics, Vol. 131, Springer,
  Berlin, 1970, pp.~167--266.

\bibitem[Ste]{stein:gen}
R.~Steinberg, {\it G\'en\'erateurs, relations et rev\^etements de
groupes
  alg\'ebriques}, Colloq. Th\'eorie des Groupes Alg\'ebriques (Bruxelles,
  1962), Librairie Universitaire, Louvain, 1962, pp.~113--127.

\bibitem[Stu]{stu:struct}
N.~Stumme, {\it The structure of locally finite split {L}ie
algebras}, J.
  Algebra {\bf 220} (1999), no.~2, 664--693.

\bibitem[Ti]{tits:JA}
J.~Tits, {\it Une classe d'alg{\`e}bres de {L}ie en relation avec
les
  alg{\`e}bres de {J}ordan}, Indag. Math. {\bf 24} (1962), 530--535.

\bibitem[vdK]{vdK}
W.~L.~J. van~der Kallen, {\it Infinitesimally central extensions of
{C}hevalley
  groups}, Springer-Verlag, Berlin, 1973, Lecture Notes in Mathematics, Vol.
  356.

\bibitem[We]{wei}
C.~Weibel, {\it An introduction to homological algebra}, Cambridge
studies in
  advanced mathematics, vol.~38, Cambridge University Press, 1994.

\bibitem[Wi]{wilson}
R.~L. Wilson, {\it Euclidean {L}ie algebras are universal central
extensions},
  Lie algebras and related topics (New Brunswick, N.J., 1981), Lecture Notes in
  Math., vol. 933, Springer, Berlin, 1982, pp.~210--213.

\bibitem[Yo1]{y2}
Y.~Yoshii, {\it Root-graded {L}ie algebras with compatible grading},
Comm.
  Algebra {\bf 29} (2001), no.~8, 3365--3391.

\bibitem[Yo2]{y:alt}
\bysame, {\it Classification of division {$\ZZ^n$}-graded
alternative
  algebras}, J. Algebra {\bf 256} (2002), no.~1, 28--50.

\bibitem[Yo3]{y:inv}
\bysame, {\it Classification of quantum tori with involution},
Canad. Math.
  Bull. {\bf 45} (2002), no.~4, 711--731.

\bibitem[Yo4]{y:ext}
\bysame, {\it Root systems extended by an abelian group and their
{L}ie
  algebras}, J. Lie Theory {\bf 14} (2004), no.~2, 371--394.

\bibitem[Yo5]{y:lie}
\bysame, {\it Lie tori---a simple characterization of extended
affine {L}ie
  algebras}, Publ. Res. Inst. Math. Sci. {\bf 42} (2006), no.~3, 739--762.

\bibitem[Yo6]{y:LARS}
\bysame, {\it Locally extended affine root systems},
preprint, April 2008.

\bibitem[Z]{zel:primeII}
E.~I. Zel{\cprime}manov, {\it Prime {J}ordan algebras. {II}},
Sibirsk. Mat. Zh.
  {\bf 24} (1983), no.~1, 89--104, 192.

\end{thebibliography}
\end{document}